\documentclass[12pt,sort&compress]{elsarticle}

\usepackage{setspace}
\usepackage[top=1in,bottom=1in,left=1in,right=1in]{geometry}
\usepackage{amssymb}
\usepackage{amsmath}
\usepackage{amsfonts}
\usepackage{amsthm}
\usepackage{newtxtext, newtxmath}
\usepackage{color}
\usepackage[skip=4pt]{caption}
\usepackage{subcaption}
\usepackage{float}
\usepackage{url}
\usepackage{graphicx}

\usepackage[pdfborder={0 0 0},colorlinks,allcolors=blue]{hyperref}

\theoremstyle{definition}
\newtheorem{remark}{Remark}

\graphicspath{{figures_arxiv/}}

\newcommand{\triplebar}{\vert\kern-0.2ex\vert\kern-0.2ex\vert}

\journal{Mathematical Models and Methods in Applied Sciences}

\begin{document}
\begin{frontmatter}

\title{A review of weakly enforced Dirichlet boundary conditions in computational flow analysis}

\author[inst1]{Ming-Chen~Hsu\corref{cor}}
\ead{jmchsu@iastate.edu}
\author[inst1]{Monu~Jaiswal}
\author[inst2,inst3,inst4]{Yuri~Bazilevs\corref{cor}}
\ead{yuri.bazilevs@vanderbilt.edu}

\cortext[cor]{Corresponding authors}

\address[inst1]{Department of Mechanical Engineering, Iowa State University, Ames, Iowa 50011, USA}
\address[inst2]{Department of Civil and Environmental Engineering, Vanderbilt University, Nashville, Tennessee 37235, USA}
\address[inst3]{Department of Mechanical Engineering, Vanderbilt University, Nashville, Tennessee 37240, USA}
\address[inst4]{Faculty of Science and Engineering, Graduate School of Creative Science and Engineering, Waseda University, Shinjuku City, Tokyo 169-8050, Japan}

\begin{abstract}
Strongly enforced Dirichlet boundary conditions require highly refined near-wall meshes to resolve steep velocity and thermal gradients. This introduces high computational costs, especially for practical flow simulations. Weakly enforced boundary conditions alleviate this burden by acting as a variationally consistent near-wall model. By allowing a controlled slip at the solid wall, weak enforcement recovers accurate flow quantities on coarse boundary-layer meshes across both incompressible and compressible regimes. Furthermore, weak boundary conditions serve as the fundamental enabling technology for immersogeometric analysis. Because the weak operator evaluates boundary integrals independently of the background mesh, high-fidelity flow analysis can be performed directly on complex geometries without fitting a mesh to the surface. This capability has facilitated direct geometry-to-analysis workflows for boundary-representation CAD models, raw point clouds, photogrammetric reconstructions, and segmented medical images. This review examines the unified mathematical development of the weak boundary condition framework from scalar advection--diffusion equations to the full Navier--Stokes equations, illustrating its versatility and robustness across incompressible and compressible flows, whether in traditional boundary-fitted, sliding-interface, or advanced immersogeometric applications.
\end{abstract}


\end{frontmatter}

\tableofcontents

\section{Introduction}
A major factor in the computational cost of flow simulation is the thin boundary layer that develops near a solid wall. Accurately resolving the sharp velocity gradients within this layer requires a mesh that is highly refined in the wall-normal direction, which produces a large number of degrees of freedom (DoFs) and high-aspect-ratio elements. In finite element formulations~\cite{Hughes87Finit}, Dirichlet boundary conditions such as the no-slip condition are typically imposed by prescribing the nodal values of the velocity, or the control variables in the case of Isogeometric Analysis (IGA)~\cite{Hughes05Isoge}. This is referred to as strong enforcement of essential boundary conditions. When the near-wall mesh is too coarse to resolve the boundary layer, strong enforcement constrains the discrete velocity to vanish at the wall even though the adjacent element lacks the resolution needed to capture the steep velocity gradient. As a result, the discrete solution poorly approximates the near-wall flow and produces inaccurate mean-flow quantities throughout the domain~\cite{Bazilevs07Weak1, Bazilevs07Weak2, Hsu12Wind, Xu19resid, Jaiswal26Weak}. For three-dimensional (3D) flows of industrial interest, where optimal mesh design is rarely feasible and resolving the boundary layers is expensive, this limitation is significant.

\subsection{Weak enforcement of essential boundary conditions: core formulations and computations}

Weak imposition of Dirichlet boundary conditions (hereafter referred to as weak BCs) was originally introduced for incompressible flows by Bazilevs and Hughes~\cite{Bazilevs07Weak1} to address this challenge in the context of stabilized and variational multiscale (VMS) finite element methods~\cite{Bazilevs07Varia}. Rather than setting the discrete solution to satisfy the boundary conditions exactly, the weak BC methodology augments the variational formulation with an operator, defined through integrals over the boundary surface, that enforces the Dirichlet data as Euler--Lagrange conditions of the formulation~\cite{Bazilevs07Weak1}. The construction draws on Nitsche's method~\cite{Nitsche71Uber, Stenberg95Onsom, Hansbo02unfit} and on the treatment of boundary and interface conditions in discontinuous Galerkin methods~\cite{Baumann99disco, Riviere01apri, Arnold02Unifi}. Weak enforcement allows the flow to slip at the wall by a controlled amount so that the discretization is not required to resolve the sharp gradients within the boundary layer. The amount of the slip is governed by the local mesh size, velocity, and viscosity~\cite{Bazilevs07Weak1, Bazilevs07Weak2}. The weak BC formulation reduces to its strong counterpart in the limit of vanishing wall-normal mesh size~\cite{Bazilevs07Weak2, Golshan15Large}, allowing mesh refinement to recover the exact Dirichlet condition. 

Numerical studies have shown that weak BCs deliver accuracy comparable to or exceeding that of strongly enforced boundary conditions on substantially coarser boundary-layer meshes~\cite{Bazilevs07Weak1, Bazilevs07Weak2, Hsu12Wind, Xu19resid, Jaiswal26Weak}. This behavior was observed for turbulent channel flow, where the convergence of the mean velocity profile suggested that weak no-slip conditions behave similarly to a wall-function model, although the boundary operator was designed on numerical rather than physical or empirical grounds~\cite{Bazilevs07Weak1}. The connection to near-wall modeling~\cite{Piomelli02Walll, Wilcox06Turbu, Bose18Wallm} was later made explicit by incorporating the law of the wall into the weak operator~\cite{Bazilevs07Weak2, Bazilevs10Isoge} and studied further in the context of large-eddy simulation (LES)~\cite{Golshan15Large}. For incompressible flows, weak BCs combined with residual-based variational multiscale turbulence modeling~\cite{Bazilevs07Varia} and its arbitrary Lagrangian--Eulerian (ALE) extension have enabled accurate, full-scale aerodynamic simulations of wind turbines~\cite{Hsu12Wind, Hsu14Finit}. Extended to the compressible regime~\cite{Xu17Compr, Jaiswal26Weak}, weak BCs have also been successfully applied to aircraft~\cite{Rajanna22Finit, Jaiswal26Weak} and gas-turbine~\cite{Xu17Compr, Kozak20High, Kozak20Optim} aerodynamics, as well as hypersonic flows~\cite{Codoni24Heat}. In each case, weak enforcement predicts quantities of engineering interest on boundary-layer meshes that are far coarser than those required by strong enforcement.

\subsection{Weak BCs in sliding/slip-interface and space--time methods}

Nitsche's method and the weak BC methodology also underpin a family of techniques for coupling non-matching meshes, in which the conditions imposed weakly are those at an interface between two mesh zones~\cite{Becker03finit, Hansbo05Nitsc}. In the ALE setting, this is the sliding-interface method~\cite{Bazilevs08NURBS}, introduced in connection with the ALE-based variational multiscale (ALE--VMS) formulation. Its counterpart in the space--time (ST) framework is the ST Slip Interface (ST-SI) method, introduced for incompressible flows~\cite{Takizawa15Space}, incompressible flows with thermal coupling~\cite{Takizawa16Compu}, and compressible flows~\cite{Takizawa17Poros}. It allows the core ST computational flow analysis methods~\cite{Tezduyar25chron1, Tezduyar25chron2} to be used as moving-mesh methods even in the presence of rotating solid surfaces, such as a car tire. Across a fluid--fluid SI, the mesh on one side rotates with the solid surface, retaining the high-resolution boundary-layer representation, while the mesh on the other side remains independent of the rotation. A fluid--solid SI~\cite{Takizawa15Space, Takizawa16Compu, Takizawa17Poros} is the ST counterpart of the weak BC formulation for non-matching fluid--structure interaction (FSI)~\cite{Hansbo03Nitsc, Bazilevs12Isoge}, and a porosity SI~\cite{Takizawa15Space, Takizawa17Poros} is for a thin porous structure with fluid on its two sides.
 
The ST-SI has since been combined with other ST methods to broaden its scope. Its synthesis with the ST Topology Change (ST-TC) method~\cite{Takizawa14Space} is the ST-SI-TC~\cite{Takizawa16Space}, which enables moving-mesh computation with actual contact between solid surfaces while retaining the high-resolution boundary-layer representation. The first 3D computation with the ST-SI-TC was for tire aerodynamics with transverse grooves, road contact, and tire deformation~\cite{Takizawa16Space}. Adding the ST-IGA~\cite{Takizawa11Multi}, which is IGA discretization in space, in time, or in both, gives the ST-SI-TC-IGA~\cite{Takizawa17Heart}. Problems of this class often involve narrow, curved gaps between solid surfaces and elements with high aspect ratios. The IGA discretization achieves good accuracy in such gaps at a reasonable element density, and therefore at a reasonable computational cost, while also being more robust in computations with high-aspect-ratio elements. With the feature introduced in Ref.~\cite{Kuraishi19Space}, the ST-SI-TC-IGA has a built-in Reynolds-equation limit for fluid films between surfaces coming into contact. It has been used for aortic valve flow with leaflet contact~\cite{Takizawa17Heart}, for tire aerodynamics with near-actual geometry and, more recently, with a complex asymmetric tread pattern~\cite{Kuraishi19Tire, Kuraishi25Space}, and for aortic-valve-to-aorta flow~\cite{Terahara26Space}. Dirichlet-boundary stabilization methods, for both weak and strong enforcement, have also been derived in the context of 1D elastodynamics to improve the accuracy of the boundary fluxes computed with higher-order functions in the ST framework~\cite{Taniguchi26Diric}.

\subsection{Weak BCs in immersed methods}\label{sec:intro_immersed}

Another role for weak BCs arises in immersed methods~\cite{Peskin02immer, Mittal05Immer, Burman15CutFE, Schillinger15finit, Griffith20Immer}. In the fluid setting, the Navier--Stokes equations are solved on a background mesh that does not conform to the geometry~\cite{Ruberg12Subdi, Schott14new, Kamensky15immer, Xu16tetra, Main18shift2, Karki25Direc}. Immersed methods relax the mesh-conforming constraint and simplify mesh generation around geometrically complex objects, but strong enforcement of wall boundary conditions becomes infeasible. Because the immersed boundary cuts arbitrarily through element interiors, the Dirichlet data cannot be imposed by directly prescribing nodal values, as is done on boundary-fitted meshes. The weak BC methodology resolves this issue directly, since its boundary operator is expressed in terms of the surface traces of the discrete fields and is independent of the mesh. However, standard numerical integration over these arbitrarily intersected elements remains a problem~\cite{Parvizian07Finit, Duster08finit}. Immersogeometric analysis (IMGA)~\cite{Kamensky15immer, Xu16tetra} was proposed to address both challenges by combining weak BCs with geometry-aware integration over the cut elements, as resolving them with sufficient geometric fidelity was shown to be essential for accurate flow solutions. IMGA was later driven directly by boundary-representation (B-rep) CAD models~\cite{Hsu16Direc, Wang17Rapid}, removing the labor-intensive geometry cleanup and meshing steps from the design-through-analysis workflow~\cite{Schillinger12isoge, Wassermann19Integ}. The framework has since been extended to compressible flows and rotorcraft~\cite{Xu19Immer1}, free-surface and marine flows~\cite{Zhu20immer}, moving-object problems~\cite{Xu19Immer2, Xu21octre}, and industrial-scale simulations~\cite{Saurabh21Indus}.

Because the weak BC operator requires only a geometric description, rather than a discrete surface mesh, to evaluate boundary traces, IMGA has recently been applied directly to sampled representations of physical objects. In this setting, a sufficiently dense point cloud of the surface is used to enforce the no-slip condition weakly, and the simulation proceeds without reconstructing an analysis-suitable CAD model or surface mesh~\cite{Balu23Direc}. This capability has enabled computational flow analysis directly on point clouds obtained from photogrammetry for in-use civil structures~\cite{Wang23Photo}, from scans of real-world objects subject to noise, incompleteness, and irregular sampling~\cite{Jaiswal24Mesh}, and from segmented medical images for patient-specific cardiovascular simulation~\cite{Corpuz25Direc}. In each of these developments, the weak imposition of the no-slip condition is the underlying mechanism that makes flow analysis possible on geometries that are never meshed and, in some cases, never fully reconstructed. 

This review unifies these developments and presents them as successive applications of a single core concept. We begin with a derivation of weak boundary conditions for the scalar advection and diffusion problems, combine them into the advection--diffusion equation, and illustrate the practical behavior of this baseline formulation with numerical examples. We then extend this mathematical structure to the incompressible Navier--Stokes equations, reviewing its application in both boundary-fitted and immersogeometric contexts. Following this, we introduce a generalized nonlinear flux framework, which serves as the foundation for our weak boundary condition formulation for compressible flows. Alongside these mathematical developments, we review the relevant literature, organizing our discussion around the two primary roles established earlier: acting as a near-wall model and serving as an enabling technology for direct geometry-to-analysis workflows. We focus throughout on computational flow analysis and compare weakly and strongly enforced boundary conditions to demonstrate the practical advantages of these formulations. For a systematic construction of the underlying Nitsche's method across general variational problems, we refer the reader to Ref.~\cite{Benzaken24Const}.

The remainder of the paper is organized as follows. Section~\ref{sec:scalar} develops the weak BC operator for scalar model problems and reviews its behavior on the advection--diffusion equation. Section~\ref{sec:incomp} presents the incompressible Navier--Stokes operator and reviews its application to wall-bounded turbulence, wind-turbine aerodynamics, and immersogeometric analysis. Section~\ref{sec:comp} details the nonlinear flux formulation, applies it to compressible flows, and reviews relevant boundary-fitted and immersogeometric results. Section~\ref{sec:conclusions} draws conclusions.

\section{Weak imposition of boundary conditions: scalar model problems}
\label{sec:scalar}

In this section, we derive the weak BC operator for scalar model problems, beginning with the pure advection and pure diffusion equations and then combining them into the advection--diffusion problem. The two equations contribute different components. The advection equation contributes a flux term that is active on the inflow boundary, stabilizing the formulation where the flow enters the domain. The diffusion equation contributes the natural consistency flux together with the adjoint-consistency and penalty/stabilization terms of Nitsche's method~\cite{Benzaken24Const}. In the combined advection--diffusion equation, these terms act together, and the local mesh size and the relative strength of advection and diffusion control how closely the discrete solution satisfies $u=g$ on the boundary. 

\subsection{Advection equation}
\label{sec:advection}
Consider the pure advection equation
\begin{align}
\label{eq:advection-eqn}
\nabla \cdot \left(\mathbf{a} u \right) - f &= 0 \quad \text{in}\;\; \Omega \text{ ,}\\
\label{eq:inflow}
u &= g \quad \text{on}\;\; \Gamma_{-} \text{ ,}
\end{align}
where $\mathbf{a}$ denotes a prescribed advective velocity field (assumed solenoidal), $f$ is a source term, and $u$ is the unknown scalar quantity being transported on a domain $\Omega \subset \mathbb{R}^2$ or $\mathbb{R}^3$ with boundary $\partial\Omega = \Gamma$. The inflow boundary is defined as $\Gamma_{-} = \left\{\,\mathbf{x} \in \Gamma \mid \mathbf{a} \cdot \mathbf{n} < 0 \, \right\} $, where $\mathbf{x}$ denotes a spatial point and $\mathbf{n}$ is the outward unit normal to $\Gamma$. The outflow boundary is given by $\Gamma_{+} = \Gamma \setminus \Gamma_{-}$. Accordingly, $\{\cdot\}_{-}$ and $\{\cdot\}_{+}$ denote the negative and positive parts of their argument, respectively.

To derive the weak formulation, we multiply Eq.~\eqref{eq:advection-eqn} by a test function $w$ and integrate over the domain $\Omega$:
\begin{align}
\label{eq:adv-weak-form-1}
\int_{\Omega} w \, \nabla \cdot \left(\mathbf{a} u \right) \,\mathrm{d}\Omega - \int_{\Omega} wf \,\mathrm{d}\Omega = 0 \text{ .}
\end{align}
The advection term in Eq.~\eqref{eq:adv-weak-form-1} can be integrated by parts, yielding 
\begin{align}
\label{eq:adv-ibp}
\int_{\Omega} w \, \nabla \cdot \left(\mathbf{a} u \right) \,\mathrm{d}\Omega =
- \int_{\Omega} \nabla w \cdot \mathbf{a} u \,\mathrm{d}\Omega 
+ \int_{\Gamma} w \left( \mathbf{a} \cdot \mathbf{n} \right) u \,\mathrm{d}\Gamma \text{ .}
\end{align}
Substituting this result back into Eq.~\eqref{eq:adv-weak-form-1} leads to the weak form%
\begin{align}
\label{eq:adv-weak-form-2}
- \int_{\Omega} \nabla w \cdot \mathbf{a} u \,\mathrm{d}\Omega 
+ \int_{\Gamma} w \left( \mathbf{a} \cdot \mathbf{n} \right) u \,\mathrm{d}\Gamma 
- \int_{\Omega}wf \,\mathrm{d}\Omega = 0 \text{ .}
\end{align}
On the boundary $\Gamma$, the second term in Eq.~\eqref{eq:adv-weak-form-2} naturally separates into contributions from the inflow and outflow regions: 
\begin{align}
\label{eq:adv-weak-form-3}
- \int_{\Omega} \nabla w \cdot \mathbf{a} u \,\mathrm{d}\Omega 
+ \int_{\Gamma_{-}} w \left\{ \mathbf{a} \cdot \mathbf{n} \right\}_{-} u \,\mathrm{d}\Gamma 
+ \int_{\Gamma_{+}} w \left\{ \mathbf{a} \cdot \mathbf{n} \right\}_{+} u \,\mathrm{d}\Gamma 
- \int_{\Omega}wf \,\mathrm{d}\Omega = 0 \text{ .}
\end{align}
Substituting Eq.~\eqref{eq:inflow} into Eq.~\eqref{eq:adv-weak-form-3} on the inflow boundary $\Gamma_{-}$ yields
\begin{align}
\label{eq:adv-weak-form-4}
- \int_{\Omega} \nabla w \cdot \mathbf{a} u \,\mathrm{d}\Omega 
+ \int_{\Gamma_{-}} w \left\{ \mathbf{a} \cdot \mathbf{n} \right\}_{-} g \,\mathrm{d}\Gamma 
+ \int_{\Gamma_{+}} w \left\{ \mathbf{a} \cdot \mathbf{n} \right\}_{+} u \,\mathrm{d}\Gamma 
- \int_{\Omega}wf \,\mathrm{d}\Omega = 0 \text{ .}
\end{align}

To retain the advection term in its original divergence form, $\nabla \cdot \left( \mathbf{a} u \right)$, we reverse the integration by parts in Eq.~\eqref{eq:adv-weak-form-4} using Eq.~\eqref{eq:adv-ibp}:
\begin{align}
\label{eq:adv-weak-form-5}
&\int_{\Omega} w \, \nabla \cdot \left(\mathbf{a} u \right) \,\mathrm{d}\Omega
- \int_{\Gamma} w \left( \mathbf{a} \cdot \mathbf{n} \right) u \,\mathrm{d}\Gamma
+ \int_{\Gamma_{-}} w \left\{ \mathbf{a} \cdot \mathbf{n} \right\}_{-} g \,\mathrm{d}\Gamma 
+ \int_{\Gamma_{+}} w \left\{ \mathbf{a} \cdot \mathbf{n} \right\}_{+} u \,\mathrm{d}\Gamma 
\nonumber \\
&\quad
- \int_{\Omega}wf \,\mathrm{d}\Omega = 0 \text{ .}
\end{align}
Splitting the boundary integral over $\Gamma$ in the second term into its inflow and outflow contributions, the outflow terms cancel exactly, and only the inflow contribution remains:
\begin{align}
\label{eq:adv-weak-form-final}
\int_{\Omega} w \, \nabla \cdot \left(\mathbf{a} u \right) \,\mathrm{d}\Omega
- \int_{\Gamma_{-}} w \left\{ \mathbf{a} \cdot \mathbf{n} \right\}_{-} \left( u - g \right) \,\mathrm{d}\Gamma 
- \int_{\Omega}wf \,\mathrm{d}\Omega = 0 \text{ .}
\end{align}
\begin{remark}
The final weak formulation in Eq.~\eqref{eq:adv-weak-form-final} is consistent with the strong form, Eq.~\eqref{eq:advection-eqn}: the inflow term vanishes whenever the Dirichlet condition $u=g$ is satisfied on the inflow boundary $\Gamma_{-}$. The inflow term adds stability: setting $w=u$ in the homogeneous case ($g=0$), its contribution to the bilinear form is $-\int_{\Gamma_{-}} \left\{ \mathbf{a} \cdot \mathbf{n} \right\}_{-} u^2 \,\mathrm{d}\Gamma \geq 0$, since $\mathbf{a} \cdot \mathbf{n} < 0$ on $\Gamma_{-}$. The term therefore adds a non-negative, dissipative contribution to the energy estimate on the inflow boundary, supplying the control over the inflow data that the advection operator alone does not provide.
\end{remark}

\subsection{Diffusion equation}
Consider the steady-state diffusion problem
\begin{align}
\label{eq:diff-eqn}
- \nabla \cdot \left( \kappa \nabla u \right) - f &= 0 \quad \text{in}\;\; \Omega \text{ ,} \\
u &=g \quad \text{on}\;\; \Gamma \text{ ,}
\end{align}
where $\Omega \subset \mathbb{R}^2$ or $\mathbb{R}^3$ is a bounded domain with boundary $\partial\Omega = \Gamma$, $\kappa$ is the diffusivity, and $f$ is a source term. For elliptic problems with $\kappa > 0$, it is standard practice to impose Dirichlet conditions on the entire boundary $\Gamma$. However, rather than enforcing the nodal values of $u$ directly, we impose these conditions weakly. To derive the weak formulation, we multiply Eq.~\eqref{eq:diff-eqn} by a test function $w$, integrate over $\Omega$, and apply integration by parts:
\begin{align}
\label{eq:diff-weak-form-1}
\int_{\Omega} \nabla w \cdot \kappa \nabla u \,\mathrm{d}\Omega 
- \int_{\Omega} w f \,\mathrm{d}\Omega 
- \int_{\Gamma} w \kappa \nabla u \cdot \mathbf{n} \,\mathrm{d}\Gamma = 0 \text{ ,}
\end{align}
where $\mathbf{n}$ is the outward unit normal to $\Gamma$. The boundary integral in Eq.~\eqref{eq:diff-weak-form-1} is the natural (Neumann) boundary contribution that arises from integration by parts. To impose the Dirichlet condition $u=g$ weakly while preserving consistency, we add an adjoint-consistency term and a stabilization term, which yields
\begin{align}
\label{eq:diff-weak-form-nitsche}
&\int_{\Omega} \nabla w \cdot \kappa \nabla u \,\mathrm{d}\Omega 
- \int_{\Omega} w f \,\mathrm{d}\Omega 
- \int_{\Gamma} w \kappa \nabla u \cdot \mathbf{n} \,\mathrm{d}\Gamma 
- \int_{\Gamma} \kappa \nabla w \cdot \mathbf{n} \left( u - g \right) \,\mathrm{d}\Gamma 
\nonumber \\
&\quad
+ \int_{\Gamma} w \tau_\mathrm{B} \left( u - g \right) \,\mathrm{d}\Gamma = 0 \text{ .}
\end{align}
Here, the first three terms are the weak form of the partial differential equation (PDE), while the last two terms impose the Dirichlet condition weakly. The second-to-last term is the adjoint-consistency term, which renders the formulation symmetric. The last term is a penalty/stabilization term that penalizes violation of the Dirichlet condition and is essential for the coercivity of the bilinear form, with $\tau_\mathrm{B}$ the stabilization parameter. The choice of $\tau_\mathrm{B}$ is a compromise: the parameter must be large enough to guarantee stability, but not so large that the operator behaves like a pure penalty method, overshadowing the variational consistency that distinguishes Nitsche's method from a penalty method and degrading the conditioning of the resulting linear equation systems~\cite{Nitsche71Uber, Embar10Impos, dePrenter18note}.

\subsection{Selection of $\tau_\mathrm{B}$ using stability analysis}
\label{sec:stability-analysis}
To determine an appropriate scaling of the stabilization parameter $\tau_\mathrm{B}$, we examine the stability of the formulation in Eq.~\eqref{eq:diff-weak-form-nitsche}. To this end, we reformulate the problem as: find $u$ such that for all $w$,
\begin{equation}
    \label{eq:bilinear-form}
    B(w,u) - F(w) = 0 \text{ ,}
\end{equation}
where $B$ is a bilinear form and $F$ a linear functional, given by
\begin{align}
B(w, u)
&= \int_{\Omega} \nabla w \cdot \kappa \nabla u \,\mathrm{d}\Omega 
- \int_{\Gamma} w \kappa \nabla u \cdot \mathbf{n} \,\mathrm{d}\Gamma 
-  \int_{\Gamma} \kappa \nabla w \cdot \mathbf{n} u \,\mathrm{d}\Gamma 
+ \int_{\Gamma} w \tau_\mathrm{B} u \,\mathrm{d}\Gamma \text{ ,}
\label{eq:diff-B-form}
\\[4pt]
F(w)
&= \int_{\Omega} w f \,\mathrm{d}\Omega 
- \int_{\Gamma} \kappa \nabla w \cdot \mathbf{n} g \,\mathrm{d}\Gamma
+ \int_{\Gamma} w \tau_\mathrm{B} g \,\mathrm{d}\Gamma \text{ .}
\label{eq:diff-F-form}
\end{align}

The well-posedness of the problem relies on the coercivity of $B(w,u)$. Specifically, we seek a norm $\triplebar \cdot \triplebar$ and a mesh-independent constant $C > 0$ such that
\begin{equation}
\label{eq:coercivity}
B(u,u) \geq C \, \triplebar u \triplebar^2 \text{ .}
\end{equation}
The choice of an energy-equivalent norm is guided by the structure of $B(w,u)$. Setting $w = u$ in Eq.~\eqref{eq:diff-B-form} yields
\begin{align}
\label{eq:diff-B-norm}
B(u,u) &= \underbrace{\int_{\Omega} \kappa \, \vert \nabla u \vert^2 \,\mathrm{d}\Omega}_{(\mathrm{i})} 
        - \underbrace{2 \int_{\Gamma} u \, \kappa \nabla u \cdot \mathbf{n} \,\mathrm{d}\Gamma}_{(\mathrm{ii})} 
        + \underbrace{\int_{\Gamma} \tau_{\mathrm{B}} \, u^2 \,\mathrm{d}\Gamma}_{(\mathrm{iii})} \text{ .}
\end{align}
Provided $\kappa > 0$ and $\tau_\mathrm{B} > 0$, terms~$(\mathrm{i})$ and~$(\mathrm{iii})$ are non-negative and contribute directly to coercivity. This motivates defining the discrete energy norm
\begin{equation}
\label{eq:triple-norm}
  \triplebar u \triplebar^2 = \int_{\Omega} \kappa \, |\nabla u|^2 \,\mathrm{d}\Omega + \int_{\Gamma} \tau_{\mathrm{B}} \, u^2 \,\mathrm{d}\Gamma \text{ .}
\end{equation}
Note that Eq.~\eqref{eq:triple-norm} defines a valid norm as it satisfies all the defining properties of a norm. The remaining term~$(\mathrm{ii})$, however, lacks a definite sign. We can perform the following steps to demonstrate that term~$(\mathrm{ii})$ can be bounded by terms~$(\mathrm{i})$ and~$(\mathrm{iii})$.

\paragraph{Step 1: Application of Young's inequality}
Applying Young's inequality
\begin{equation}
2\vert a b \vert \leq a^2/\varepsilon + \varepsilon b^2
\end{equation}
for any $\varepsilon > 0$ pointwise on $\Gamma$ with $a = u$ and $b = \nabla u \cdot \mathbf{n}$, we establish the lower bound
\begin{equation}
\label{eq:young}
-2 \int_{\Gamma} u \, \kappa \nabla u \cdot \mathbf{n} \,\mathrm{d}\Gamma \geq \,
-\underbrace{\int_{\Gamma} \frac{\kappa}{\varepsilon} \, u^2 \,\mathrm{d}\Gamma}_{(\mathrm{iv})} \,
-\underbrace{\int_{\Gamma} \kappa \, \varepsilon \, (\nabla u \cdot \mathbf{n})^2 \,\mathrm{d}\Gamma}_{(\mathrm{v})} \text{ .}
\end{equation}
Term~$(\mathrm{iv})$ shares the same functional form as the penalty term~$(\mathrm{iii})$ and can be directly absorbed. However, term~$(\mathrm{v})$ involves the normal flux on the boundary, which is not bounded by the norm in Eq.~\eqref{eq:triple-norm}.

\paragraph{Step 2: Trace inverse estimate}
To bound term~$(\mathrm{v})$, we employ a trace inverse inequality. For piecewise polynomial functions, there exists a mesh-independent constant $C_\mathrm{I} > 0$, dependent only on the polynomial degree and element type, such that
\begin{equation}
\label{eq:inverse-estimate}
\int_{\Gamma} (\nabla u \cdot \mathbf{n})^2 \,\mathrm{d}\Gamma \leq \frac{C_\mathrm{I}}{h} \int_{\Omega} \vert \nabla u \vert^2 \,\mathrm{d}\Omega \text{ ,}
\end{equation}
where $h$ denotes the characteristic element size, evaluated element-wise; Eq.~\eqref{eq:inverse-estimate} is stated in global form for notational simplicity.

Combining inequalities~\eqref{eq:young} and~\eqref{eq:inverse-estimate} yields
\begin{equation}
\label{eq:combined-bound}
-2 \int_\Gamma u \, \kappa \nabla u \cdot \mathbf{n} \,\mathrm{d}\Gamma \geq -\int_\Gamma \frac{\kappa}{\varepsilon} \, u^2 \,\mathrm{d}\Gamma - \varepsilon \, \frac{C_\mathrm{I}}{h} \int_\Omega \kappa \, \vert \nabla u \vert^2 \,\mathrm{d}\Omega \text{ .}
\end{equation}
Substituting Eq.~\eqref{eq:combined-bound} into Eq.~\eqref{eq:diff-B-norm} and grouping like terms gives
\begin{equation}
\label{eq:absorbed}
B(u, u) \geq \int_\Omega \kappa \left( 1 - \varepsilon \, \frac{C_\mathrm{I}}{h} \right) \vert \nabla u \vert^2 \,\mathrm{d}\Omega + 
\int_\Gamma \left( \tau_\mathrm{B} - \frac{\kappa}{\varepsilon} \right) u^2 \,\mathrm{d}\Gamma \text{ .}
\end{equation}
To ensure coercivity, we therefore need
\begin{equation}
\varepsilon < \frac{h}{C_\mathrm{I}}
\end{equation}
and
\begin{equation}
\tau_\mathrm{B} > \frac{\kappa}{\varepsilon} \text{ ,}
\end{equation}
which implies that
\begin{equation}
\tau_\mathrm{B} > \frac{C_\mathrm{I} \, \kappa}{h} \text{ .}
\end{equation}
Therefore, setting $\tau_\mathrm{B} = \dfrac{C_\mathrm{P} \, \kappa}{h}$ with $C_\mathrm{P} > C_\mathrm{I}$ guarantees stability.

\subsection{Combining the advection and diffusion equations}
The strong form of the steady-state advection--diffusion problem is given as
\begin{align}
\label{eq:adv-diff-strong}
\nabla \cdot \left(\mathbf{a} u \right) - \nabla \cdot \left( \kappa \nabla u \right) - f &= 0 \quad \text{in}\;\; \Omega \text{ ,}\\
u &=g \quad \text{on}\;\; \Gamma \text{ .}
\end{align}
Combining Eqs.~\eqref{eq:adv-weak-form-final} and~\eqref{eq:diff-weak-form-nitsche}, the weak formulation for the advection--diffusion equation is obtained as follows: find $u$ such that for all $w$,
\begin{align}
\label{eq:adv-diff-weak}
& \int_{\Omega} w \left( \nabla \cdot \left( \mathbf{a} u \right) - f \right) \,\mathrm{d}\Omega 
 + \int_{\Omega} \nabla w \cdot \kappa \nabla u \,\mathrm{d}\Omega 
 - \int_{\Gamma} w \kappa \nabla u \cdot \mathbf{n} \,\mathrm{d}\Gamma 
\nonumber \\
& \quad
- \int_{\Gamma_{-}} \left( w \left\{ \mathbf{a} \cdot \mathbf{n} \right\}_{-} + \kappa \nabla w \cdot \mathbf{n} \right) \left( u - g \right) \,\mathrm{d}\Gamma 
- \int_{\Gamma_{+}} \kappa \nabla w \cdot \mathbf{n} \left( u - g \right) \,\mathrm{d}\Gamma
\nonumber \\
& \quad
+ \int_{\Gamma} w \frac{C_\mathrm{P}\, \kappa}{h} \left( u - g \right) \,\mathrm{d}\Gamma = 0 \text{ .}
\end{align}

\begin{remark}
The combined operator in Eq.~\eqref{eq:adv-diff-weak} adjusts the strength with which the Dirichlet condition $u=g$ is enforced according to the local mesh size $h$ and the relative importance of advection and diffusion, so that the appropriate balance is selected automatically rather than prescribed manually. The constant $C_\mathrm{P}$ must be taken large enough to retain stability and coercivity, with its lower bound computed from an element-level inverse estimate~\cite{Ciarlet78Finit, Harari92What, Warburton03Oncon, Embar10Impos, Evans13Expli}. On the inflow boundary $\Gamma_{-}$, the advective contribution adds to the adjoint-consistency and penalty terms, so the boundary data is enforced more strongly where the flow enters the domain. In the advective limit $\kappa \to 0$, the diffusive adjoint-consistency and penalty terms vanish, the formulation reduces smoothly to the pure-advection operator of Eq.~\eqref{eq:adv-weak-form-final}, and the condition $u=g$ is no longer enforced on the outflow boundary $\Gamma_{+}$.
\end{remark}

\subsection{Numerical examples: the advection--diffusion equation}
\label{sec:adv-diff-results}

We begin by examining the behavior of the weak BC operator on the scalar advection--diffusion equation (see Eq.~\eqref{eq:adv-diff-weak}). The original study of Ref.~\cite{Bazilevs07Weak1} compared the weakly and strongly enforced formulations on a sequence of 1D, 2D, and 3D boundary-layer problems. The 1D and 2D examples consider the transport of a scalar whose exact solution has a thin boundary layer near the wall that the mesh is too coarse to resolve. With strong enforcement, requiring the discrete solution to attain the prescribed value on the wall over an under-resolved boundary layer produces spurious oscillations that propagate into the interior. With weak enforcement through the operator of Eq.~\eqref{eq:adv-diff-weak}, the discrete solution is allowed to depart from the boundary value by an amount controlled by the local mesh size and the relative strength of advection and diffusion. The spurious oscillations are reduced or removed, and the solution away from the boundary is more accurate as a result.

\begin{figure}[!h]
  \centering 
  \includegraphics[width=\textwidth]{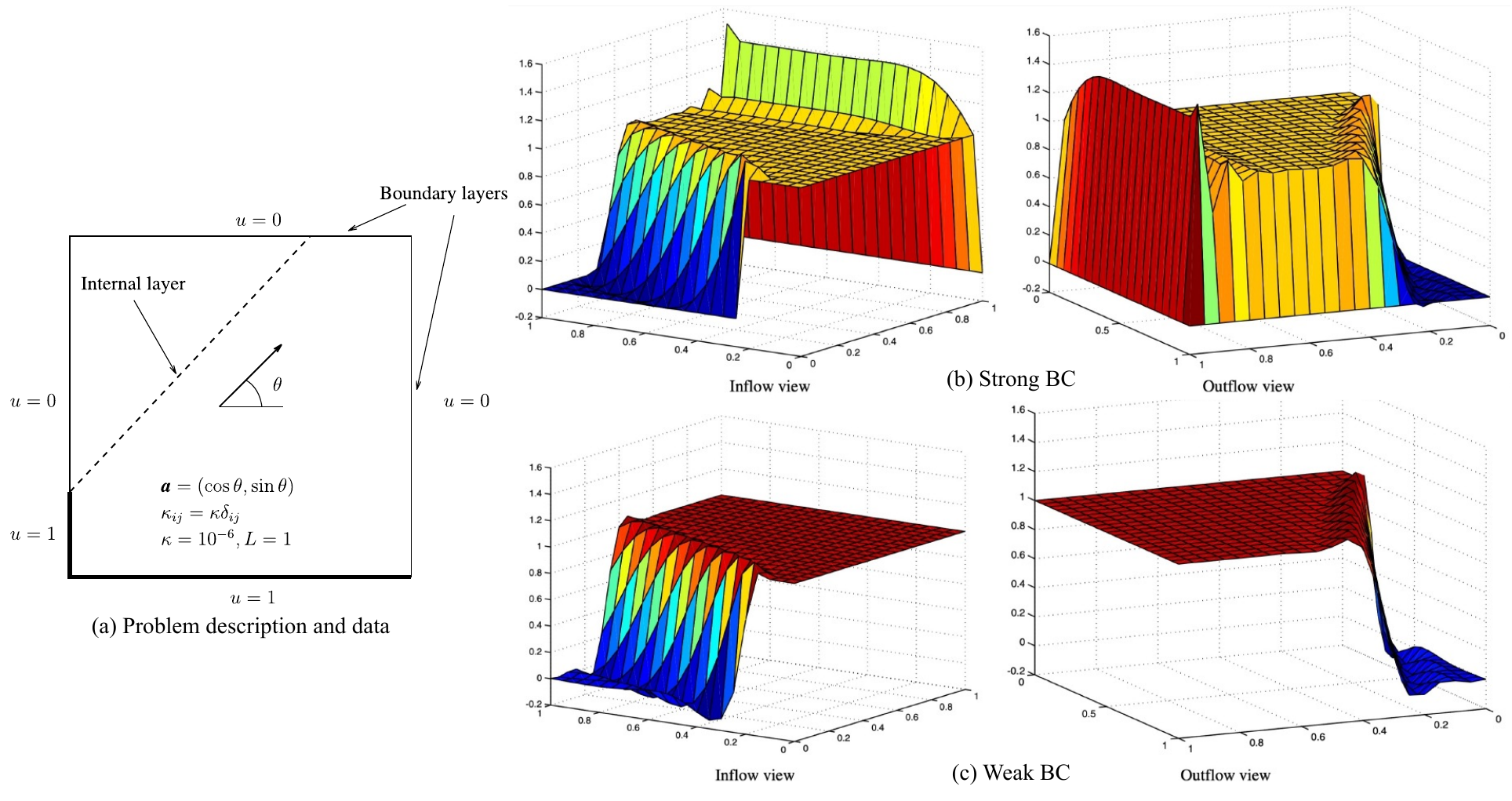}
  \caption{Advection skew to the mesh: strongly enforced versus weakly enforced Dirichlet boundary conditions. Adapted from Ref.~\cite{Bazilevs07Weak1}.}
  \label{fig:advskew}
\end{figure}

A representative 2D case is advection skew to the mesh, where the diffusivity is small, making the problem advection-dominated. The presence of unresolved interior and boundary layers poses difficulties for most existing techniques, and oscillations are often observed in their vicinity. Comparing computations in which all Dirichlet conditions are enforced strongly with those in which they are enforced weakly (Figure~\ref{fig:advskew}), we observe that strong enforcement produces a spurious overshoot at the outflow boundary, where the boundary layer is unresolved on the coarse mesh. The computed result exceeds the exact solution by more than 50\%. Weak enforcement successfully removes this overshoot at the outflow boundary, while reproducing the inflow data fairly well, with only a slight oscillation in the region of the discontinuity.

\begin{figure}[!t]
  \centering
  \includegraphics[width=\textwidth]{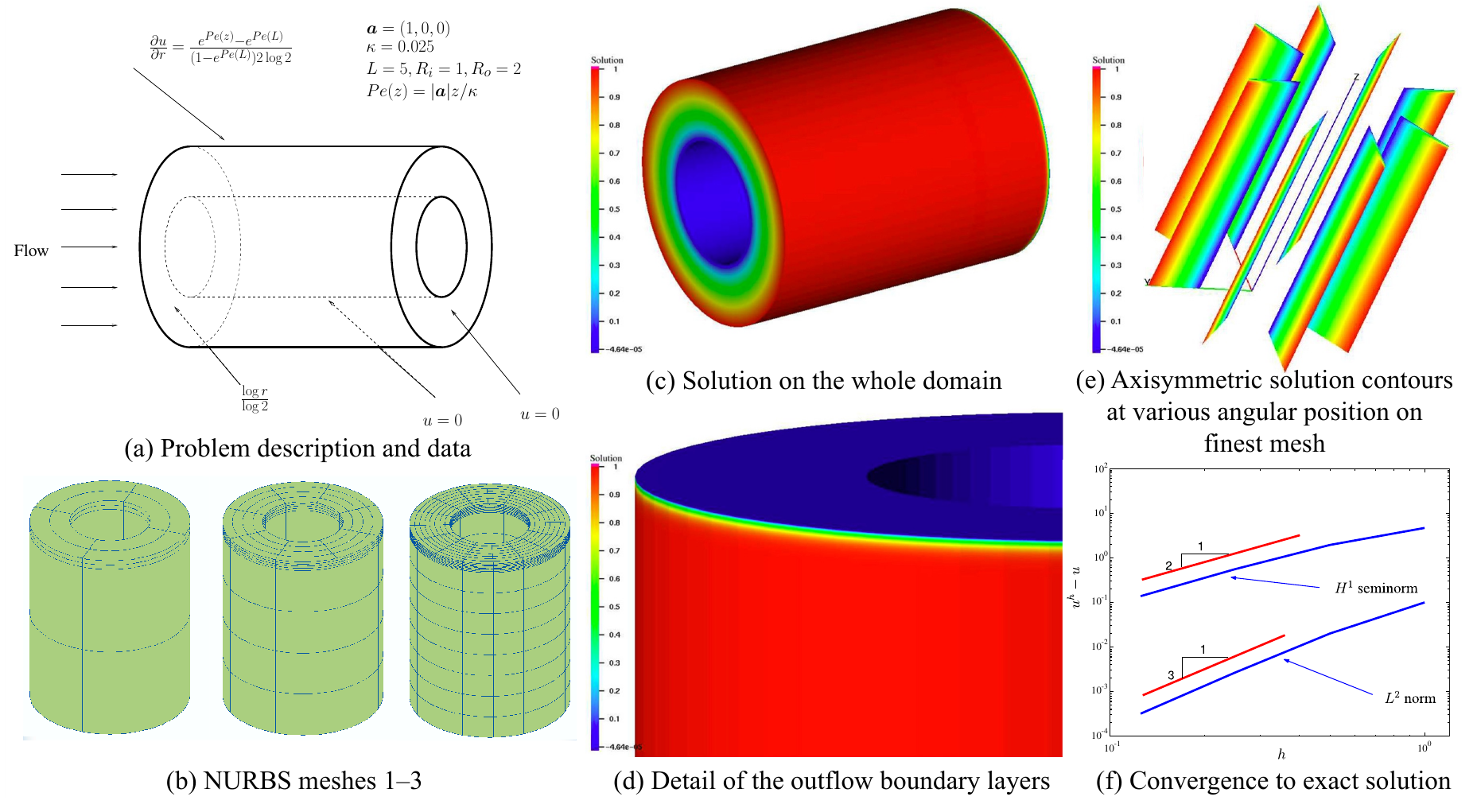}
  \caption{Advection--diffusion in an annular region: problem setup, NURBS meshes, and solutions. Adapted from Ref.~\cite{Bazilevs07Weak1}.}
  \label{fig:annular}
\end{figure}

\newpage
The 3D example is posed on an annular region (Figure~\ref{fig:annular}), for which an exact geometric description of the hollow cylinder is obtained using a quadratic NURBS basis in the isogeometric framework~\cite{Hughes05Isoge}. The Dirichlet conditions are imposed weakly on a sequence of meshes that are locally refined in the thin outflow layer. The weak formulation yields a stable solution that recovers the pointwise axisymmetric response expected for the problem, and the error measured against the analytical solution in the $L^2$-norm and $H^1$-seminorm shows optimal convergence rates under mesh refinement.

Two features of these results carry over to all subsequent formulations. First, the weak operator is variationally consistent: the exact solution satisfies the discrete formulation, so weak enforcement introduces no modeling error, and refinement recovers the strongly enforced result. Second, the advantage of weak enforcement is largest where strong enforcement performs worst---on meshes too coarse to resolve the boundary layer. Both properties hold for the wall-bounded turbulent flows introduced later, where the boundary layer is physical rather than a feature of a manufactured solution.

\section{Weak BCs for incompressible flow}
\label{sec:incomp}

The scalar developments of the previous section carry over directly to the no-slip condition for the incompressible Navier--Stokes equations, the setting in which weak BCs have seen their widest use.

\subsection{The incompressible Navier--Stokes weak BC operator}
\label{sec:incomp-ns}

Let $\Omega \subset \mathbb{R}^d$, $d=2$ or $3$, be the spatial domain with boundary $\Gamma = \partial\Omega$. Let $\Gamma^\mathrm{D} \subset \Gamma$ denote the portion on which the velocity is prescribed and $\Gamma^\mathrm{N} \subset \Gamma$ be the Neumann boundary. The strong form of the problem can be written as
\begin{align}
\label{eq:ins-strong-mom}
\rho \left( \frac{\partial\mathbf{u}}{\partial t} + \mathbf{u}\cdot\nabla\mathbf{u} - \mathbf{f} \right)
- \nabla\cdot\boldsymbol{\sigma}\left(\mathbf{u},p\right) &= \mathbf{0}
\quad \text{in}\;\; \Omega \text{ ,} \\
\label{eq:ins-strong-cont}
\nabla\cdot\mathbf{u} &= 0
\quad \text{in}\;\; \Omega \text{ ,} \\
\label{eq:ins-strong-bc}
\mathbf{u} &= \mathbf{g}
\quad \text{on}\;\; \Gamma^\mathrm{D} \text{ ,}\\
\label{eq:ins-neumann-bc}
-p\, \mathbf{n} + 2\mu \,  \boldsymbol{\varepsilon}(\mathbf{u})\, \mathbf{n}  &= \mathbf{h}
\quad \text{on}\;\; \Gamma^\mathrm{N}\text{ ,}
\end{align}
where $\mathbf{u}$ is the velocity, $p$ is the pressure, $\rho$ is the (constant) density, $\mathbf{f}$ is the body force per unit mass, $\mathbf{h}$ is the prescribed traction on $\Gamma^\mathrm{N}$, and $\boldsymbol{\sigma}\left(\mathbf{u},p\right) = -p\mathbf{I} + 2\mu\boldsymbol{\varepsilon}\left(\mathbf{u}\right)$ is the stress tensor, with $\mathbf{I}$ the identity tensor, $\mu$ the dynamic viscosity, and $\boldsymbol{\varepsilon}\left(\mathbf{u}\right) = \frac{1}{2}\left(\nabla\mathbf{u} + \left(\nabla\mathbf{u}\right)^{\mathrm{T}} \right)$ the strain-rate tensor. On a no-slip boundary, $\mathbf{g}$ is the prescribed velocity, which is zero for a stationary wall.

The discrete problem is posed over finite element (or, in IGA, NURBS) spaces and is stabilized using the residual-based VMS formulation~\cite{Bazilevs07Varia}, which on moving domains is recast in its ALE--VMS form~\cite{Bazilevs08NURBS, Hsu14Finit}. The VMS framework decomposes the velocity and pressure into coarse- and fine-scale parts, $\mathbf{u} = \mathbf{u}^h + \mathbf{u}^{\prime}$ and $p = p^h + p^{\prime}$, where the coarse scales are represented by the discrete spaces and the fine scales are modeled. Let $B^{\mathrm{VMS}}\left(\{\mathbf{w}^h,q^h\}, \{\mathbf{u}^h,p^h\}\right)$ and $F^{\mathrm{VMS}}\left(\{\mathbf{w}^h,q^h\}\right)$ denote, respectively, the semi-linear form and the linear functional of the semi-discrete VMS formulation, where $\{\mathbf{u}^h,p^h\}$ are the coarse-scale velocity and pressure trial solutions and $\{\mathbf{w}^h,q^h\}$ the corresponding test functions. The fluid domain $\Omega$ is divided into individual spatial finite element subdomains $\Omega^e$, and the detailed expressions for $B^{\mathrm{VMS}}$ and $F^{\mathrm{VMS}}$ are
\begin{align} \label{eq:B}
\nonumber & B^{\mathrm{VMS}}\left(\{\mathbf{w}^h, q^h\},\{\mathbf{u}^h, p^h\}\right) =\int_{\Omega} \mathbf{w}^h \cdot \rho\left(\frac{\partial \mathbf{u}^h}{\partial t}+\mathbf{u}^h \cdot \nabla \mathbf{u}^h\right) \,\mathrm{d}\Omega +\int_{\Omega} \boldsymbol{\varepsilon}(\mathbf{w}^h): \boldsymbol{\sigma}(\mathbf{u}^h, p^h) \,\mathrm{d}\Omega \\
\nonumber & \quad +\int_{\Omega} q^h \nabla \cdot \mathbf{u}^h \,\mathrm{d}\Omega -\sum_e \int_{\Omega^e}\left(\mathbf{u}^h \cdot \nabla \mathbf{w}^h+\frac{\nabla q^h}{\rho}\right) \cdot \mathbf{u}^{\prime} \,\mathrm{d}\Omega -\sum_e \int_{\Omega^e} p^{\prime} \nabla \cdot \mathbf{w}^h \,\mathrm{d}\Omega \\
\nonumber & \quad+\sum_e \int_{\Omega^e} \mathbf{w}^h \cdot\left(\mathbf{u}^{\prime} \cdot \nabla \mathbf{u}^h\right) \,\mathrm{d}\Omega -\sum_e \int_{\Omega^e} \frac{\nabla \mathbf{w}^h}{\rho}:\left(\mathbf{u}^{\prime} \otimes \mathbf{u}^{\prime}\right) \,\mathrm{d}\Omega \\
& \quad+\sum_e \int_{\Omega^e}\left(\mathbf{u}^{\prime} \cdot \nabla \mathbf{w}^h\right) \cdot \bar{\tau} \left(\mathbf{u}^{\prime} \cdot \nabla \mathbf{u}^h\right) \,\mathrm{d}\Omega\text{ ,}
\end{align}
and
\begin{align}\label{eq:F}
F^\mathrm{VMS}\left(\{\mathbf{w}^h, q^h\}\right) =&\int_{\Omega}\mathbf{w}^h\cdot\rho\mathbf{f}\,\mathrm{d}\Omega 
+ \int_{\Gamma^\mathrm{N}}\mathbf{w}^h\cdot\mathbf{h}\,\mathrm{d}\Gamma\text{ ,}
\end{align}
where $\mathbf{u}^{\prime}$ and $p^{\prime}$ are the fine-scale terms associated with the VMS formulation, modeled in terms of the residuals of the coarse-scale solution as
\begin{align}\label{eq:u_prime}
&\mathbf{u}^{\prime}=-\tau_{\mathrm{M}}\left(\rho\left(\frac{\partial \mathbf{u}^h}{\partial t}+\mathbf{u}^h \cdot \nabla \mathbf{u}^h-\mathbf{f}\right)-\nabla \cdot \boldsymbol{\sigma}\left(\mathbf{u}^h, p^h\right)\right)\text{ ,} \\
\label{eq:p_prime}
&p^{\prime}=-\rho \tau_{\mathrm{C}} \nabla \cdot \mathbf{u}^h\text{ .}
\end{align}
In the above, $\bar{\tau}$, $\tau_\mathrm{M}$, and $\tau_\mathrm{C}$ are the stabilization parameters, and their detailed expressions can be found in Ref.~\cite{Xu16tetra}. Other options for the stabilization parameters can be found in Refs.~\cite{Tezduyar00Finit, Hsu10Impro, Takizawa18Stabi, Takizawa23Varia}.

In the standard, strongly enforced approach, the trial and test functions are constrained to satisfy $\mathbf{u}^h=\mathbf{g}$ and $\mathbf{w}^h=\mathbf{0}$ on $\Gamma^\mathrm{D}$, and the discrete problem is given by $B^{\mathrm{VMS}} - F^{\mathrm{VMS}} = 0$. To enforce the condition weakly, the constraint is removed from the function spaces and the following weak BC operator is added to the formulation:
\begin{align}
\label{eq:ins-weak-bc}
& - \int_{\Gamma^\mathrm{D}} \mathbf{w}^h\cdot
\left( -p^h\,\mathbf{n} + 2\mu\boldsymbol{\varepsilon}\left(\mathbf{u}^h\right)\mathbf{n} \right)
\,\mathrm{d}\Gamma \nonumber \\
& - \int_{\Gamma^\mathrm{D}}
\left( 2\mu\boldsymbol{\varepsilon}\left(\mathbf{w}^h\right)\mathbf{n} + q^h\,\mathbf{n} \right)
\cdot \left( \mathbf{u}^h - \mathbf{g} \right) \,\mathrm{d}\Gamma \nonumber \\
& - \int_{\Gamma^\mathrm{D}_-} \mathbf{w}^h\cdot
\rho\left( \mathbf{u}^h\cdot\mathbf{n} \right)\left( \mathbf{u}^h - \mathbf{g} \right)
\,\mathrm{d}\Gamma \nonumber \\
& + \int_{\Gamma^\mathrm{D}} \tau_\mathrm{B}^{\mathrm{TAN}}
\left( \mathbf{w}^h - \left(\mathbf{w}^h\cdot\mathbf{n}\right)\mathbf{n} \right) \cdot
\left( \left(\mathbf{u}^h-\mathbf{g}\right)
- \left(\left(\mathbf{u}^h-\mathbf{g}\right)\cdot\mathbf{n}\right)\mathbf{n} \right)
\,\mathrm{d}\Gamma \nonumber \\
& + \int_{\Gamma^\mathrm{D}} \tau_\mathrm{B}^{\mathrm{NOR}}
\left( \mathbf{w}^h\cdot\mathbf{n} \right)
\left( \left(\mathbf{u}^h-\mathbf{g}\right)\cdot\mathbf{n} \right)
\,\mathrm{d}\Gamma \text{ ,}
\end{align}
where $\Gamma^\mathrm{D}_- = \left\{ \mathbf{x} \in \Gamma^\mathrm{D} \mid \mathbf{u}^h\cdot\mathbf{n} < 0 \right\}$ is the inflow part of $\Gamma^\mathrm{D}$, $\mathbf{n}$ is the outward unit normal, and $\tau_\mathrm{B}^{\mathrm{TAN}}$ and $\tau_\mathrm{B}^{\mathrm{NOR}}$ are stabilization parameters acting on the tangential and normal components of the velocity, respectively.

The five terms in Eq.~\eqref{eq:ins-weak-bc} are the vector, momentum-equation counterparts of the contributions derived for the scalar case. The first term is the consistency term. It cancels the traction that arises from integration by parts of the stress in the interior formulation, so that no spurious traction is imposed on the no-slip boundary, and it ensures that the exact solution satisfies the discrete formulation. The second term is the adjoint-consistency term, the symmetric counterpart of the first; its presence renders the formulation adjoint consistent, which is associated with optimal convergence in lower-order norms. Note that the continuity test function $q^h$ also appears in this term, so that the wall-normal component of the Dirichlet condition is weighted by the pressure test function as well. The third term is active only on the inflow portion $\Gamma^\mathrm{D}_-$ and further enforces the boundary condition where flow enters the domain; as shown for pure advection in Section~\ref{sec:advection}, it adds a non-negative contribution to the energy estimate and improves stability without affecting consistency or adjoint consistency. The final two terms are the Nitsche-type stabilization terms, which penalize the deviation of the discrete velocity from $\mathbf{g}$ and provide the coercivity that the consistency and adjoint-consistency terms would otherwise compromise. The tangential and normal components are penalized separately, which allows the no-penetration condition in the normal direction to be enforced more strongly than the no-slip condition in the tangential direction.

Following the analysis of the scalar diffusion equation, the penalty parameter is now given by
\begin{align}
\label{eq:tauB}
\tau_\mathrm{B}^{(\cdot)} = \frac{C_\mathrm{P}\,\mu}{h_n} \text{ ,}
\end{align}
where $h_n$ is a measure of the element size at the boundary, taken in the wall-normal direction, and $C_\mathrm{P}$ is a dimensionless constant whose lower bound follows from the inverse estimate in Eq.~\eqref{eq:inverse-estimate}.

\newpage
\begin{remark}
Under mesh refinement, the no-slip condition is recovered in the limit $h\to 0$~\cite{Bazilevs07Weak2}. On a sufficiently fine boundary-layer mesh, the weak operator yields results close to those obtained with strong enforcement. On a coarse mesh, however, the discrete solution is allowed to slip at the wall. This controlled slip acts as a model for the unresolved thin boundary layer, and accurate global flow quantities are obtained even when the wall-normal mesh size is large~\cite{Bazilevs07Weak1, Bazilevs07Weak2, Hsu12Wind}.
\end{remark}

\begin{remark}
For immersed methods, in which $\Gamma^\mathrm{D}$ cuts arbitrarily through the interiors of the background elements, the integrals in Eq.~\eqref{eq:ins-weak-bc} are evaluated over the immersed surface independently of the background mesh. Because the operator references only traces of the discrete fields on $\Gamma^\mathrm{D}$, no boundary nodes are required. On a cut element, $\tau_\mathrm{B}^{(\cdot)}$ may be prescribed directly~\cite{Xu16tetra, Kamensky15immer} or estimated from the local cut configuration~\cite{Embar10Impos, Annavarapu12robus, Ruess13Weakl}. These properties make weak enforcement applicable in the immersed setting, where the absence of boundary-fitted nodes renders strong enforcement of the boundary conditions infeasible.
\end{remark}

\subsection{Wall-bounded turbulent flows}
\label{sec:turb}

\begin{figure}[!b]
  \centering
  \vspace{-5pt}
  \includegraphics[width=\textwidth]{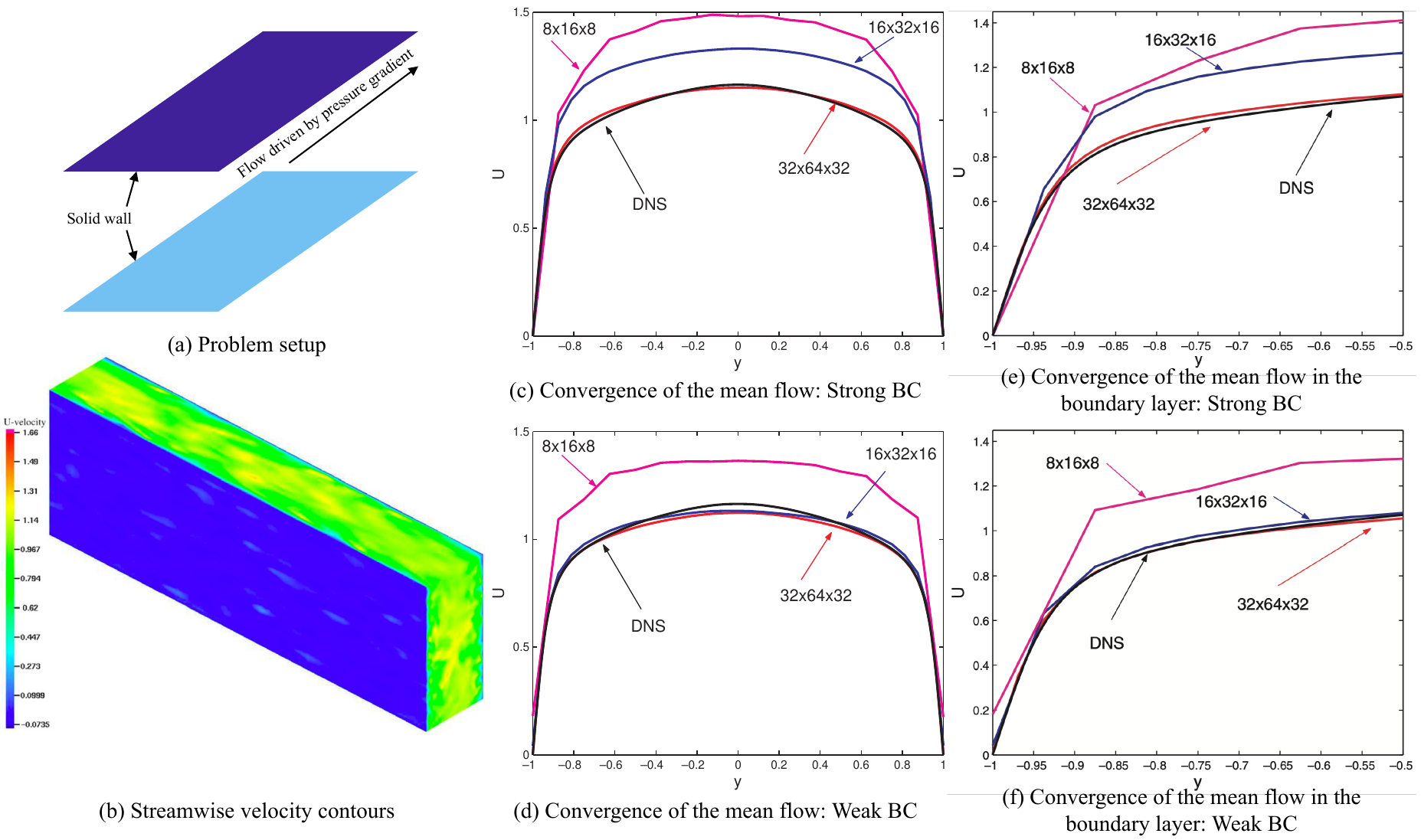}
  \caption{Turbulent channel flow at $Re_\tau = 180$: convergence of the mean velocity profile to the DNS benchmark for strongly and weakly enforced no-slip conditions on uniform meshes. Adapted from Ref.~\cite{Bazilevs07Weak1}.}
  \label{fig:channel180}
  \vspace{-10pt}
\end{figure}

Equilibrium turbulent channel flow, the canonical wall-bounded turbulence benchmark, was the first turbulent flow on which weakly enforced no-slip conditions were shown to outperform their strong counterparts~\cite{Bazilevs07Weak1}. At a friction Reynolds number of $Re_\tau = 180$, on uniform meshes that could not resolve the near-wall region, strongly enforced no-slip conditions overpredicted the mean velocity, while weakly enforced conditions recovered mean profiles in close agreement with the direct numerical simulation (DNS) data~\cite{Kim87Turbu} (Figure~\ref{fig:channel180}). Convergence of the mean profile under mesh refinement was substantially faster for the weak formulation. Even though no turbulence physics was built into the boundary term, its behavior resembled that of a wall-function model.

\begin{figure}[!t]
  \centering
  \includegraphics[width=\textwidth]{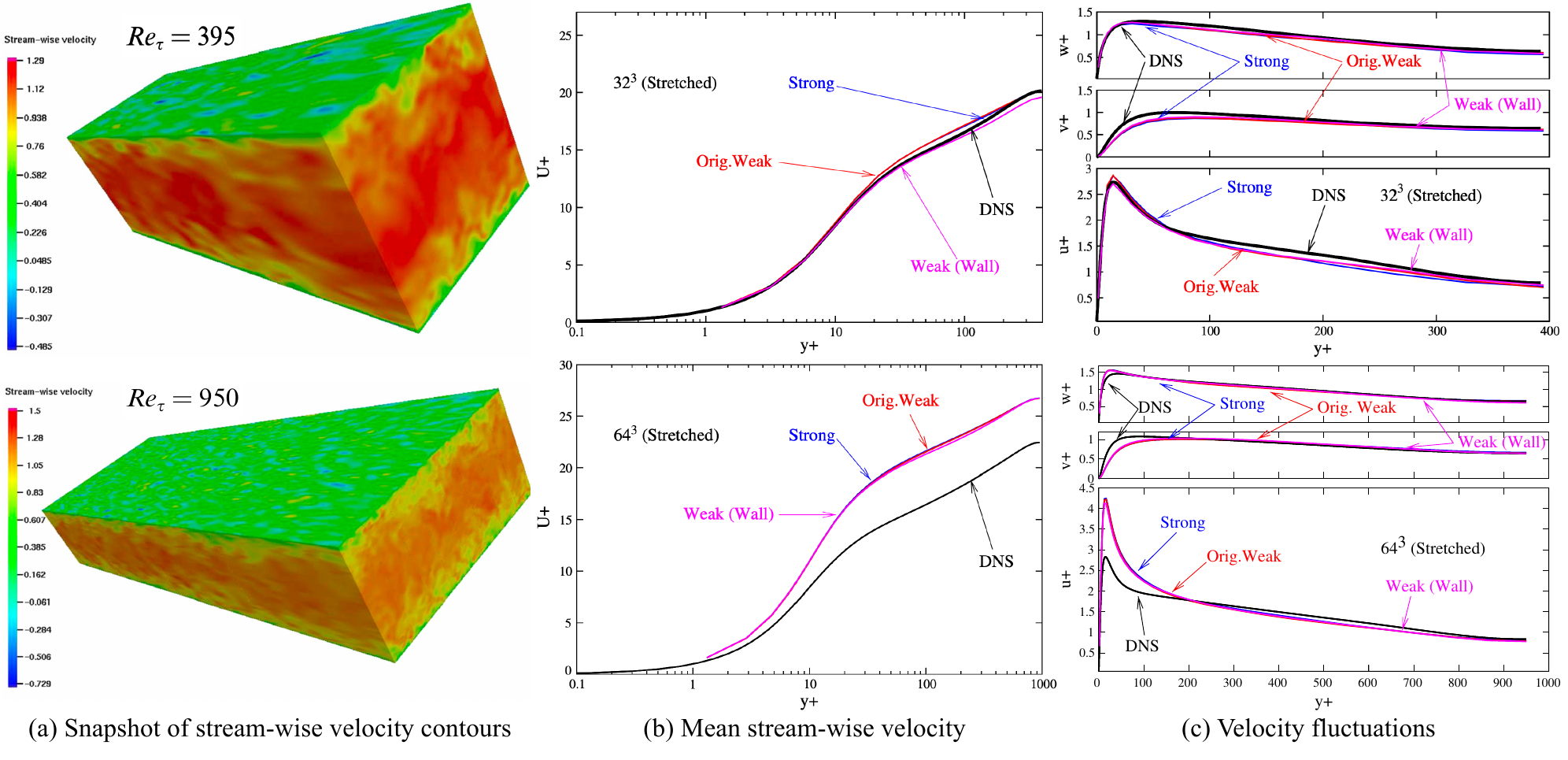}
  \caption{Turbulent channel flow at $Re_\tau = 395$ and $Re_\tau = 950$: mean velocity profiles for the strong, original weak, and law-of-the-wall weak formulations, compared with DNS. Adapted from Ref.~\cite{Bazilevs07Weak2}.}
  \label{fig:turbulent_channel}
\end{figure}

\begin{figure}[!t]
  \centering
  \includegraphics[width=\textwidth]{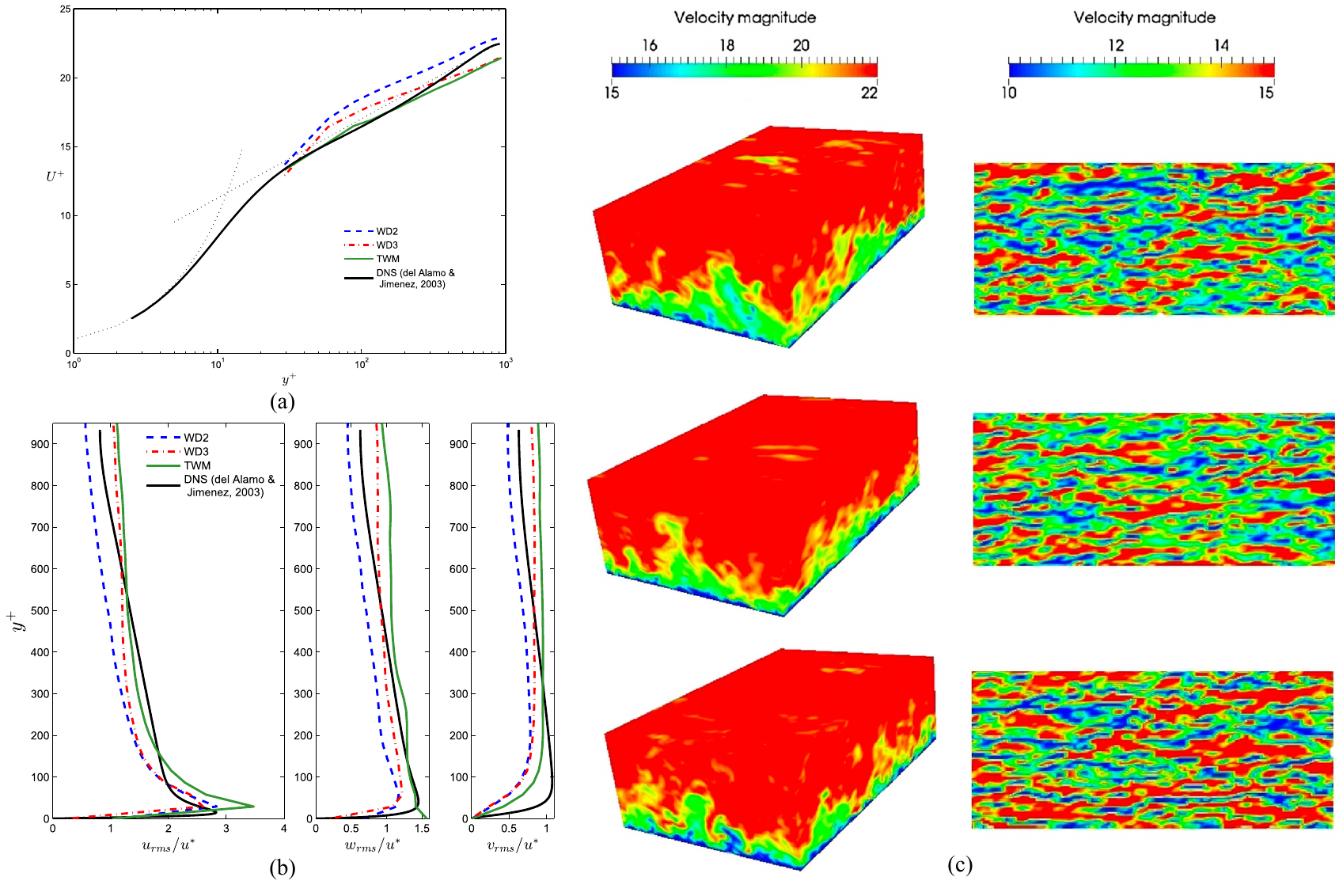}
  \caption{Large-eddy simulation of turbulent channel flow at $Re_\tau = 950$ with near-wall modeling through weakly enforced no-slip conditions: (a)~mean velocity in log wall units, (b)~root-mean-square velocity profiles, both compared with DNS (and the law of the wall in (a)), and (c)~instantaneous flow-speed contours. Adapted from Ref.~\cite{Golshan15Large}.}
  \label{fig:channel-les}
\end{figure}

This observation motivated a reformulation that incorporated the law of the wall consistently into the weak operator, which further improved the prediction of mean-flow quantities at $Re_\tau = 395$ and $950$ (Figure~\ref{fig:turbulent_channel})~\cite{Bazilevs07Weak2}. The isogeometric variational multiscale formulation with weak BCs was subsequently assessed on unstretched meshes~\cite{Bazilevs10Isoge}, and the connection of weak BCs to near-wall modeling was studied further in the context of large-eddy simulation (Figure~\ref{fig:channel-les})~\cite{Golshan15Large}. Because the weak BC formulation reduces to strong enforcement as the wall-normal mesh size vanishes, refinement toward the wall recovers the exact no-slip condition. The additional cost of weak enforcement is negligible, since the added integrals are evaluated only over the wall boundary~\cite{Bazilevs07Weak2}.

\begin{figure}[!t]
  \centering
  \vspace{2pt}
  \includegraphics[width=\textwidth]{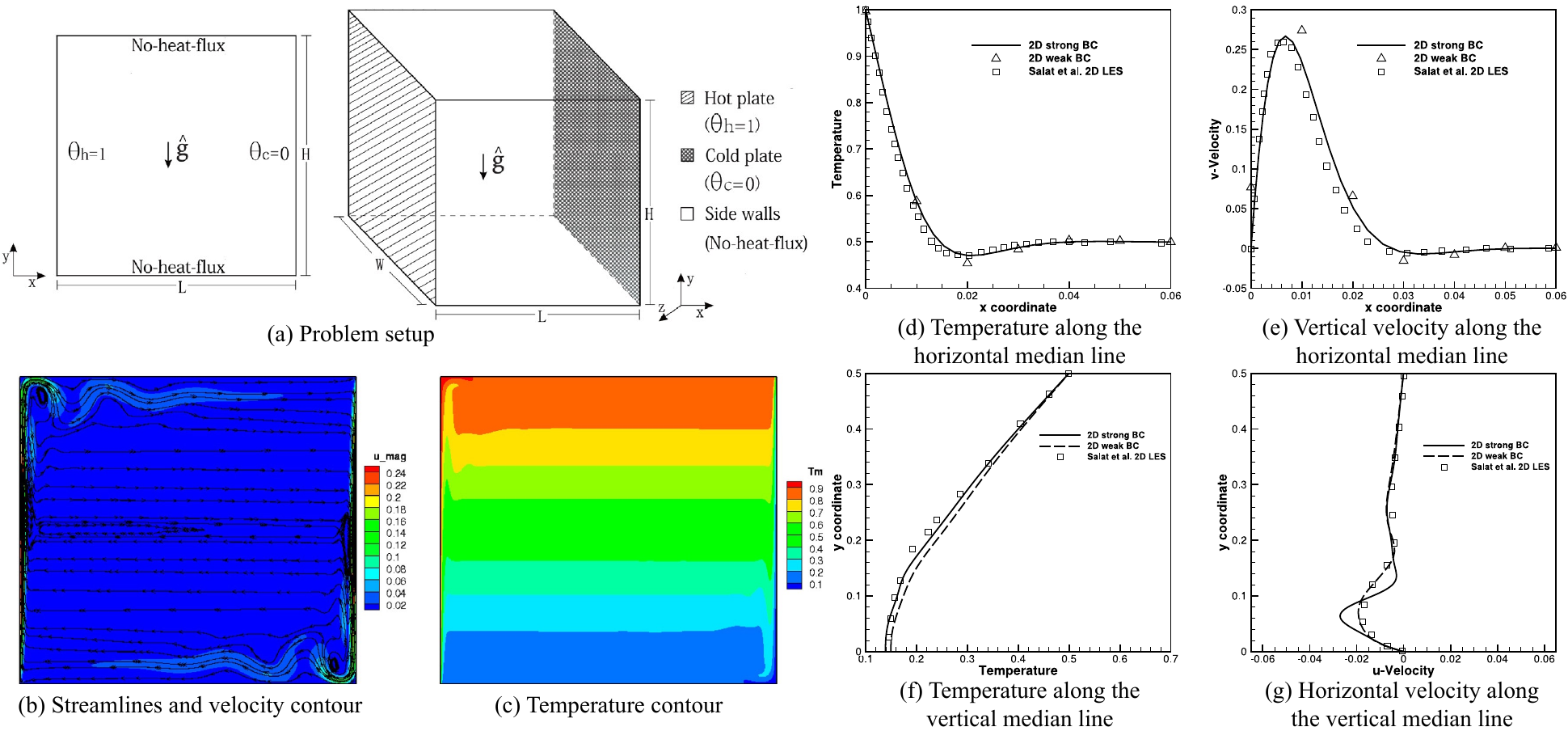}
  \caption{Natural convection in a differentially heated cavity at $Ra = 10^9$: problem setup, flow visualization, and mean temperature and velocity profiles along the median lines for weakly enforced velocity and temperature conditions on a coarse $100\times100$ mesh and strongly enforced conditions on a fine $600\times600$ mesh, compared with reference LES data. Adapted from Ref.~\cite{Xu19resid}.}
  \label{fig:cavity}
  \vspace{-15pt}
\end{figure}

The same near-wall mechanism extends to thermal boundary layers. For natural convection in a differentially heated cavity at high Rayleigh number, where the thermal and velocity boundary layers on the hot and cold walls are both thin, weakly imposing the velocity and temperature Dirichlet conditions on coarse, boundary-fitted meshes improved the prediction of the wall Nusselt number relative to strong enforcement~\cite{Xu19resid}. On a coarse, uniform mesh, weak enforcement also recovered temperature and velocity profiles that matched both the reference data and the strongly enforced results on fine meshes (Figure~\ref{fig:cavity}). This establishes that the coarse-mesh advantage observed for the no-slip condition carries over to thermal Dirichlet data, a property that becomes important in the compressible-flow setting of Section~\ref{sec:comp}, where the wall-temperature condition is also imposed weakly.

\subsection{Aerodynamics of wind turbines}
\label{sec:windturbine}

\begin{figure}[!t]
  \centering
  \includegraphics[width=\textwidth]{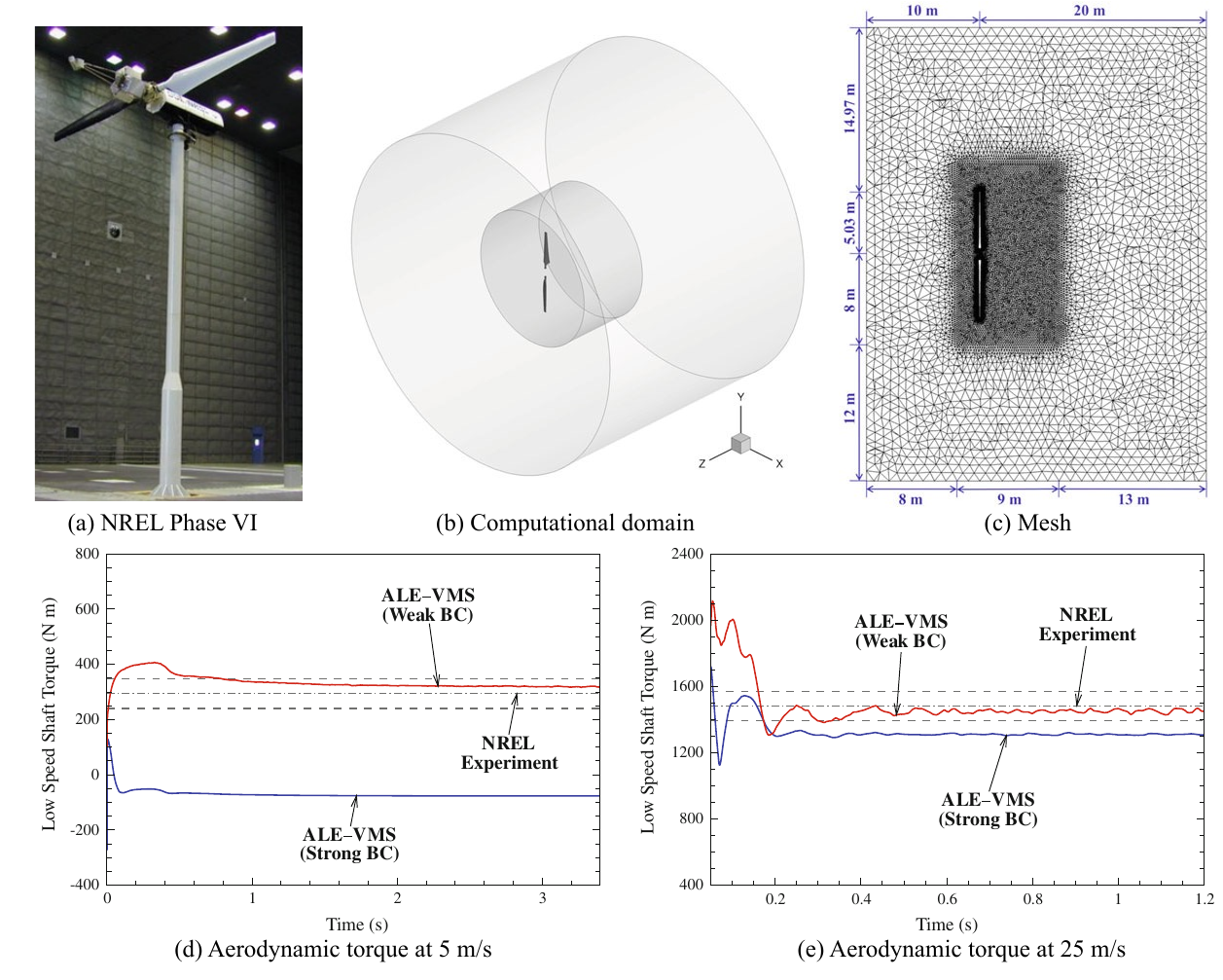}
  \caption{NREL Phase~VI rotor: aerodynamic (low-speed-shaft) torque for weakly and strongly enforced no-slip conditions, compared with experimental data, at 5~m/s and 25~m/s. The strong-BC computation stalls and underpredicts the torque on the coarse boundary-layer mesh. Adapted from Ref.~\cite{Hsu12Wind}.}
  \label{fig:nrelVI-torque}
  \vspace{-2pt}
\end{figure}

Wind-turbine aerodynamics involves thin boundary layers developing over large-scale geometries. Fully resolving the boundary layer over the rotor is often computationally expensive and impractical. Obtaining accurate solutions on coarse meshes is therefore highly desirable, and this is the regime in which weak BCs are most effective. Combined with the ALE--VMS formulation on a rotating mesh, weak BCs reproduced the measured aerodynamic (low-speed-shaft) torque and the blade-surface pressure distribution of the full-scale NREL Phase~VI rotor on meshes with only coarse boundary-layer refinement~\cite{Hsu12Wind}. On the same meshes, strong enforcement produced spurious flow separation and underpredicted the torque (Figure~\ref{fig:nrelVI-torque}).

\begin{figure}[!t]
  \centering
  \includegraphics[width=\textwidth]{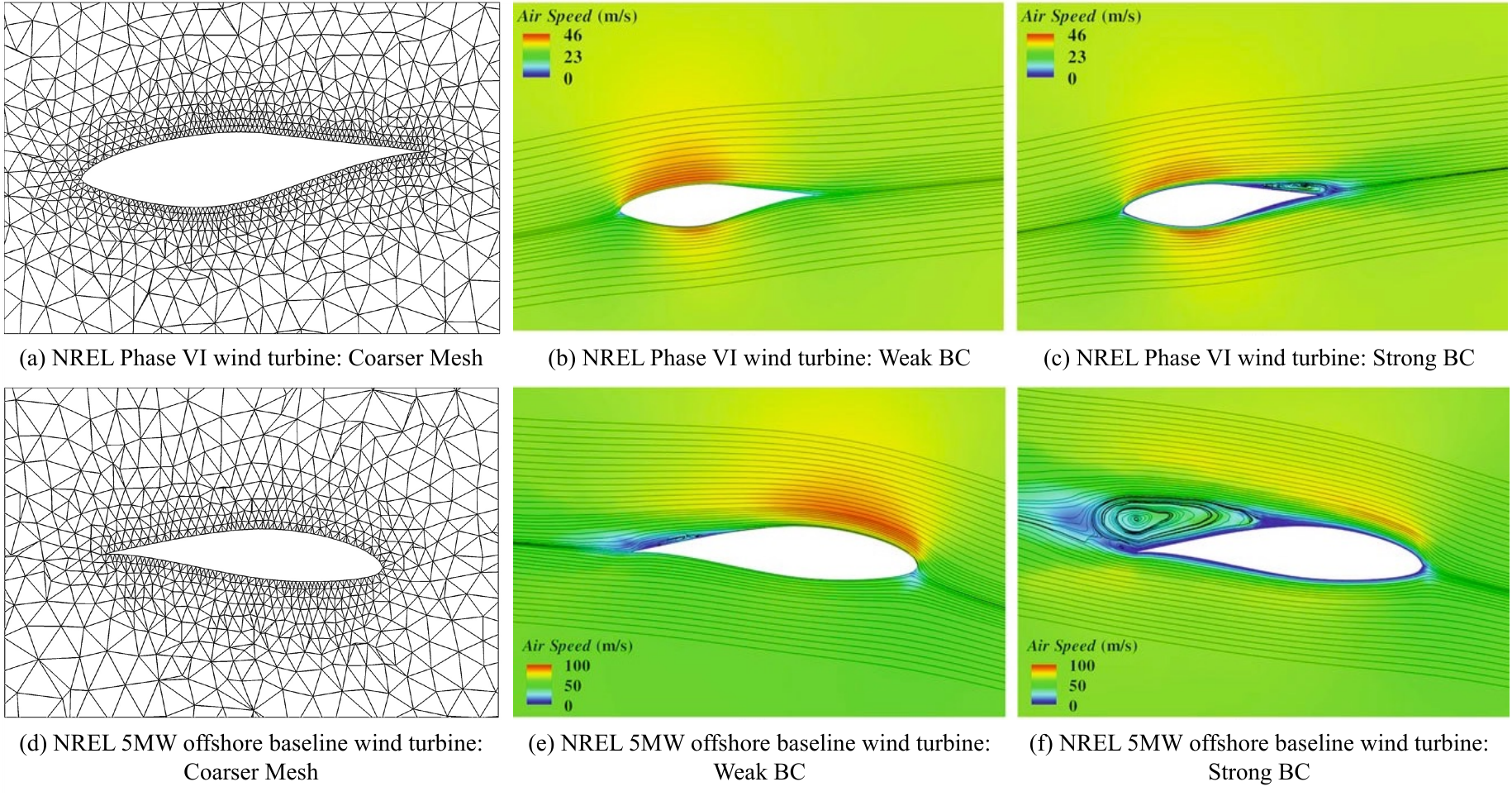}
  \caption{Weakly versus strongly enforced no-slip conditions on coarse boundary-layer meshes, shown through air-speed contours with velocity streamlines at a blade cross-section: (a)--(c)~NREL Phase~VI rotor at 5~m/s ($0.8R$) and (d)--(f)~NREL 5MW offshore baseline rotor ($0.75R$), each with weak and strong enforcement. Strong enforcement produces a thick, under-resolved boundary layer and premature separation at both scales. Adapted from Ref.~\cite{Hsu12Wind}.}
  \label{fig:nrelVI-contours}
\end{figure}

The cause of the torque error is visible in the near-wall flow. On the coarse mesh, the strongly enforced computation develops a thick, under-resolved boundary layer that separates prematurely, whereas the weakly enforced computation allows the flow to remain attached and recovers the correct surface pressure. The same comparison, with no change to the formulation, was repeated on the much larger NREL 5MW offshore wind-turbine rotor~\cite{Hsu12Wind}. The behavior was similar at both scales: weak enforcement on a coarse boundary-layer mesh reproduced the reference torque, while strong enforcement produced a thick spurious boundary layer, unphysical separation, and an inaccurate torque. The cross-section contour and streamline plots (Figure~\ref{fig:nrelVI-contours}) clearly illustrate this difference in flow separation for both rotors.

\begin{figure}[!t]
  \centering
  \includegraphics[width=\textwidth]{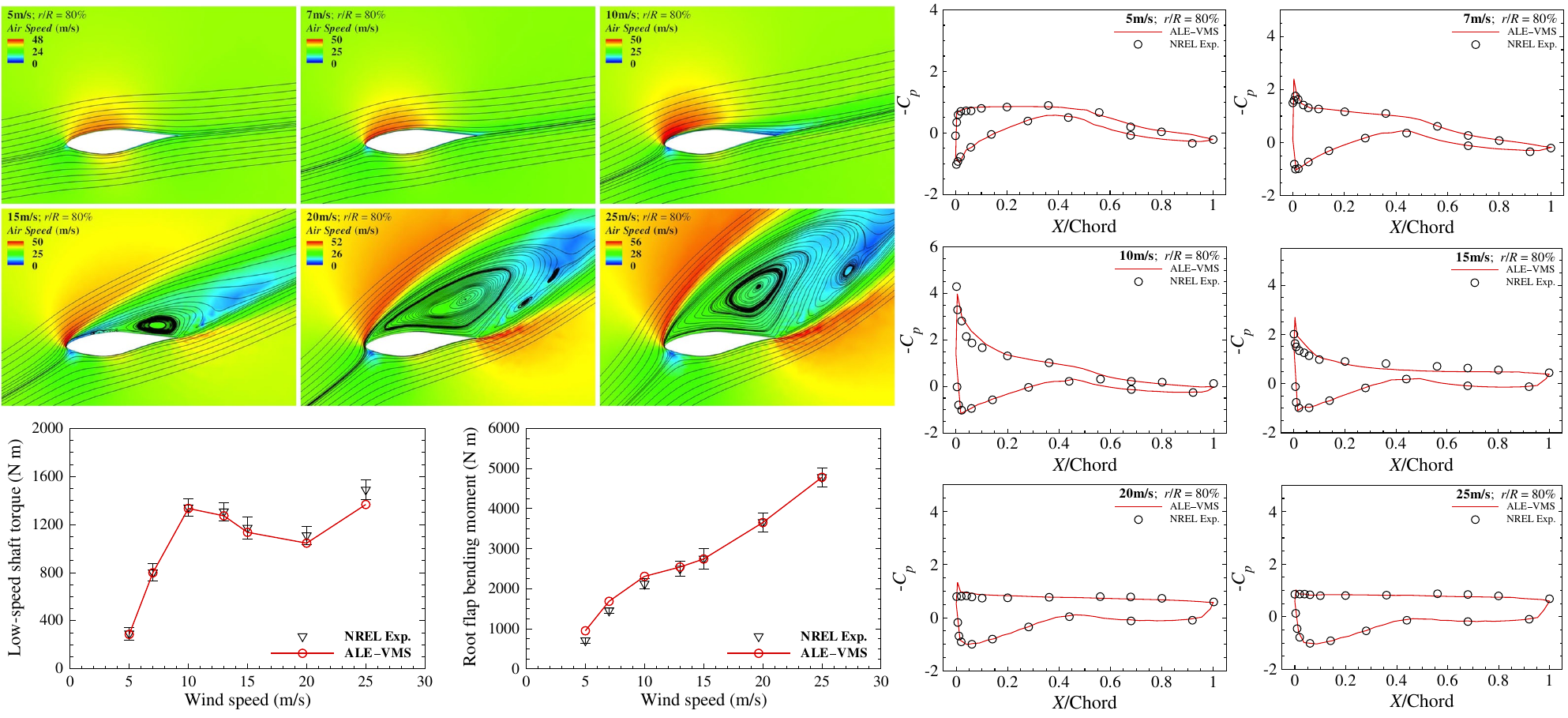}
  \caption{NREL Phase~VI rotor-only validation with weakly enforced no-slip conditions, across wind speeds from 5 to 25~m/s: air-speed contours and velocity streamlines, low-speed-shaft torque and root flap bending moment versus wind speed, and sectional pressure coefficient. The latter two are compared with the NREL experimental data. Adapted from Ref.~\cite{Hsu14Finit}.}
  \label{fig:nrelVI-rotor}
\end{figure}

A comprehensive validation of the weakly enforced ALE--VMS formulation against the NREL Phase~VI experiment was subsequently carried out across the full measured wind-speed range, from 5 to 25~m/s~\cite{Hsu14Finit}. The rotor-only simulations reproduced the measured low-speed-shaft torque, root flap bending moment, sectional force coefficients, and blade-surface pressure distributions, again on boundary-layer meshes far coarser than strong enforcement would require. The same weak BC operator and the same parameter values were used throughout, even though the flow regime varied substantially with wind speed, ranging from fully attached at 5~m/s to massively separated at 25~m/s (Figure~\ref{fig:nrelVI-rotor}). This insensitivity to the flow regime contrasts with classical eddy-viscosity wall models, which typically require recalibration for each flow condition.

\begin{figure}[!t]
  \centering
  \includegraphics[width=\textwidth]{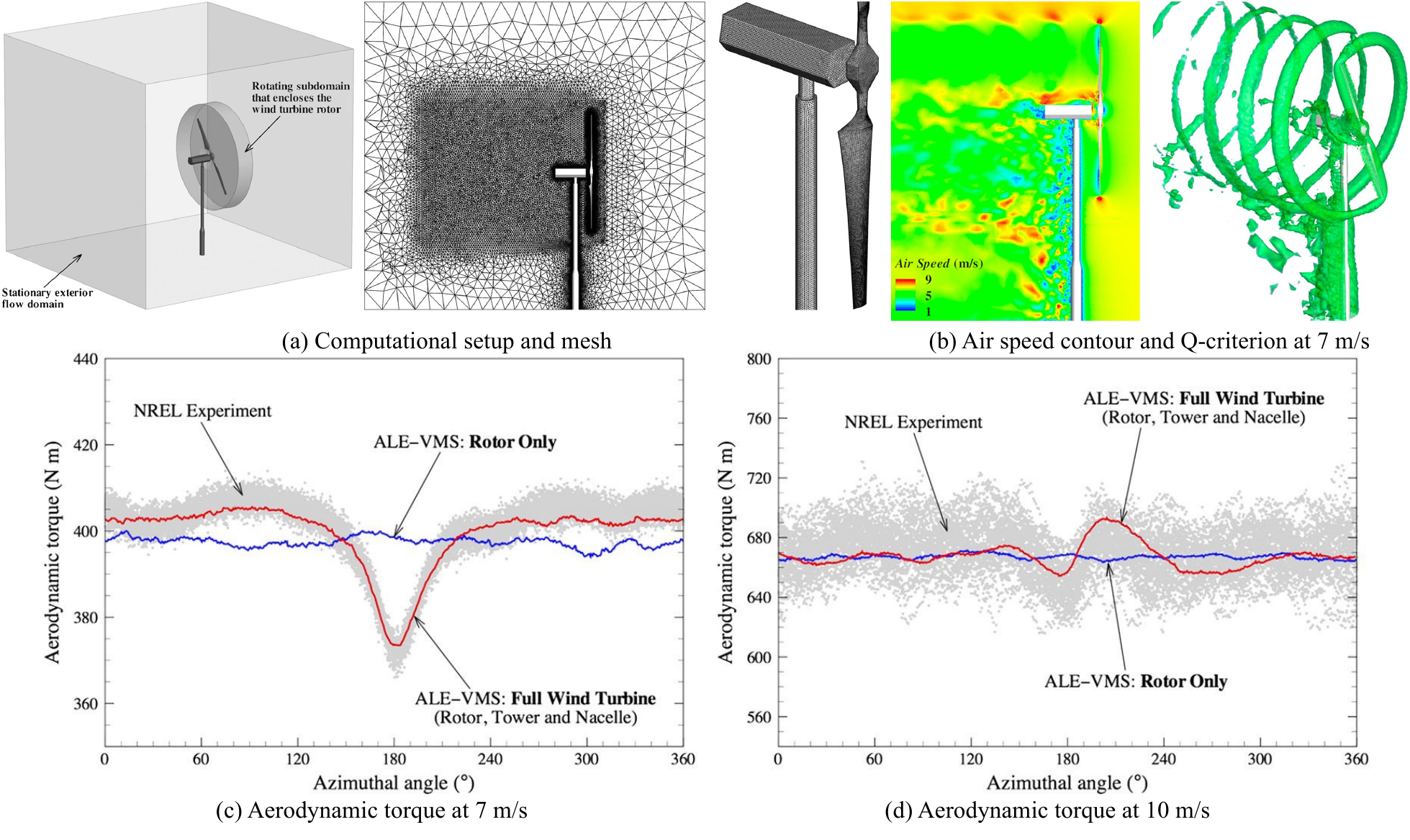}
  \caption{NREL Phase~VI full-wind-turbine simulation with the sliding-interface method at 7~m/s and 10~m/s: the full-tower model closely matches the experimental measurements of the torque drop as the blade passes in front of the tower ($180^\circ$), in contrast to the rotor-only result. Adapted from Ref.~\cite{Hsu14Finit}.}
  \label{fig:nrelVI-fullturbine}
  \vspace{-11pt}
\end{figure}

To capture the influence of the tower, the study was extended to the full wind-turbine configuration of rotor, nacelle, and tower~\cite{Hsu14Finit}. The rotating (rotor) and stationary (tower) subdomains were coupled across a non-matching sliding interface, with kinematic and traction continuity enforced weakly through Nitsche-type interface terms analogous to the weak wall BC operators~\cite{Bazilevs08NURBS, Hsu14Finit}. The full-machine results reproduced the measured blade--tower interaction: as a blade passes in front of the tower, the single-blade aerodynamic torque at the 7~m/s condition drops by about 8\%, an effect the rotor-only model cannot capture (Figure~\ref{fig:nrelVI-fullturbine}).

\subsection{Immersogeometric analysis}
\label{sec:imga}

The previous applications demonstrate how weak boundary conditions reduce computational cost in boundary-fitted frameworks. Here, we show that these operators are essential in non-boundary-fitted discretizations. As discussed in Section~\ref{sec:intro_immersed}, immersed methods embed the geometry within a background mesh that does not conform to the physical boundaries. The interface arbitrarily intersects the fluid elements, making the strong imposition of wall conditions infeasible. Immersogeometric analysis instead couples a weak boundary condition operator with geometrically accurate quadrature over the intersected elements, ensuring both the weak imposition of Dirichlet boundary conditions and a faithful representation of the immersed surface without a boundary-fitted mesh.

\subsubsection{The tetrahedral finite cell method}

The immersogeometric concept was established for incompressible flow around complex geometries using a tetrahedral finite cell method, in which the immersed object is embedded in a non-boundary-fitted tetrahedral mesh and the cut elements are integrated by sub-cell-based adaptive quadrature~\cite{Xu16tetra}. Tetrahedral background elements were chosen for their flexibility in providing local refinement around complex geometries, and the intersected elements were recursively subdivided into quadrature sub-cells to faithfully capture the volume integral on the fluid side of the immersed surface without modifying the underlying mesh. 

\begin{figure}[!b]
  \centering
  \includegraphics[width=\textwidth]{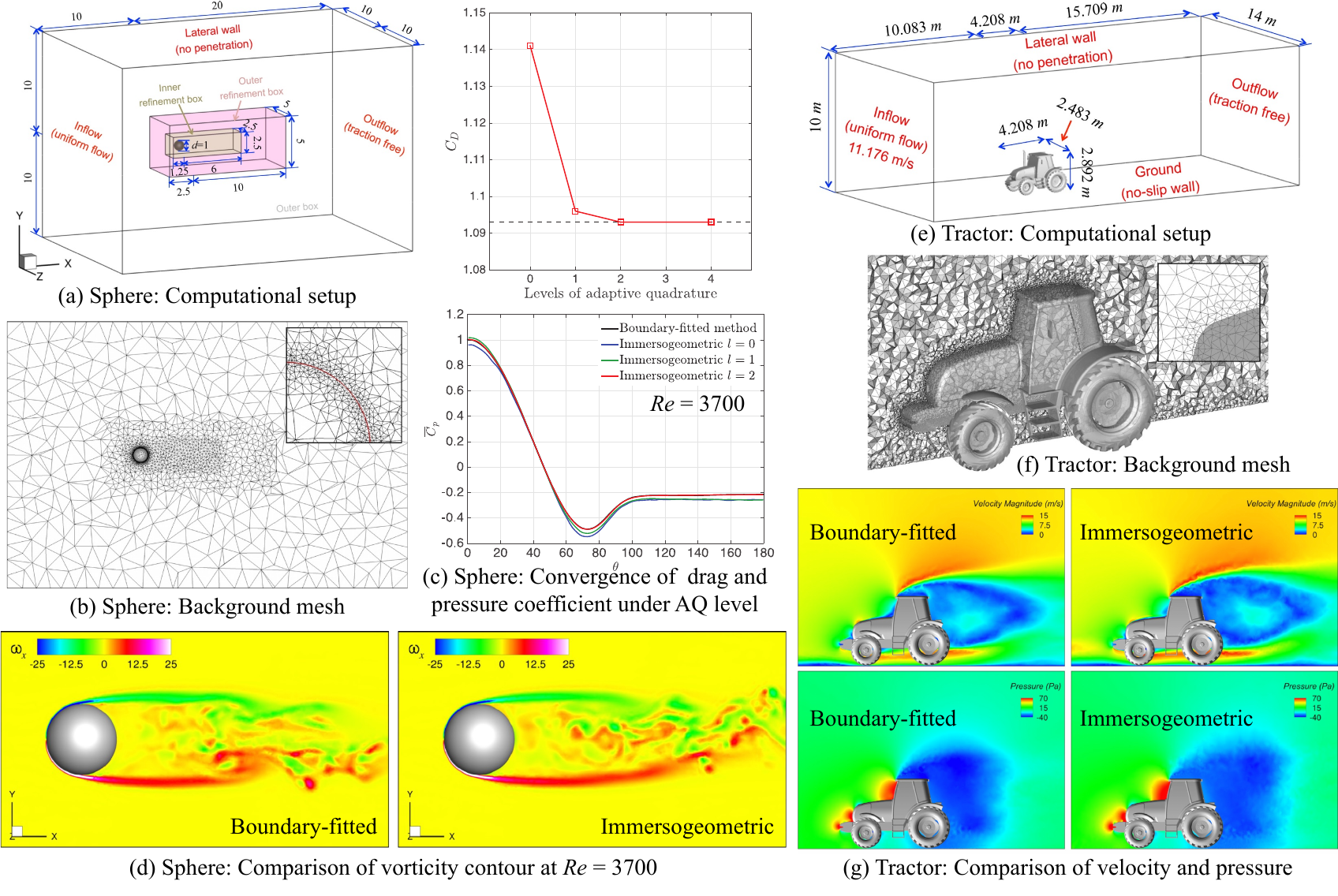}
  \caption{Tetrahedral finite cell method: immersed sphere benchmark and industrial-scale tractor aerodynamics on non-boundary-fitted meshes. Adapted from Ref.~\cite{Xu16tetra}.}
  \label{fig:tetfcm}
  \vspace{-10pt}
\end{figure}

The method was assessed on flow past a sphere at Reynolds numbers of $Re = 100$, $300$, and $3700$. This range covers steady and unsteady laminar regimes, as well as the transitional and turbulent wake features present at $Re=3700$. The drag coefficient, Strouhal number, recirculation length, and surface pressure distribution were all in excellent agreement with boundary-fitted reference values on meshes of comparable resolution. At $Re=3700$, the computed time-averaged quantities and wake statistics also compared favorably to DNS data from the literature (Figure~\ref{fig:tetfcm}). A systematic study showed that the immersogeometric results converge to the boundary-fitted reference under subdivision of the adaptive quadrature, without requiring the background fluid mesh to be refined beyond the resolution of its boundary-fitted counterpart. The method was then applied to the aerodynamic analysis of a full-scale agricultural tractor at $Re \approx 3\times 10^6$. The tractor contains fine geometric features that would typically require extensive defeaturing for boundary-fitted analysis. However, the immersogeometric workflow efficiently generated a background mesh that intersects with the tractor surface without requiring geometric simplification. The time-averaged pressure distribution and surface drag were captured with an accuracy comparable to that of a boundary-fitted reference on a mesh of similar size (Figure~\ref{fig:tetfcm}).

\subsubsection{Direct flow analysis on CAD models}

\begin{figure}[!b]
  \centering
  \vspace{-20pt}
  \includegraphics[width=.89\textwidth]{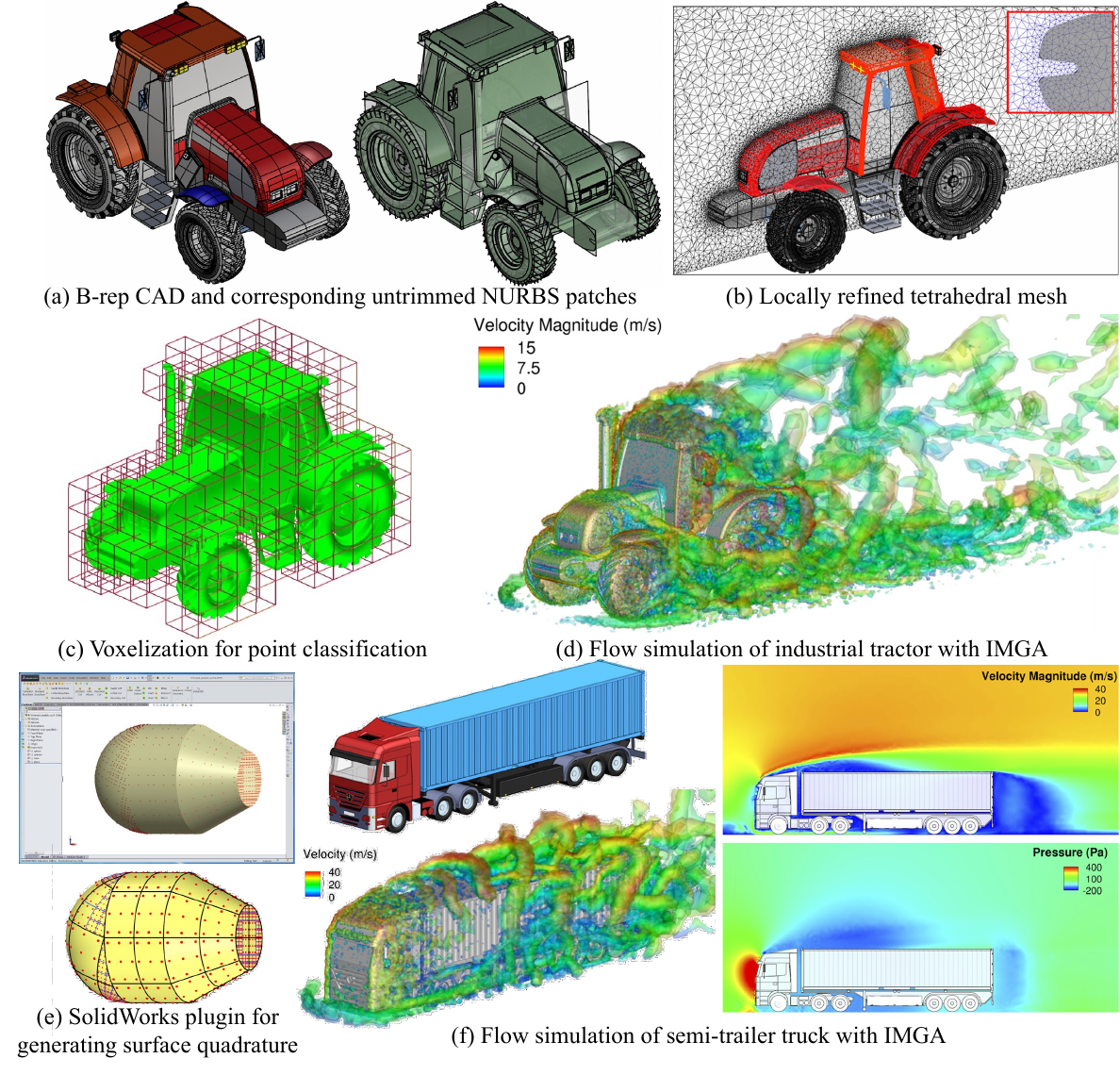}
  \caption{Immersogeometric analysis driven directly by B-rep CAD models. Adapted from Refs.~\cite{Hsu16Direc, Wang17Rapid}.}
  \label{fig:brep}
  \vspace{-20pt}
\end{figure}

Because IMGA enforces boundary conditions directly on the immersed geometry through the weak BC operator rather than on a discrete surface mesh, it can natively utilize B-rep CAD models. Two geometric operations drive the analysis: a surface quadrature on the immersed B-rep for evaluating the weak BC integrals, and a point membership classification of background-mesh quadrature points to identify those inside or outside the immersed object~\cite{Hsu16Direc}. Surface quadrature is generated directly on the trimmed B-rep surfaces, and the point membership classification is performed efficiently through a high-resolution voxelization obtained from GPU rendering of the immersed surfaces. The framework was initially developed for trimmed NURBS patches~\cite{Hsu16Direc} and subsequently extended to analytic surfaces such as planes, cones, spheres, and tori that would otherwise require conversion to NURBS for analysis~\cite{Wang17Rapid}. Validation on flow past a sphere using three boundary representations (triangulated, untrimmed NURBS, and trimmed NURBS) showed that the drag coefficient converges to values in close agreement with reference data under refinement of the surface quadrature density. The framework was then applied to the aerodynamic analysis of an agricultural tractor and a semi-trailer truck, both immersed directly as B-rep CAD models into a background mesh and bypassing the labor-intensive steps of defeaturing and surface meshing (Figure~\ref{fig:brep}).

\subsubsection{Free-surface, moving-object, and industrial-scale flows}

\begin{figure}[!b]
  \centering
  \includegraphics[width=\textwidth]{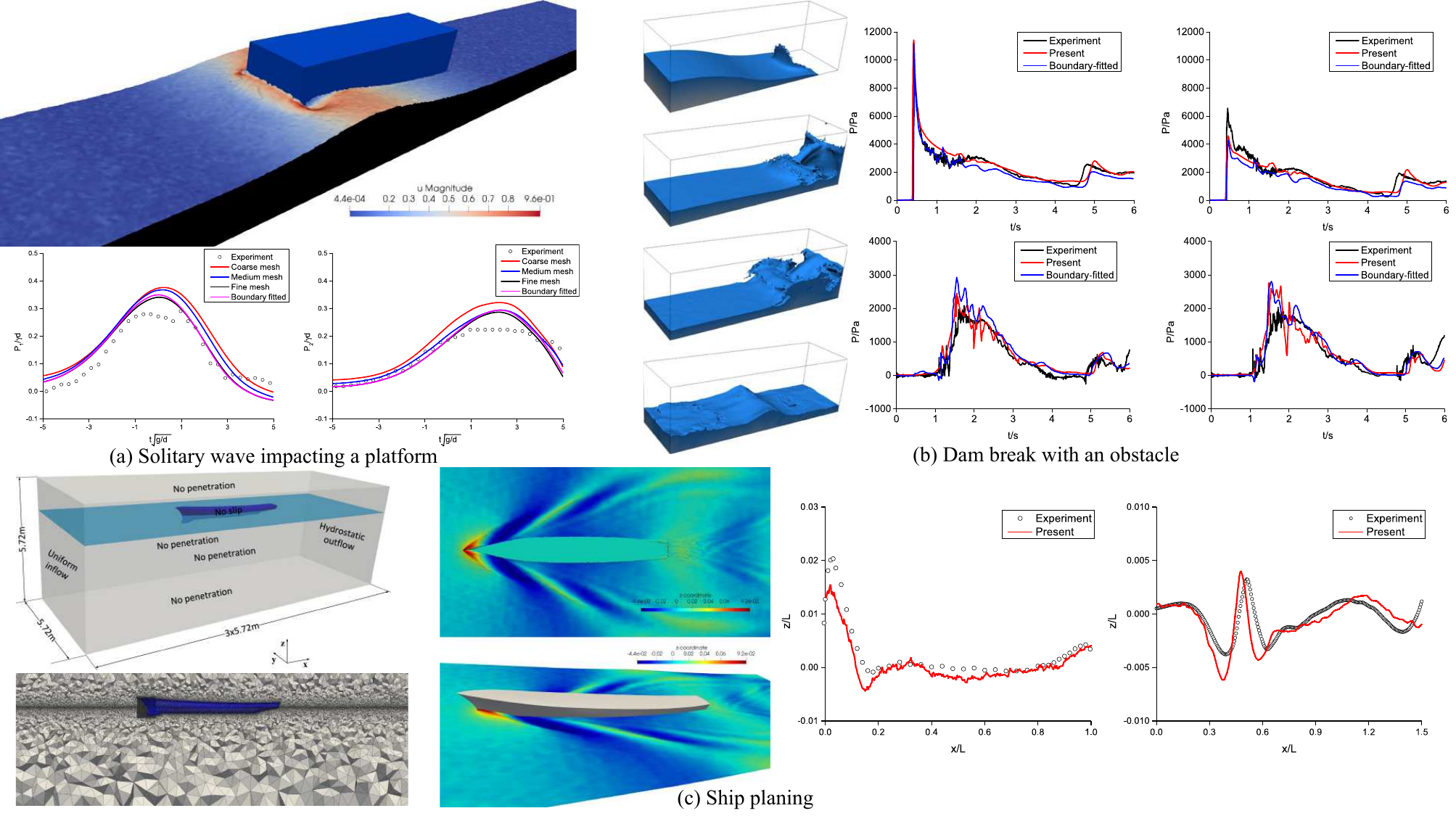}
  \caption{Immersogeometric analysis of free-surface flows: solitary wave impacting a platform, dam break with an obstacle, and ship planing. Adapted from Ref.~\cite{Zhu20immer}.}
  \label{fig:imga-free}
\end{figure}

The immersogeometric framework was extended to free-surface flows of interest in marine engineering~\cite{Zhu20immer}. The weak BC operator on the immersed structure was coupled with a level-set treatment of the air--water interface, and a residual-based VMS formulation was used to stabilize the coupled Navier--Stokes and level-set convection equations. The method was assessed on three benchmark problems of increasing complexity: a solitary wave impacting a stationary platform, a dam break with an obstacle, and the planing of a DTMB 5415 ship hull. The predicted pressure histories and free-surface elevations were in good agreement with experimental measurements and boundary-fitted reference simulations (Figure~\ref{fig:imga-free}).

\begin{figure}[!t]
  \centering
  \includegraphics[width=\textwidth]{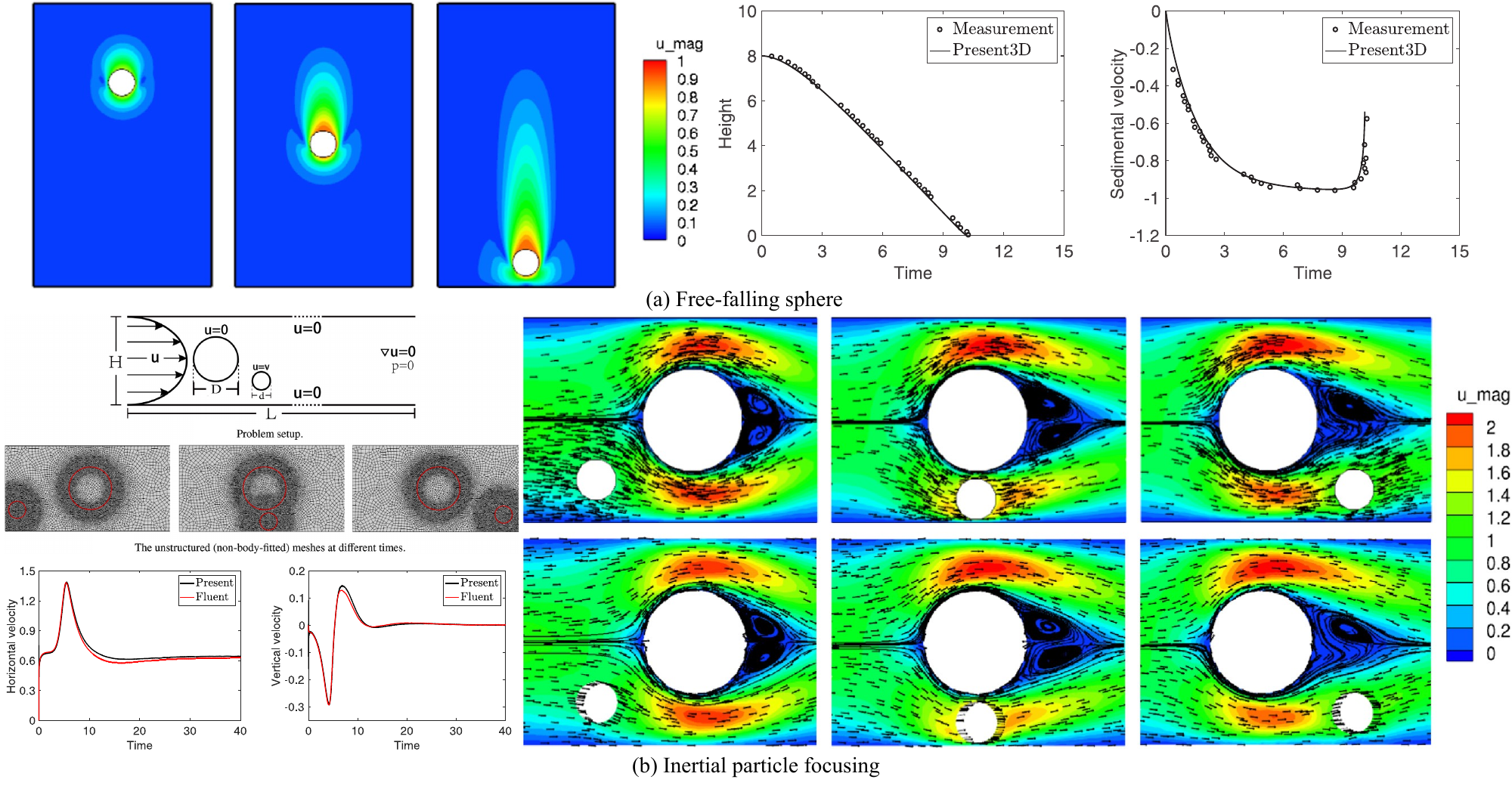}
  \caption{Immersogeometric analysis of moving objects in incompressible flows: free-falling sphere in a viscous fluid and inertial particle focusing in a pillar-obstructed microfluidic channel, with the no-slip condition enforced weakly on the moving immersed surface. Adapted from Ref.~\cite{Xu19Immer2}.}
  \label{fig:imga-moving}
\end{figure}

IMGA was also extended to moving-object problems, in which a rigid body translates and rotates through a non-boundary-fitted background mesh under the action of integrated surface forces and external body forces~\cite{Xu19Immer2}. Particular attention was given to the freshly cleared nodes, namely background nodes that were inside the object at one time step and in the fluid domain at the next. Because these nodes lacked a field history, their values were interpolated from neighboring fluid nodes and points on the object surface. The method was validated on free-falling 2D cylinders and 3D spheres in viscous fluids, with terminal velocities and trajectories matching analytical and experimental references. It was then applied to 2D inertial particle focusing in straight and pillar-obstructed microfluidic channels (Figure~\ref{fig:imga-moving}). At comparable accuracy, the moving IMGA simulations were approximately eight times faster than a boundary-fitted commercial code at the same mesh density. The framework was later extended to utilize adaptively refined, octree-based background meshes to enable scalable 3D simulations of inertial particle migration in straight and converging-diverging microchannels~\cite{Xu21octre}.

\begin{figure}[!t]
  \centering
  \includegraphics[width=\textwidth]{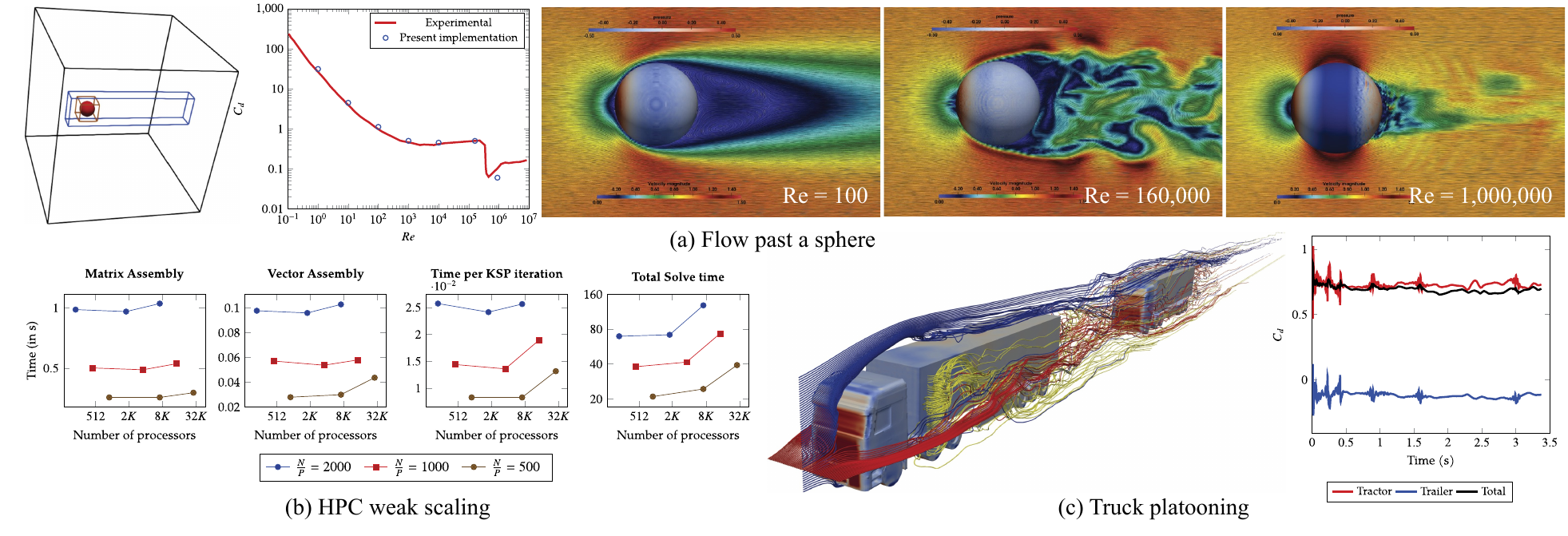}
  \caption{Immersogeometric analysis of industrial-scale flows: sphere drag crisis and semi-truck platooning. Adapted from Ref.~\cite{Saurabh21Indus}.}
  \label{fig:imga-ind}
\end{figure}

To support industrial-scale large-eddy simulations, IMGA was implemented on adaptively refined octree background meshes with a parallel framework that demonstrated scalability up to 32K processors~\cite{Saurabh21Indus}. Efficiency at this scale was further enabled by rapid point membership classification and tensorized matrix assembly. The framework was applied to the flow over a sphere across Reynolds numbers from 1 to $10^6$, reproducing the drag crisis observed in experiments. It was also used to characterize the drag reduction of the trailing vehicle in a two-truck semi-trailer platoon (Figure~\ref{fig:imga-ind}).

\subsubsection{Direct flow analysis on point clouds, photogrammetric reconstructions, and medical images}

\begin{figure}[!t]
  \centering
  \includegraphics[width=\textwidth]{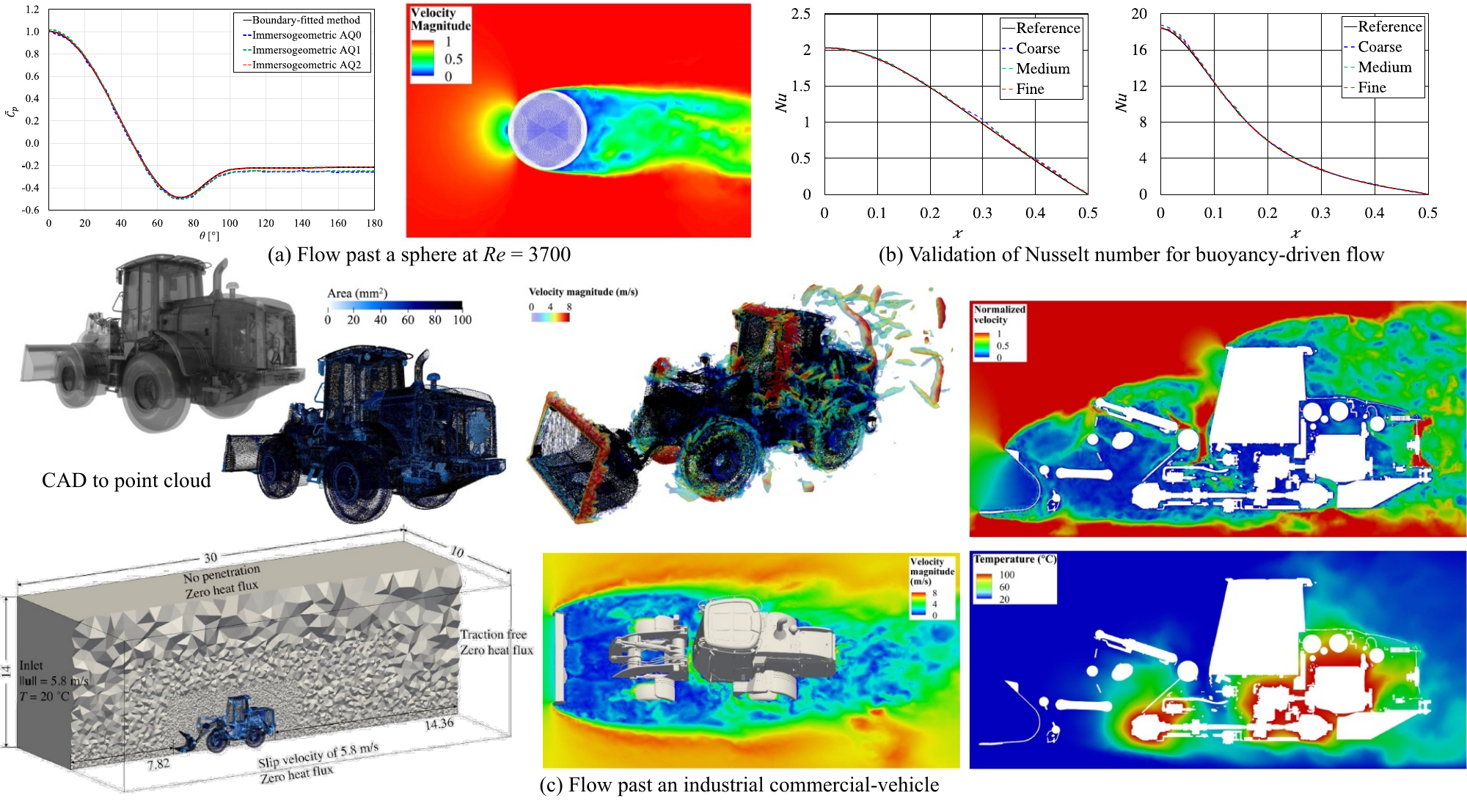}
  \caption{Direct immersogeometric flow and heat transfer analysis on point clouds, with the no-slip condition enforced weakly on the points: fluid and thermal validations against boundary-fitted references and an industrial construction-vehicle application. Adapted from Ref.~\cite{Balu23Direc}.}
  \label{fig:pointcloud-direct}
\end{figure}

Alongside direct point-cloud analysis in solid mechanics with the finite cell method~\cite{Kudela20Direc, Hartmann22Enfor}, IMGA has recently been applied directly to sampled representations of physical objects, bypassing the geometric model altogether. A sufficiently dense point cloud of the object boundary provides the geometric traces needed to evaluate the weak BC operator, enabling flow analysis without reconstructing a CAD model or a discrete surface mesh~\cite{Balu23Direc, Wang23Photo, Jaiswal24Mesh, Corpuz25Direc}. Three geometric operations drive the analysis directly from the point cloud: inside--outside classification of the background-mesh quadrature points, estimation of the surface normal at each point, and computation of the Jacobian determinant used for surface integration. The classification uses the winding number, which extends naturally to non-manifold and non-watertight inputs and therefore does not require a watertight surface representation. The framework was assessed on benchmark problems across a wide range of Reynolds numbers, including incompressible flow past a sphere, with surface pressure and flow quantities in excellent agreement with boundary-fitted references. It was further extended to heat transfer, with thermal Dirichlet conditions imposed weakly on the points alongside the no-slip condition. The method was then applied to a John Deere 544K Wheel Loader, sampled from the as-designed CAD assembly as a point cloud of over 12 million points (Figure~\ref{fig:pointcloud-direct}). The wheel loader is composed largely of finite-thickness shells with gaps, intersections, and collocated faces that would require extensive cleanup and defeaturing for boundary-fitted analysis. Using the point-cloud representation circumvents this geometry-simplification step entirely. The small geometric details are fully preserved, and their impact on the flow is clearly visible (Figure~\ref{fig:pointcloud-direct}).

\begin{figure}[!t]
  \centering
  \includegraphics[width=\textwidth]{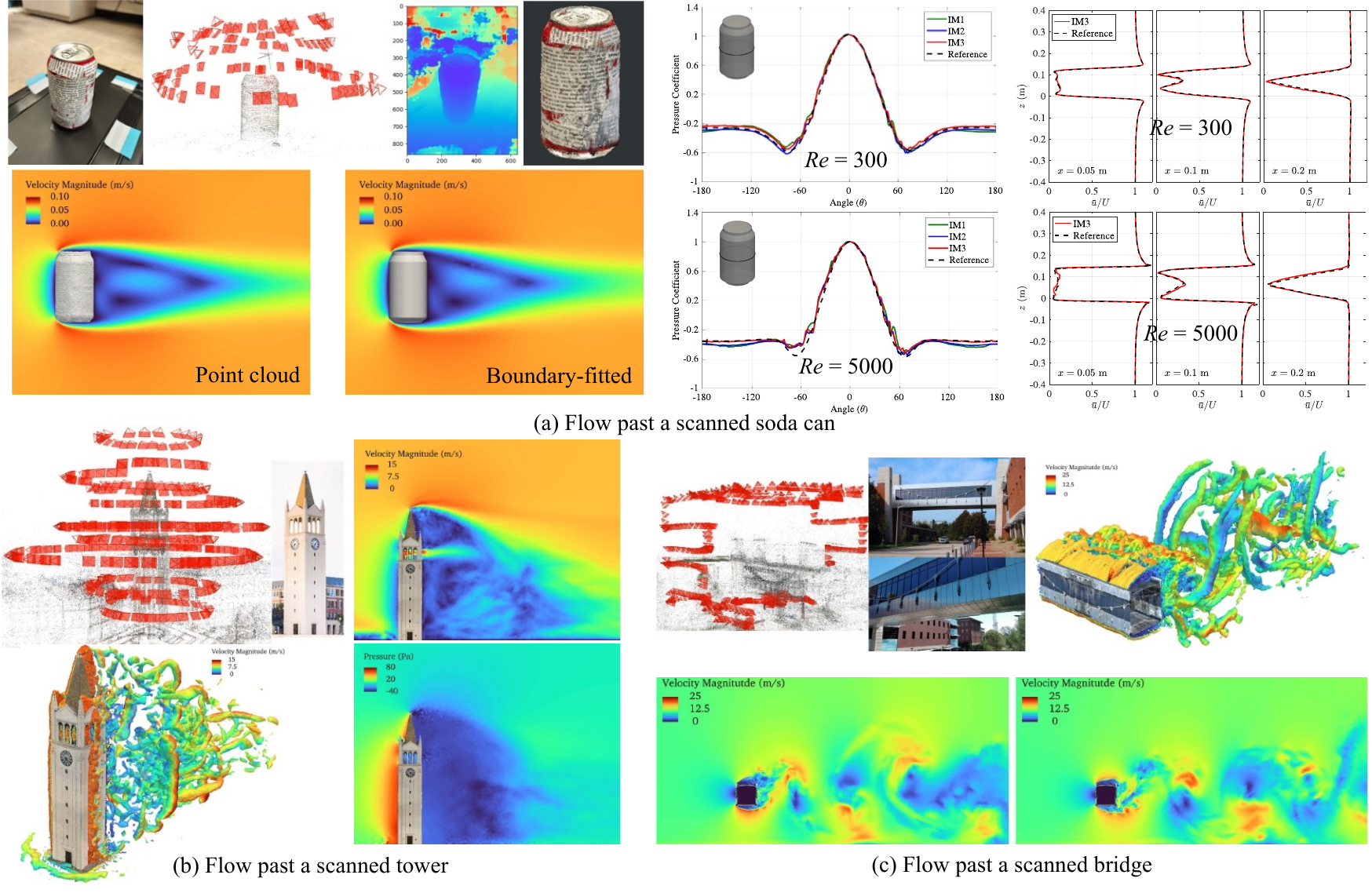}
  \caption{Photogrammetry-based CFD on point clouds of in-use civil structures, with the no-slip condition enforced weakly on the reconstructed points. Adapted from Ref.~\cite{Wang23Photo}.}
  \label{fig:photogrammetry}
\end{figure}

Photogrammetry-based CFD~\cite{Wang23Photo} extends this capability to in-use civil structures using geometry acquired from photographic surveys. The reconstruction pipeline begins with a structure-from-motion (SfM) step to estimate camera poses from a set of overlapping 2D images. Next, a learned multi-view stereo step infers a dense depth map for each image and fuses these maps into a dense point cloud. An octree-based subsampling method then generates a more uniformly distributed point cloud suitable for immersogeometric analysis. The framework was validated on the flow past a standard 12~oz soda can reconstructed from cellphone imagery, at Reynolds numbers of 300 and 5000, with the surface pressure coefficients and the wake velocity profiles matching the boundary-fitted references computed on an idealized CAD model of the same can~\cite{Wang23Photo}. The approach was then applied to two in-use civil structures captured by drone surveys: a bell tower and a pedestrian bridge (Figure~\ref{fig:photogrammetry}). In both cases, the flow analysis was carried out directly on the reconstructed point clouds. These applications demonstrated the robustness of the photogrammetry-based high-fidelity CFD framework for evaluating the aerodynamics of complex, real-world infrastructure.

\begin{figure}[!t]
  \centering
  \includegraphics[width=\textwidth]{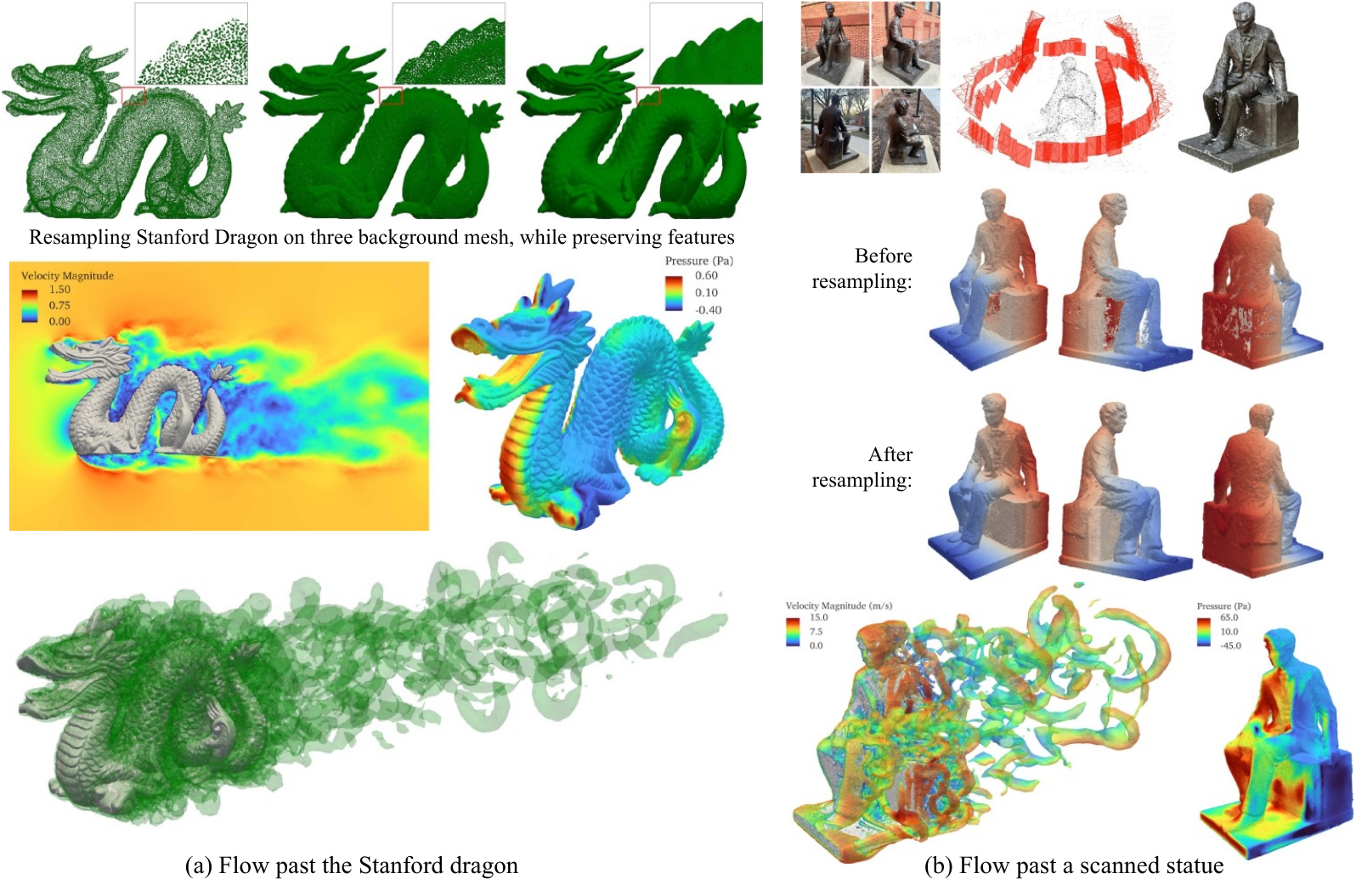}
  \caption{Mesh-driven resampling and regularization for robust point-cloud-based flow analysis on scanned objects. Adapted from Ref.~\cite{Jaiswal24Mesh}.}
  \label{fig:pointcloud-nimbus}
\end{figure}

Real-world point clouds are frequently characterized by noise, incompleteness, and irregular sampling due to sensor inaccuracies, surface reflections, and limited scanning coverage. To address these issues, Jaiswal et al.~\cite{Jaiswal24Mesh} introduced NIMBUS, an integrated framework combining immersogeometric CFD, photogrammetric reconstruction, and mesh-driven resampling. The resampling step uses the background fluid mesh to define the target point distribution. A local quadratic surface is fitted to the original points nearest to an intersected element. The medians of the element are then intersected with this surface to generate exactly one resampled point per intersected element. This guarantees that each intersected background element contains at least one point, which is the local requirement for effective weak enforcement of the no-slip condition. Additionally, a ghost penalty operator was applied to address both the stability and ill-conditioning issues caused by small cut elements inherent to immersed methods~\cite{Burman10Ghost, Burman12Ficti, Burman14Ficti, dePrenter17Condi, Badia18aggre, Larsson22finit, dePrenter23Stabi}. The approach was validated on a randomly sampled sphere to demonstrate hole-filling and noise robustness, and on the Stanford dragon to show the retention of fine geometric details while filling sparse regions. For the Stanford bunny, resampling the point cloud to match each background mesh enabled a direct mesh convergence study; without it, the simulation would have suffered from flow leakage through under-sampled cut elements. The methodology was benchmarked against boundary-fitted results for a scanned soda can: the resampled point cloud reduced the spurious pressure oscillations observed with a prior subsampled point cloud, bringing the results closer to the reference solutions. The complete pipeline was then applied to real-world scans containing severe geometric defects, including the Stanford dragon and an in-house photogrammetric scan of a statue with sparse, incomplete regions (Figure~\ref{fig:pointcloud-nimbus}). Ultimately, the framework eliminates the need for labor-intensive CAD reconstruction and enables stable flow simulations directly on imperfect scanned geometries.

\begin{figure}[!t]
  \centering
  \includegraphics[width=\textwidth]{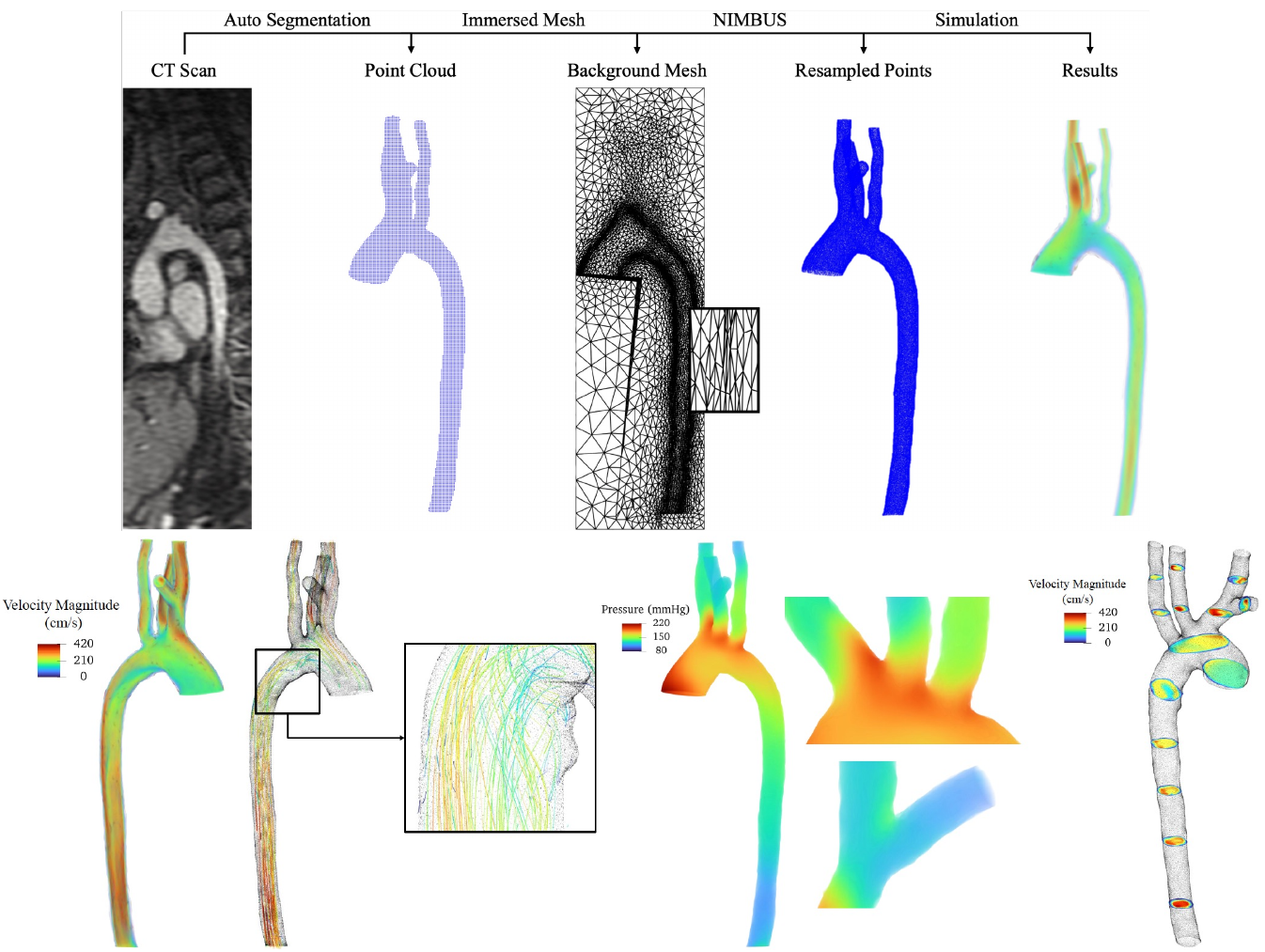}
  \caption{Direct medical image to flow analysis using auto-segmentation and point-cloud-based CFD. Adapted from Ref.~\cite{Corpuz25Direc}.}
  \label{fig:pointcloud-medical}
  \vspace{-4pt}
\end{figure}

Applied to point clouds obtained from auto-segmentation of medical images, the same framework yields patient-specific cardiovascular simulations directly from CT data~\cite{Corpuz25Direc}. The auto-segmentation uses a dual-stream 3D deep learning model that combines a U-Net with a gradient-based shape stream to capture both the main aortic vessel and its smaller branches. It produces a sparse dataset, with one point sampled at the center of each segmented voxel. NIMBUS then resamples these sparse points into an analysis-suitable point cloud that preserves the lumen geometry observed in the scan. To efficiently resolve near-wall hemodynamics, the framework uses an anisotropic background mesh built from a point-cloud-derived metric tensor that clusters elements in the wall-normal direction and relaxes them tangentially~\cite{Corpuz25Direc}. This delivers boundary-layer-like refinement around the immersed point cloud without the excessive element counts that isotropic refinement would produce. The framework was validated on benchmark hemodynamic problems with analytical or established reference solutions, including Hagen--Poiseuille pipe flow and steady flow through a bifurcation. It was then applied to a patient-specific aorta segmented from a CT scan, using a steady inflow corresponding to peak systole and resistance outlet conditions tuned to the peak-systolic pressures reported for each branch (Figure~\ref{fig:pointcloud-medical}). The entire pipeline, from imaging data to CFD results, proceeds with minimal manual intervention: no CAD reconstruction, no boundary-fitted mesh generation, and no manual segmentation. This makes it well suited for digital-twin and large-cohort cardiovascular workflows.

Across all of these developments, from CAD models to point clouds to medical images, the weak imposition of the no-slip condition is the underlying mechanism that makes flow analysis possible on geometries that are never meshed and, in several cases, never fully reconstructed.

\section{Weak BCs for compressible flow}
\label{sec:comp}

The majority of the weak BC literature focuses on incompressible flows. Extending these methods to the compressible regime introduces complexities due to the coupled mass, momentum, and energy fluxes at the wall. The first weak BC formulation for compressible flows adapted the incompressible weak BC velocity operator and incorporated the thermal condition as a standalone term~\cite{Xu17Compr, Xu19Immer1}. While this decoupled treatment performed well in many cases, instabilities emerged in certain flow scenarios. To maintain stability, these cases often relied on hybrid approaches requiring the strong enforcement of specific variables, such as strong temperature enforcement in the hypersonic regime~\cite{Codoni24Heat} or strong normal-velocity enforcement at adiabatic boundaries in transonic cases~\cite{Rajanna22Finit}. This limitation was later resolved by unifying the velocity, temperature, and heat-flux operators within the mathematical framework of advective--diffusive systems~\cite{Jaiswal26Weak}. The following sections review this unified formulation: the scalar model from Section~\ref{sec:scalar} is first generalized to a nonlinear flux, which is then applied to the compressible Navier--Stokes equations to systematically derive the weak wall BC operators.

\subsection{Nonlinear (general) flux}
Let $\Omega \subset \mathbb{R}^d$ ($d = 2$ or $3$) define a fixed spatial domain with a boundary $\Gamma = \partial \Omega$. Let $\mathbf{F}\left(u\right) : \Omega \rightarrow \mathbb{R}^d$ represent a nonlinear (general) flux and $S$ a source term, and let $\Gamma_{*} \subset \Gamma$ denote the portion of the boundary on which Dirichlet data $g : \Gamma_{*} \rightarrow \mathbb{R}$ is prescribed. The boundary value problem consists of finding $u : \overline{\Omega} \rightarrow \mathbb{R}$ such that
\begin{align}
\label{eq:general-flux-strong}
\nabla \cdot \mathbf{F}\left( u \right) - S &= 0 \quad \text{ in } \Omega \text{ ,} \\
u &= g \quad \text{ on } \Gamma_{*} \text{ .}
\end{align}
To derive the weak formulation, we multiply Eq.~\eqref{eq:general-flux-strong} by a test function $w$, integrate over the domain $\Omega$, and apply integration by parts to obtain
\begin{align}
\label{eq:general-flux-weak-1}
- \int_{\Omega} \nabla w \cdot \mathbf{F}\left( u \right) \,\mathrm{d}\Omega 
- \int_{\Omega} w S \,\mathrm{d}\Omega
+ \int_{\Gamma} w \mathbf{F}\left( u \right) \cdot \mathbf{n} \,\mathrm{d}\Gamma = 0 \text{ .} 
\end{align}
Following the same procedure as in Section~\ref{sec:advection} to enforce $u=g$ on $\Gamma_{*}$, we obtain the weak formulation: find $u$ such that for all $w$,
\begin{align}
\label{eq:general-flux-weak-2}
\int_{\Omega} w \left( \nabla \cdot \mathbf{F}\left( u \right) - S \right) \,\mathrm{d}\Omega
- \int_{\Gamma_{*}} w \left( \mathbf{F}\left( u \right) - \mathbf{F}\left( g \right)  \right) \cdot \mathbf{n} \,\mathrm{d}\Gamma = 0 \text{ .} 
\end{align}
For the boundary integral, we apply a Taylor expansion about the prescribed state $g$ and retain the linear term:
\begin{align}
\label{eq:taylor-expansion}
\mathbf{F}\left( u \right) - \mathbf{F}\left( g \right) \approx \left. \frac{\partial \mathbf{F}}{\partial u} \right|_{u=g} \left( u - g \right) \text{ .}
\end{align}
Hence, the final weak form is given by
\begin{align}
\label{eq:general-flux-weak-final}
\int_{\Omega} w \left( \nabla \cdot \mathbf{F}\left( u \right) - S \right) \,\mathrm{d}\Omega
- \int_{\Gamma_{*}} w \left( \left. \frac{\partial \mathbf{F}}{\partial u} \right|_{u=g} \cdot \mathbf{n} \right) \left( u - g \right) \,\mathrm{d}\Gamma = 0 \text{ .} 
\end{align}

\begin{remark}
In the special case $\mathbf{F}\left( u \right) = \mathbf{a}u$ and $\Gamma_{*}=\Gamma_{-}$, the formulation reduces to the weak form presented in Eq.~\eqref{eq:adv-weak-form-final} for the pure advection problem.
\end{remark}

\begin{remark}
The Jacobian $\partial \mathbf{F}/\partial u$ evaluated at $u=g$ plays the role of a solution-dependent convective velocity, generalizing the advection field $\mathbf{a}$ of Section~\ref{sec:advection} to the nonlinear setting.
\end{remark}

\subsection{The compressible Navier--Stokes equations and weak BC operators}
This section derives the formulation for weak enforcement of wall boundary conditions for the compressible Navier--Stokes equations using the design developed in the previous section. All equations are written in the Eulerian frame. Throughout the section, $(\cdot)_{,t}$ denotes a partial time derivative taken with respect to a fixed spatial coordinate and $(\cdot)_{,i}$ denotes a spatial gradient, where $i = 1, \ldots, d$ and $d=2,3$ is the space dimension. The Einstein summation convention for repeated indices is also used.

\subsubsection{Strong form}
\label{sec:cns-strong-form}
The reduced form of the Navier--Stokes equations of compressible flows from Refs.~\cite{Xu17Compr, Xu19Immer1} is taken as a starting point of the developments that follow. We introduce the vectors of conservation ($\mathbf{U}$) and pressure-primitive ($\mathbf{Y}$) variables as
\begin{align}
\label{conservation_var}
\mathbf{U}  = 
\begin{bmatrix}
 \rho \\
 \rho u_1 \\
 \rho u_2 \\
 \rho u_3 \\
 \rho e \\
\end{bmatrix}
\end{align}
and
\begin{align}
\label{primitive_var}
\mathbf{Y} = 
\begin{bmatrix}
 p \\
 u_1 \\
 u_2 \\
 u_3 \\
 T \\
\end{bmatrix}\text{ ,}
\end{align}
where $\mathbf{U}$ is in a reduced form that utilizes internal rather than total energy~\cite{Codoni21Stabi}. In the above, $\rho$ is the density, $u_i$ is the $i$th velocity component, $p$ is the pressure, and $T$ is the temperature. In addition, $e = c_\mathrm{v} T$ is the specific internal energy of the fluid, where $c_\mathrm{v} = R/(\gamma-1)$ is the specific heat at constant volume, $R$ is the ideal gas constant, and $\gamma$ is the heat capacity ratio. Finally, the pressure, density, and temperature unknowns are related through the ideal gas equation of state, $p = \rho R T$.

Using conservation variables, the Navier--Stokes equations of compressible flows hold in a fixed spatial domain $\Omega \subset \mathbb{R}^d$ and may be written as 
\begin{align}
\label{NS-conservation}
\mathbf{U}_{,t}+\mathbf{F}^{\mathrm{adv} \backslash p}_{i,i} + \mathbf{F}^p_{i,i} + \mathbf{F}^{\mathrm{sp}} - \mathbf{F}^{\mathrm{diff}}_{i,i} - \mathbf{S} = \mathbf{0}\text{ ,}
\end{align}
where the convective flux (without the pressure) is given by
\begin{align}
\label{eq:convective_flux}
\mathbf{F}^{\mathrm{adv} \backslash p}_{i} =
\begin{bmatrix}
\rho u_i\\
\rho u_i u_1\\
\rho u_i u_2\\
\rho u_i u_3 \\
\rho  u_i e  \\
\end{bmatrix}\text{ ,}
\end{align}
the pressure flux is given by
\begin{align}
\label{eq:pressure_flux}
\mathbf{F}^p_{i} =
\begin{bmatrix}
0\\
p\delta_{1i}\\
p\delta_{2i}\\
p\delta_{3i}\\
0\\
\end{bmatrix}\text{ ,}
\end{align}
the diffusive flux is given by
\begin{align}
\label{eq:diff_flux}
\mathbf{F}^\mathrm{diff}_i = 
\begin{bmatrix}
0\\
\tau_{1i}\\
\tau_{2i}\\
\tau_{3i}\\
-\phi^q_i\\
\end{bmatrix}\text{ ,}
\end{align}
the contribution of stress-power in the energy equation is given by
\begin{align}
\label{eq:stress_power}
\mathbf{F}^\mathrm{sp} &= 
\begin{bmatrix}
0 \\
0 \\
0 \\
0\\
p u_{i,i} - \tau_{ij}u_{j,i}\\
\end{bmatrix}\text{ ,}
\end{align}
and $\mathbf{S}$ is the source term. In Eqs.~\eqref{eq:pressure_flux}--\eqref{eq:stress_power}, $\delta_{ij}$ is the Kronecker delta, and $\tau_{ij}$ and $\phi^q_i$ are the viscous stress and heat flux, respectively, given by
$
\tau_{ij} = \lambda u_{k,k}\delta_{ij}+\mu\left(u_{i,j}+u_{j,i}\right)
$
and
$
 \phi^q_i = - \kappa T_{,i}
$, 
where $\mu$ is the dynamic viscosity, $\lambda = - 2 \mu /3$, and $\kappa$ is the thermal conductivity.

\subsubsection{Weak form}
\label{sec:cns-weak-form}
A weak form of the compressible-flow equations may be stated using conservation variables as follows: find $\mathbf{U}$ such that for all vector-valued test functions $\mathbf{W}$,
\begin{align}
\label{weak_form_U}
&\int_{\Omega}\mathbf{W}\cdot\mathbf{U}_{,t}\,\mathrm{d}\Omega + \int_{\Omega}\mathbf{W}\cdot\mathbf{F}^{\mathrm{adv} \backslash p}_{i,i}\,\mathrm{d}\Omega - \int_{\Omega}\mathbf{W}_{,i}\cdot \mathbf{F}^{p}_{i}\,\mathrm{d}\Omega 
 + \int_{\Omega}\mathbf{W}_{,i}\cdot \mathbf{F}^{\mathrm{diff}}_{i}\,\mathrm{d}\Omega 
 \nonumber \\
&\quad
 + \int_{\Omega}\mathbf{W} \cdot \mathbf{F}^{\mathrm{sp}}\,\mathrm{d}\Omega 
- \int_{\Omega}\mathbf{W} \cdot \mathbf{S}\,\mathrm{d}\Omega - \int_{\Gamma^H}\mathbf{W} \cdot \mathbf{H}\,\mathrm{d}\Gamma = 0 \text{ ,}
\end{align}
where $\Gamma^H$ is the part of the boundary of $\Omega$ with prescribed traction and heat-flux boundary conditions, whose values are collected in the vector $\mathbf{H}$. The vector $\mathbf{H}$ is given by
\begin{align}
\label{eq:cns_H}
\mathbf{H} &=
\begin{bmatrix}
    0 \\
    -p n_1 + \tau_{1i} n_i \\
    -p n_2 + \tau_{2i} n_i \\
    -p n_3 + \tau_{3i} n_i \\
    -\phi^q_i n_i
\end{bmatrix}\text{ ,}
\end{align}
where $n_i$ is the $i$th component of the outward unit surface normal $\mathbf{n}$. Note that the pressure and diffusive fluxes are integrated by parts in the above formulation to enable natural imposition of traction and heat-flux BCs. 

While popular, conservation variables are perhaps not the most convenient option for computation (see, e.g., the discussion in Ref.~\cite{Hauke98compa}). A more natural choice is presented by the pressure-primitive variables. To reformulate the weak form of the compressible-flow problem given by Eq.~\eqref{weak_form_U} using pressure-primitive variables, we first use the chain rule of differentiation to compute
\begin{align}
\mathbf{U}_{,t} = \mathbf{U}_{,\mathbf{Y}} \mathbf{Y}_{,t} = \mathbf{A}_0\mathbf{Y}_{,t}
\end{align}
and
\begin{align}
\mathbf{F}^{\mathrm{adv} \backslash p}_{i,i} = \mathbf{F}^{\mathrm{adv} \backslash p}_{i,\mathbf{Y}}\,\mathbf{Y}_{,i} = \mathbf{A}^{\mathrm{adv} \backslash p}_i\,\mathbf{Y}_{,i}\text{ .}
\end{align}
In addition, a direct computation gives
\begin{align}
\label{pflux}
\mathbf{F}^{p}_{i} = \mathbf{A}^{p}_i\,\mathbf{Y}\text{ ,}
\end{align}
and
\begin{align}
\mathbf{F}^{\mathrm{diff}}_{i,i} = (\mathbf{K}_{ij}\mathbf{Y}_{,j})_{,i} \text{ .}
\end{align}
Introducing the above relations into the weak form given by Eq.~\eqref{weak_form_U}, we arrive at the weak form for the pressure-primitive variables: find $\mathbf{Y} = [p,\mathbf{u}, T]^\mathrm{T}$ such that for all vector-valued test functions $\mathbf{W} = [\delta p, \delta \mathbf{u}, \delta T]^\mathrm{T}$,
\begin{align}
\label{weak_form_P}
&\int_{\Omega}\mathbf{W}\cdot\mathbf{A}_0\mathbf{Y}_{,t}\,\mathrm{d}\Omega
+ \int_{\Omega}\mathbf{W}\cdot\mathbf{A}^{\mathrm{adv} \backslash p}_i\,\mathbf{Y}_{,i}\,\mathrm{d}\Omega - \int_{\Omega}\mathbf{W}_{,i}\cdot \mathbf{A}^{p}_i\,\mathbf{Y}\,\mathrm{d}\Omega 
 + \int_{\Omega}\mathbf{W}_{,i}\cdot \mathbf{K}_{ij}\mathbf{Y}_{,j}\,\mathrm{d}\Omega 
\nonumber \\
& \quad 
+ \int_{\Omega}\mathbf{W} \cdot \mathbf{F}^{\mathrm{sp}}\,\mathrm{d}\Omega 
- \int_{\Omega}\mathbf{W} \cdot \mathbf{S}\,\mathrm{d}\Omega - \int_{\Gamma^H}\mathbf{W} \cdot \mathbf{H}\,\mathrm{d}\Gamma = 0\text{ .}  
\end{align}
Detailed expressions for the matrices $\mathbf{A}_0$, $\mathbf{A}^{\mathrm{adv} \backslash p}_i$, $\mathbf{A}^{p}_i$, and $\mathbf{K}_{ij}$ appearing in Eq.~\eqref{weak_form_P} may be found in Ref.~\cite[Appendix A]{Xu17Compr}.

\subsubsection{Design of weak wall BC operators}
\label{sec:weak_BC}

Let $\mathbf{G}$ be a vector of prescribed wall BC values and $\Gamma^G$ be the part of the boundary of $\Omega$ where $\mathbf{G}$ is prescribed. To design a weak wall BC operator on $\Gamma^G$, we focus on the convective and pressure flux terms. Note that the convective-flux terms are not integrated by parts, whereas the pressure-flux terms are. The following terms are added to Eq.~\eqref{weak_form_P} or Eq.~\eqref{weak_form_U}:
\begin{align}
\label{weak_BC_adv_pres}
- &\int_{\Gamma^G} \mathbf{W} \cdot \left( \mathbf{F}^{\mathrm{adv} \backslash p}_{i}(\mathbf{Y}) - \mathbf{F}^{\mathrm{adv} \backslash p}_{i}(\mathbf{G})\right)\,n_i\,\mathrm{d}\Gamma \nonumber \\
+ &\int_{\Gamma^G} \mathbf{W} \cdot \mathbf{F}^{p}_{i}(\mathbf{G})\,n_i\,\mathrm{d}\Gamma\text{ .}
\end{align}
For pressure-primitive variables, $\mathbf{G}$ collects the wall values of $[p,\mathbf{u},T]^\mathrm{T}$ (specific cases are considered in the following sections). Expanding the difference in the first term of Eq.~\eqref{weak_BC_adv_pres} in a Taylor series and retaining only the linear terms yields
\begin{align}
\label{Taylor}
\mathbf{F}^{\mathrm{adv} \backslash p}_{i}(\mathbf{Y}) & \approx \mathbf{F}^{\mathrm{adv} \backslash p}_{i}(\mathbf{G}) + \mathbf{F}^{\mathrm{adv} \backslash p}_{i,\mathbf{Y}} |_{\mathbf{G}}\,(\mathbf{Y}-\mathbf{G}) \nonumber \\
& = \mathbf{F}^{\mathrm{adv} \backslash p}_{i}(\mathbf{G}) + \mathbf{A}^{\mathrm{adv} \backslash p}_i |_{\mathbf{G}}\,(\mathbf{Y}-\mathbf{G})\text{ .}
\end{align}
Substituting the above result together with Eq.~\eqref{pflux} into Eq.~\eqref{weak_BC_adv_pres} gives the weak wall BC operator
\begin{align}
\label{weak_BC_adv_pres_lin}
- &\int_{\Gamma^G} \mathbf{W} \cdot [\mathbf{A}^{\mathrm{adv} \backslash p}_i |_{\mathbf{G}}\,n_i]\,(\mathbf{Y}-\mathbf{G})\,\mathrm{d}\Gamma \nonumber \\
+ &\int_{\Gamma^G} \mathbf{W} \cdot [\mathbf{A}^{p}_i\,n_i]\,\mathbf{G}\,\mathrm{d}\Gamma \text{ .}
\end{align}
We now examine two wall BC cases: a prescribed temperature and a prescribed heat flux. We refer to the constant-temperature condition as isothermal and the zero-heat-flux condition as adiabatic. In both cases, the flow velocity at the wall is zero, and the pressure and density are unknown and are computed from an equation of state (the ideal gas law is assumed in this work).

\subsubsection{Temperature wall BCs}
\label{sec:temp-bc}
For the temperature wall BC, the vector $\mathbf{G}$ is given as follows:
\begin{align}
\label{G_Temp_BC}
\mathbf{G} = 
\begin{bmatrix}
 p \\
 0 \\
 0 \\
 0 \\
 T_g \\
\end{bmatrix}\text{ ,}
\end{align}
where $T_g$ is the prescribed temperature, which may vary in space and time, the zeros in the middle reflect zero flow velocity, and $p$ is an unknown taken from the solution. With these values, the matrices $\mathbf{A}^{\mathrm{adv} \backslash p}_i |_{\mathbf{G}}\,n_i$ and $\mathbf{A}^{p}_i\,n_i$ become

\begin{align}
\label{Conv_Mat_Temp_BC}
\mathbf{A}^{\mathrm{adv} \backslash p}_i |_{\mathbf{G}}\,n_i = 
\begin{bmatrix}
 0 & \rho n_1 & \rho n_2 & \rho n_3 & 0 \\
 0 & 0 & 0 & 0 & 0 \\
 0 & 0 & 0 & 0 & 0 \\
 0 & 0 & 0 & 0 & 0 \\
 0 & \rho c_\mathrm{v} T_g n_1 & \rho c_\mathrm{v} T_g n_2 & \rho c_\mathrm{v} T_g n_3 & 0 \\
\end{bmatrix}
\end{align}
and
\begin{align}
\label{Pres_Mat_Temp_BC}
\mathbf{A}^{p}_i\,n_i = 
\begin{bmatrix}
 0 & 0 & 0 & 0 & 0 \\
 n_1 & 0 & 0 & 0 & 0 \\
 n_2 & 0 & 0 & 0 & 0 \\
 n_3 & 0 & 0 & 0 & 0 \\
 0 & 0 & 0 & 0 & 0 \\
\end{bmatrix}\text{ .}
\end{align}
The products $[\mathbf{A}^{\mathrm{adv} \backslash p}_i |_{\mathbf{G}}\,n_i]\,(\mathbf{Y}-\mathbf{G})$ and $[\mathbf{A}^{p}_i\,n_i]\,\mathbf{G}$ are obtained from a direct computation as
\begin{align}
\label{Prod_Conv_TempBC}
[\mathbf{A}^{\mathrm{adv} \backslash p}_i |_{\mathbf{G}}\,n_i]\,(\mathbf{Y}-\mathbf{G}) = 
\begin{bmatrix}
 \rho\,(u_in_i - 0)\\
 0 \\
 0 \\
 0 \\
 \rho c_\mathrm{v} T_g\,(u_in_i-0) \\
\end{bmatrix}
\end{align}
and
\begin{align}
\label{Prod_Pres_TempBC}
[\mathbf{A}^{p}_i\,n_i]\,\mathbf{G} = 
\begin{bmatrix}
 0 \\
 p n_1 \\
 p n_2 \\
 p n_3 \\
 0 \\
\end{bmatrix}\text{ .}
\end{align}
Introducing the above results into Eq.~\eqref{weak_BC_adv_pres_lin} gives an explicit expression for the weak wall BC operator
\begin{align}
\label{weak_BC_adv_pres_lin_exp_temp}
-  \int_{\Gamma^G} \delta p\,\rho\,(u_in_i-0)\,\mathrm{d}\Gamma 
+  \int_{\Gamma^G} \delta u_i n_i\,p\,\mathrm{d}\Gamma 
-  \int_{\Gamma^G} \delta T\,\rho c_\mathrm{v} T_g\,(u_i n_i - 0)\,\mathrm{d}\Gamma\text{ .}
\end{align}
We note that the first two integrals in Eq.~\eqref{weak_BC_adv_pres_lin_exp_temp} represent the ``adjoint'' consistency and consistency terms for the pressure, respectively. While initially introduced in Ref.~\cite{Xu17Compr} in an \textit{ad hoc} fashion, these terms are derived systematically in Ref.~\cite{Jaiswal26Weak} using the machinery of weak BCs for advective--diffusive systems. This mathematical framework also naturally yields the last integral, a modification that further enhances the robustness of the formulation across a wider range of flow conditions, as evidenced by the numerical results shown in later sections.

We now state the complete weak temperature wall BC operator, adding the convective stabilization contributions not covered by the Taylor expansion approach and the contributions from the viscous and heat-flux terms:
\begin{align}
\label{wbc_temp}
 - & \int_{\Gamma^G} \delta p\,\rho\,(u_in_i - 0)\,\mathrm{d}\Gamma 
+  \int_{\Gamma^G} \delta u_i n_i\,p\,\mathrm{d}\Gamma 
-  \int_{\Gamma^G} \delta u_i\,\rho\{u_jn_j\}_{-}\,(u_i - 0_i)\,\mathrm{d}\Gamma \nonumber \\
- & \int_{\Gamma^G} \delta u_i\,2\mu u_{(i,j)}n_j\,\mathrm{d}\Gamma 
-  \int_{\Gamma^G} \delta u_in_i\,\lambda u_{j,j}\,\mathrm{d}\Gamma \nonumber \\
- & \int_{\Gamma^G} 2\mu\,\delta u_{(i,j)} n_j\,(u_i - 0_i)\,\mathrm{d}\Gamma 
-  \int_{\Gamma^G} \lambda\,\delta u_{i,i}\,(u_jn_j - 0)\,\mathrm{d}\Gamma 
+  \int_{\Gamma^G} \delta u_{i} \frac{C_\mathrm{P}\,\mu}{h}\,(u_i - 0_i)\,\mathrm{d}\Gamma \nonumber \\
- & \int_{\Gamma^G} \delta T\,\rho c_\mathrm{v} T_g\,(u_i n_i - 0)\,\mathrm{d}\Gamma 
-  \int_{\Gamma^G} \delta T\,\rho c_\mathrm{v} \{u_jn_j\}_{-}\,(T-T_g)\,\mathrm{d}\Gamma \nonumber \\
- & \int_{\Gamma^G} \delta T\,\kappa\,T_{,i}n_i\,\mathrm{d}\Gamma 
-  \int_{\Gamma^G} \kappa\,\delta T_{,i}n_i\,(T-T_g)\,\mathrm{d}\Gamma 
+  \int_{\Gamma^G} \delta T\,\frac{C_\mathrm{P}\,\kappa}{h}\,(T-T_g)\,\mathrm{d}\Gamma\text{ .}
\end{align}
Here, $C_\mathrm{P}$ is a positive constant taken just large enough to guarantee stability, $h$ is the local mesh size, and $(\cdot)_{(i,j)}$ indicates that the tensor $(\cdot)$ is symmetrized.

\newpage
\begin{remark}
A direct inspection of the convective matrix in Eq.~\eqref{Conv_Mat_Temp_BC} reveals that it has all zero eigenvalues. However, given the non-symmetry of the matrix, its eigenvalues alone are not sufficient to guarantee stability. The convective contribution from Eq.~\eqref{Conv_Mat_Temp_BC} to the boundary energy estimate is instead governed by the quadratic form $\left( \mathbf{Y} - \mathbf{G} \right) \cdot [ \mathbf{A}^{\mathrm{adv} \backslash p}_i |_{\mathbf{G}}\,n_i ]\,\left( \mathbf{Y} - \mathbf{G} \right)$, which, using Eq.~\eqref{Prod_Conv_TempBC}, evaluates to $\rho  c_\mathrm{v}  T_g  (T - T_g) (u_i n_i)$. This is proportional to the normal velocity $u_i n_i$, which vanishes at a no-slip wall for the exact solution. Therefore, the convective contribution carries no energy across $\Gamma^G$ and is stability-neutral. 
\end{remark}

\subsubsection{Heat-flux wall BCs}
\label{sec:flux-bc}
 
For the heat-flux wall BC, because the temperature is unknown, the vector $\mathbf{G}$ is now given by
\begin{align}
\label{G_Flux_BC}
\mathbf{G} = 
\begin{bmatrix}
 p \\
 0 \\
 0 \\
 0 \\
 T \\
\end{bmatrix}\text{ ,}
\end{align}
where the unknown temperature $T$ is substituted for the prescribed temperature $T_g$. This results in a slight modification of the product
\begin{align}
\label{Prod_Conv_FluxBC}
[\mathbf{A}^{\mathrm{adv} \backslash p}_i |_{\mathbf{G}}\,n_i]\,(\mathbf{Y}-\mathbf{G}) = 
\begin{bmatrix}
 \rho\,(u_in_i - 0)\\
 0 \\
 0 \\
 0 \\
 \rho c_\mathrm{v} T\,(u_in_i-0) \\
\end{bmatrix} \text{ ,}
\end{align}
while the pressure-flux contributions remain the same. Explicit expressions for the weak BC operator contributions from the convective and pressure-flux terms become
\begin{align}
\label{weak_BC_adv_pres_lin_exp_flux}
-  \int_{\Gamma^G} \delta p\,\rho\,(u_in_i-0)\,\mathrm{d}\Gamma 
+  \int_{\Gamma^G} \delta u_i n_i\,p\,\mathrm{d}\Gamma 
-  \int_{\Gamma^G} \delta T\,\rho c_\mathrm{v} T\,(u_i n_i - 0)\,\mathrm{d}\Gamma \text{ .}
\end{align}
The complete weak heat-flux wall BC operator, adding the convective stabilization contributions not covered by the Taylor expansion approach and the contributions from the viscous terms and the prescribed wall heat flux, becomes
\begin{align}
\label{wbc_flux}
- & \int_{\Gamma^G} \delta p\,\rho\,(u_in_i - 0)\,\mathrm{d}\Gamma 
+  \int_{\Gamma^G} \delta u_i n_i\,p\,\mathrm{d}\Gamma 
-  \int_{\Gamma^G} \delta u_i\,\rho\{u_jn_j\}_{-}\,(u_i - 0_i)\,\mathrm{d}\Gamma \nonumber \\
- & \int_{\Gamma^G} \delta u_i\,2\mu u_{(i,j)}n_j\,\mathrm{d}\Gamma 
-  \int_{\Gamma^G} \delta u_in_i\,\lambda u_{j,j}\,\mathrm{d}\Gamma \nonumber \\
- & \int_{\Gamma^G} 2\mu\,\delta u_{(i,j)} n_j\,(u_i - 0_i)\,\mathrm{d}\Gamma 
-  \int_{\Gamma^G} \lambda\,\delta u_{i,i}\,(u_jn_j - 0)\,\mathrm{d}\Gamma 
+  \int_{\Gamma^G} \delta u_{i} \frac{C_\mathrm{P}\,\mu}{h}\,(u_i - 0_i)\,\mathrm{d}\Gamma 
\nonumber \\
- & \int_{\Gamma^G} \delta T\,\rho c_\mathrm{v} T\,(u_i n_i - 0)\,\mathrm{d}\Gamma 
-  \int_{\Gamma^G} \delta T\,h_T\,\mathrm{d}\Gamma \text{ ,}
\end{align}
where $h_T$ is a known prescribed heat flux at the wall. 

\subsection{Aerodynamics on boundary-fitted meshes}
\label{sec:comp-bf}

\begin{figure}[!b]
  \centering
  \includegraphics[width=\textwidth]{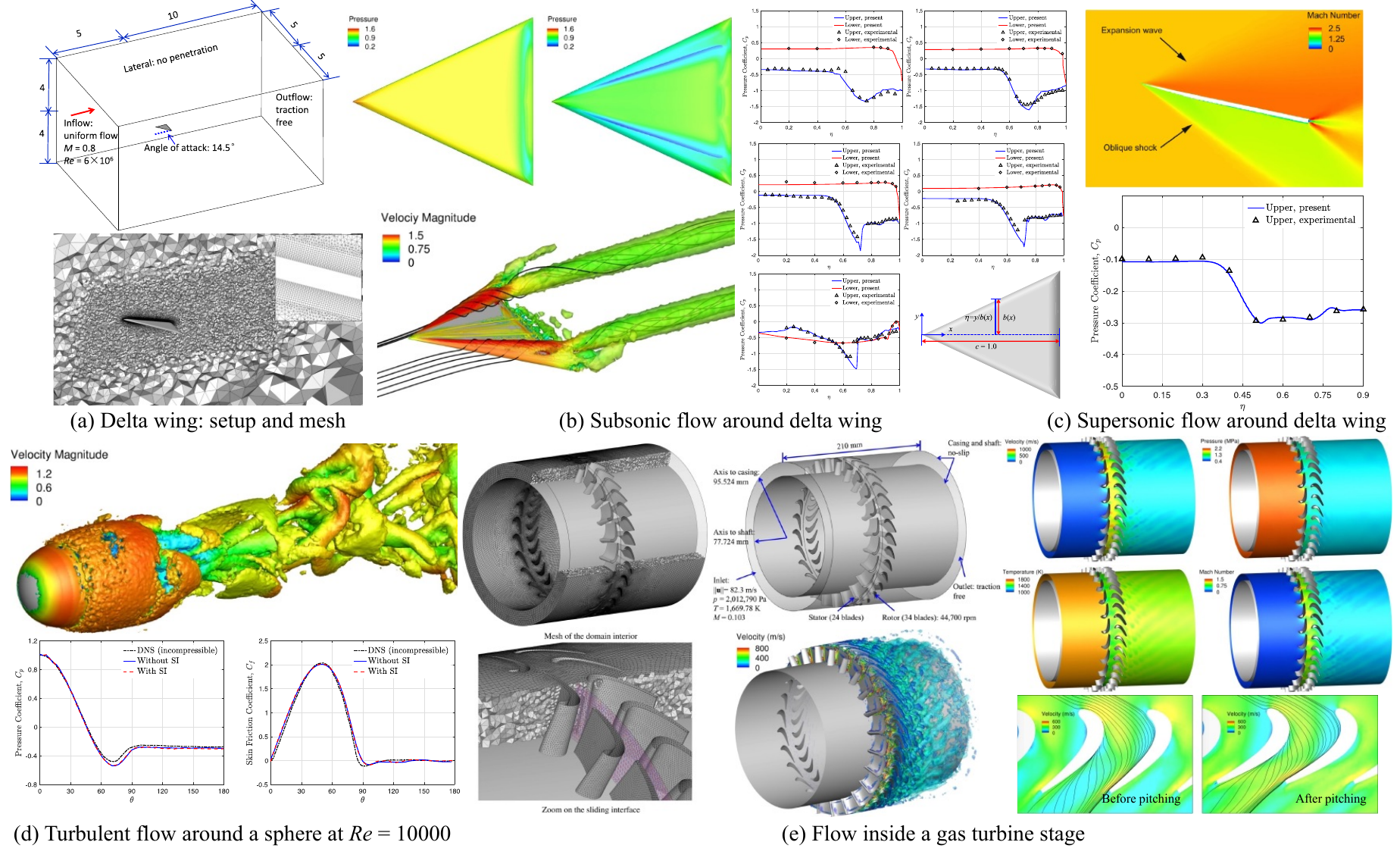}
  \caption{Weakly enforced wall conditions for compressible flow on boundary-fitted meshes: benchmark problems and gas-turbine aerodynamics. Adapted from Ref.~\cite{Xu17Compr}.}
  \label{fig:comp-bf-gas}
\end{figure}

The compressible weak BC formulation was first developed for moving-domain problems and applied, together with a sliding-interface formulation, to the aerodynamic analysis of a gas-turbine stage~\cite{Xu17Compr}. The preliminary validation studies in that work covered a wide range of Reynolds and Mach numbers, including subsonic ($M=0.8$) and supersonic ($M=2.0$) flows over a NASA delta wing. The computed surface pressure distributions agreed with the wind-tunnel data even on coarse boundary-layer meshes with a first-layer height of $y^+\approx 225$~\cite{Xu17Compr}. The sliding-interface formulation was validated on turbulent flow past a sphere at $Re=10{,}000$ and $M=0.1$. Computations with and without the sliding interface produced drag, pressure, and skin-friction quantities that were in excellent agreement. The validated method was then applied to the flow inside a gas-turbine stage. The results showed that the weak BCs deliver accurate near-wall behavior on coarse boundary-layer meshes in the compressible regime and remain consistent across non-matching, rotating subdomains (Figure~\ref{fig:comp-bf-gas}).

\begin{figure}[!t]
  \centering
  \includegraphics[width=\textwidth]{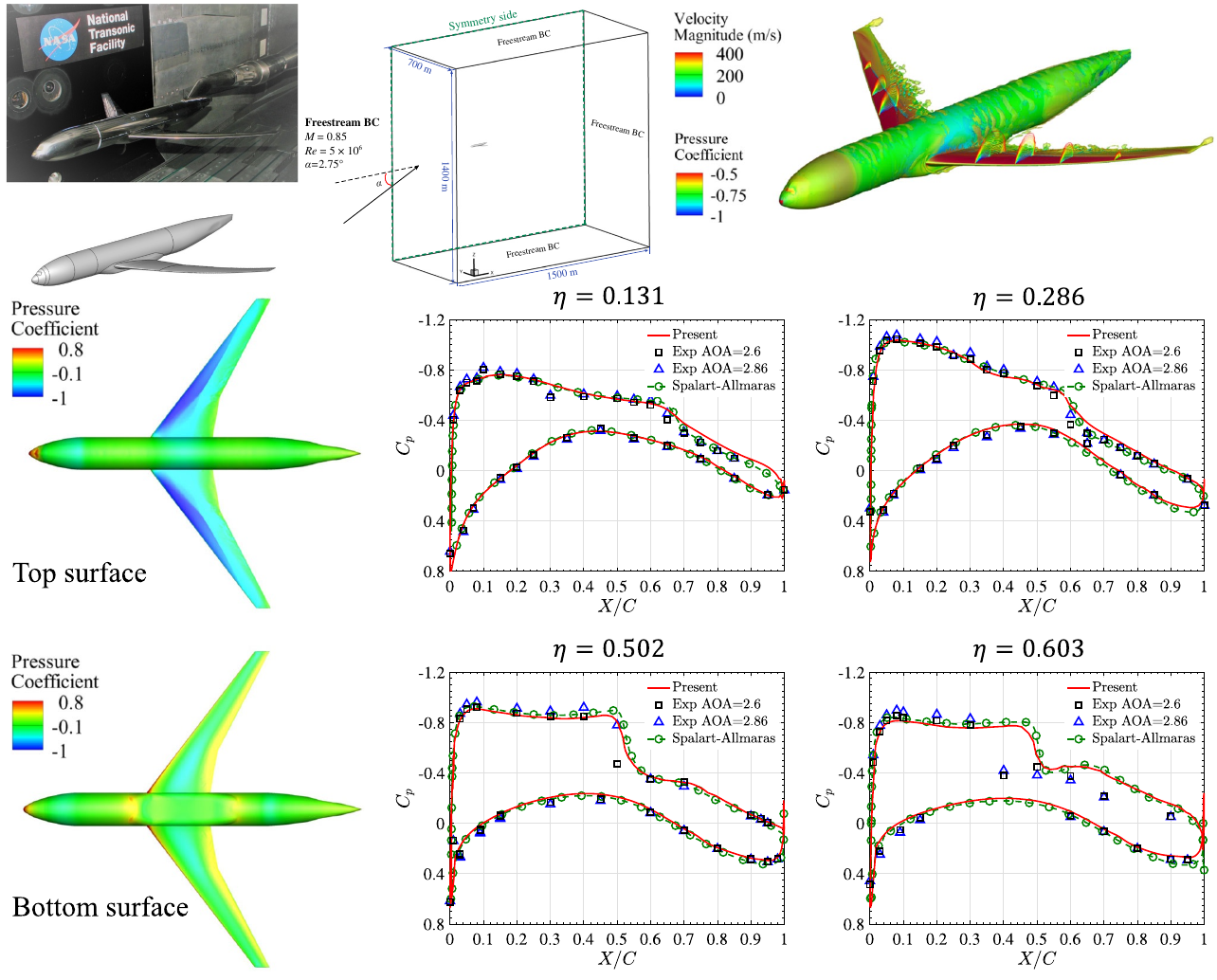}
  \caption{Weakly enforced wall conditions for compressible flow on boundary-fitted meshes: aircraft (NASA CRM) aerodynamics. Adapted from Ref.~\cite{Rajanna22Finit}.}
  \label{fig:comp-bf-air}
\end{figure}

The framework was later validated extensively for canonical aircraft aerodynamic configurations, including the NACA 0012 and RAE 2822 airfoils, the ONERA M6 wing, and the NASA Common Research Model (CRM) at various Mach and Reynolds numbers in Ref.~\cite{Rajanna22Finit}. An entropy-based discontinuity-capturing operator~\cite{Tezduyar86Disco, Hughes86Beyon, Hughes86disco, LeBeau91Finit, Almeida96Adapt, Hauke98compa, Tezduyar06Stabi, Hughes10Stabi} was combined with streamline upwind Petrov--Galerkin (SUPG) stabilization~\cite{Brooks82Strea, Hughes84Finit, Hughes86Symme, Hughes86gener, Hughes87Conve, Shakib91compr, LeBeau93SUPGf, Aliabadi93Space, Tezduyar94Massi, Mittal98unifi, Hauke01Simpl} and the weak BC operator to provide additional dissipation near shocks while preserving the consistency of the formulation. The computed pressure coefficients were all in good agreement with reference values (Figure~\ref{fig:comp-bf-air}). In all of these computations, the weak BCs behaved as a near-wall model, yielding accurate pressure and shock predictions without highly refined boundary-layer meshes. 

\begin{figure}[!t]
  \centering
  \includegraphics[width=\textwidth]{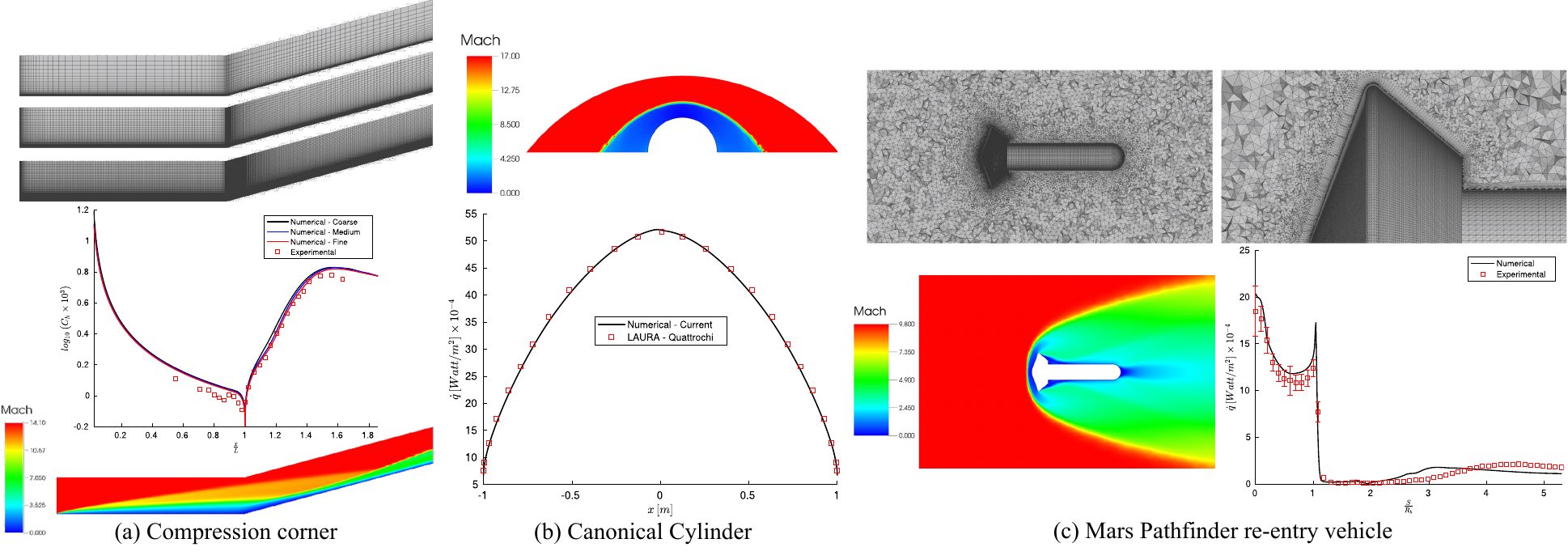}
  \caption{Weakly enforced wall conditions for hypersonic flow: surface heat-flux prediction for a compression corner, a canonical cylinder, and the Mars Pathfinder re-entry vehicle, using a stabilized formulation combining SUPG, discontinuity capturing, and weak BCs. Adapted from Ref.~\cite{Codoni24Heat}.}
  \label{fig:comp-hyper}
\end{figure}

The approach was later extended to the hypersonic regime by combining SUPG, an entropy-based discontinuity-capturing operator, and weak imposition of the no-slip condition for heat-flux prediction~\cite{Codoni24Heat}. The formulation was validated on a compression corner at $M=14.1$ with ramp angles of $15^\circ$ and $24^\circ$, the flow over a 2D cylinder at $M=17$, and the Mars Pathfinder re-entry vehicle at $M=10$. The predicted surface heat-flux distributions were in good agreement with the experimental measurements and reference numerical results across the benchmark set, including the separation region of the compression corner and the stagnation-point heating of the re-entry capsule (Figure~\ref{fig:comp-hyper}). The isothermal wall-temperature condition was imposed strongly in this work, which preserved stability across the high-Mach cases but motivated the unified treatment developed later~\cite{Jaiswal26Weak}, in which the velocity and temperature conditions were derived together within a single advective--diffusive framework.

\begin{figure}[!t]
  \centering
  \includegraphics[width=\textwidth]{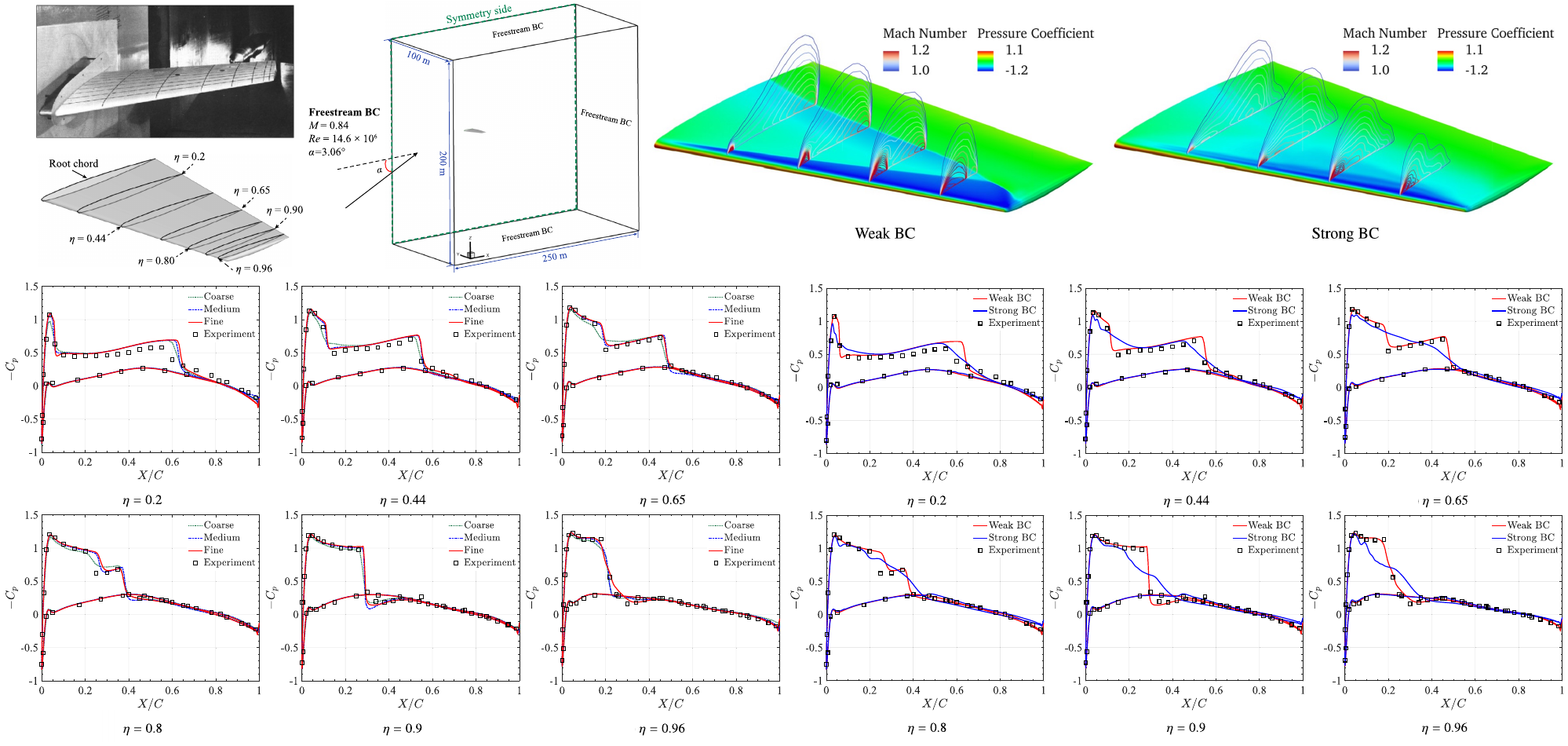}
  \caption{Weakly enforced compressible wall BCs: ONERA M6 wing pressure distribution compared with experiment, and weakly versus strongly enforced surface pressure and shock structure. Adapted from Ref.~\cite{Jaiswal26Weak}.}
  \label{fig:comp-new}
\end{figure}

The weak wall BC operators derived in Sections~\ref{sec:weak_BC}--\ref{sec:flux-bc} build on these earlier compressible formulations. By constructing the operator contributions from the advective and pressure fluxes within the framework of advective--diffusive systems, the derivation identifies the required weak BC terms systematically. This approach removes the need for the hybrid strong--weak enforcement used in earlier work for specific high-speed cases, such as the strong enforcement of the isothermal wall-temperature condition in the hypersonic regime~\cite{Codoni24Heat} and strong normal-velocity enforcement with adiabatic conditions in some transonic aircraft cases~\cite{Rajanna22Finit}. The stability and accuracy of the formulation were demonstrated on a set of 2D and 3D benchmark problems: subsonic flow over the NACA 0012 airfoil, transonic flow over the RAE 2822 airfoil, and transonic flow over the ONERA M6 wing~\cite{Jaiswal26Weak}. The predicted pressure distributions were in excellent agreement with reference experimental and numerical data. Accurate near-wall predictions were obtained without highly refined boundary-layer meshes, and on the same coarse meshes, the weak treatment captured the shock location and strength more accurately than strong enforcement (Figure~\ref{fig:comp-new}).

\subsection{Immersogeometric compressible flow}
\label{sec:comp-imga}

\begin{figure}[!b]
  \centering
  \includegraphics[width=\textwidth]{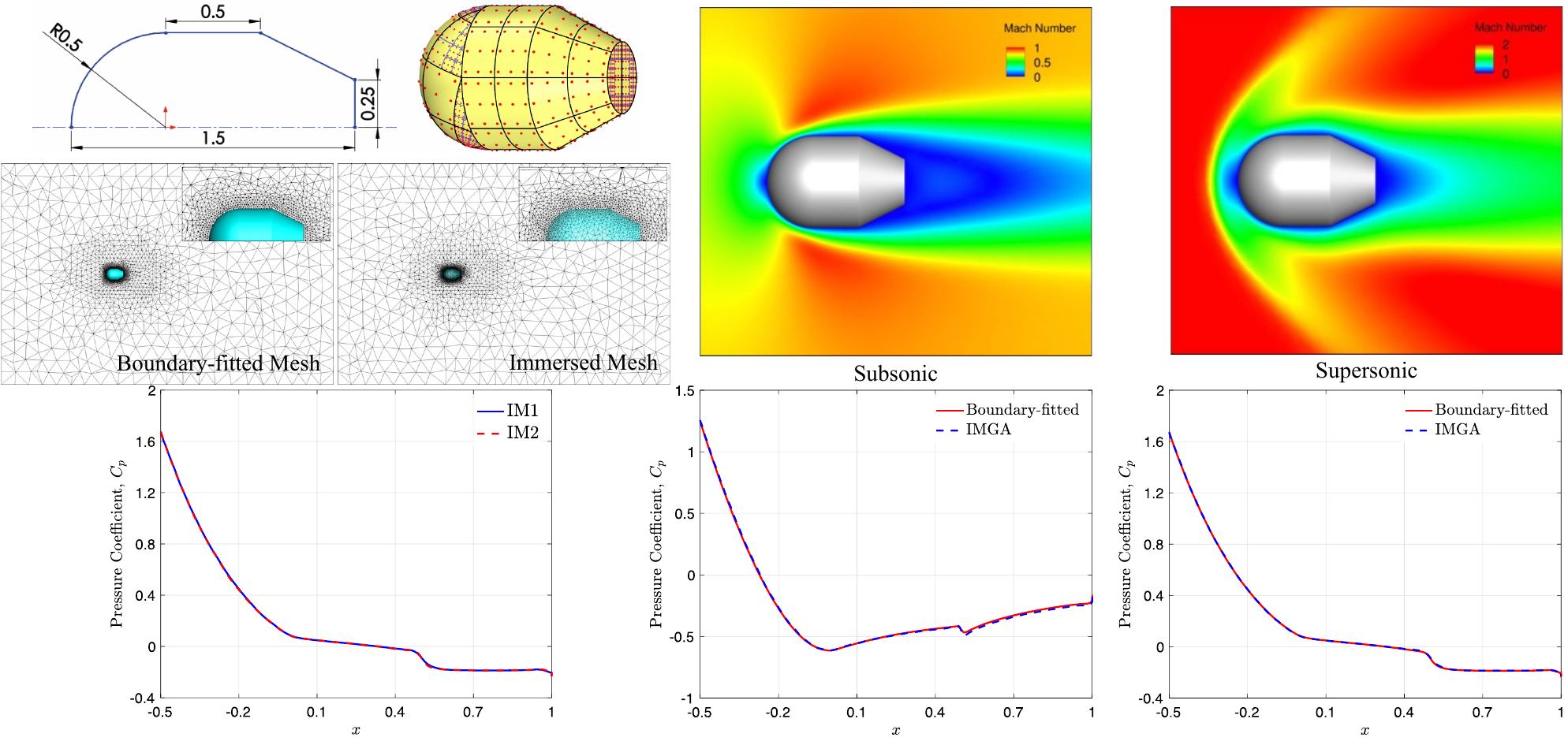}
  \caption{Immersogeometric analysis of compressible flow: subsonic and supersonic flows past a torpedo-shaped body. Adapted from Ref.~\cite{Xu19Immer1}.}
  \label{fig:comp-imga-torpedo}
\end{figure}

\begin{figure}[!t]
  \centering
  \includegraphics[width=\textwidth]{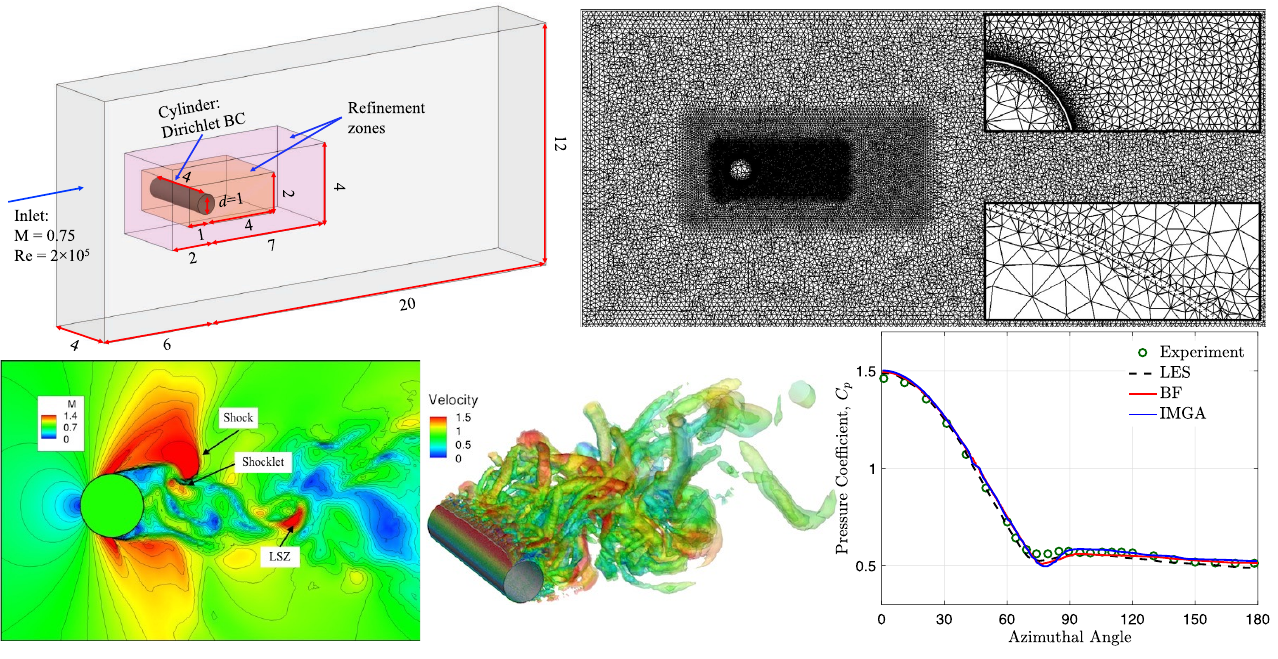}
  \caption{Immersogeometric analysis of compressible flow: transonic flow past a 3D cylinder. Adapted from Ref.~\cite{Xu19Immer1}.}
  \label{fig:comp-imga-cylinder}
\end{figure}

The compressible weak BC operator was combined with the immersogeometric framework to analyze compressible flow problems on non-boundary-fitted meshes~\cite{Xu19Immer1}. The operator was reformulated using the non-symmetric Nitsche method~\cite{Hartmann07Adjoi, Schillinger16nonsy}, which removes the formal lower bound on the penalty parameters required by the symmetric formulation and is especially convenient when these parameters are difficult to estimate in cut elements. The formulation can be used penalty-free~\cite{Burman12penal, Boiveau16penal, Dettmer16stabi}, or with a small penalty retained to improve the $L^2$ accuracy~\cite{Kirby05Selec, Heimann13unfit, Guo17param}, which was the choice made here. The formulation was validated on benchmark problems covering a wide range of Reynolds and Mach numbers, including laminar subsonic ($M=0.8$) and supersonic ($M=2.0$) flows past a 3D torpedo-shaped body and the turbulent transonic flow past a 3D cylinder at $Re=200{,}000$ and $M=0.75$~\cite{Xu19Immer1}. For the torpedo-shaped body, the IMGA framework accurately captured critical flow phenomena, such as the detached bow shock ahead of the leading edge in the supersonic regime. The predicted surface pressure distributions for both the subsonic and supersonic cases demonstrated excellent agreement with boundary-fitted reference computations using comparable mesh sizes (Figure~\ref{fig:comp-imga-torpedo}). The cylinder case features turbulent boundary-layer separation together with shocks, localized supersonic zones, and shocklets in the wake. The surface pressure distributions on the cylinder predicted by the IMGA framework agreed with the experimental data and the LES results from the literature, and with a boundary-fitted reference computed using the symmetric Nitsche method (Figure~\ref{fig:comp-imga-cylinder}).

\begin{figure}[!t]
  \centering
  \includegraphics[width=\textwidth]{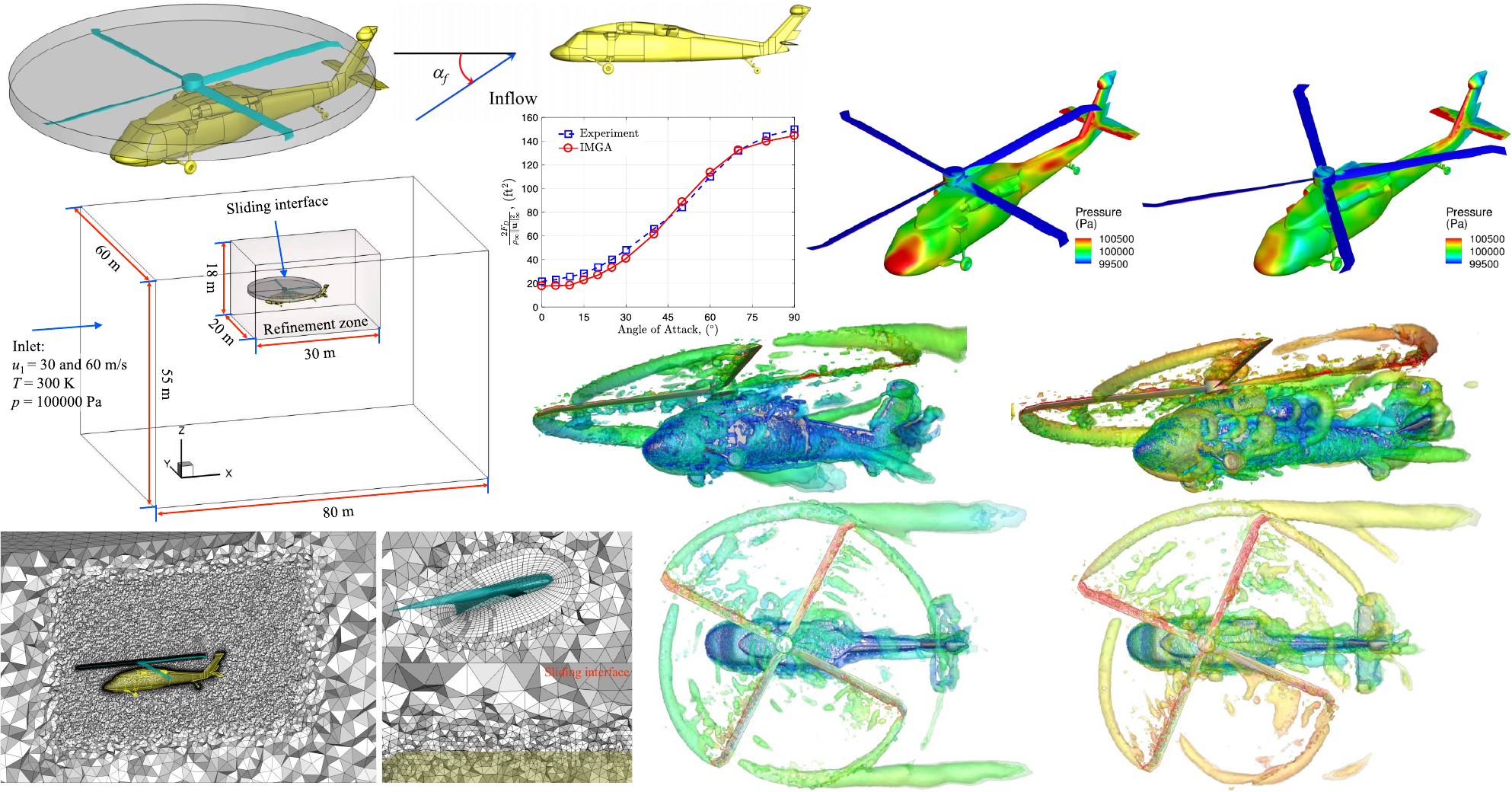}
  \caption{Immersogeometric analysis of compressible flow: rotorcraft aerodynamics in which the rotor is resolved with a boundary-fitted mesh inside a rotating sliding-interface subdomain, while the fuselage and landing gear are immersed in a non-boundary-fitted background mesh. Wall conditions are enforced weakly on both. Adapted from Ref.~\cite{Xu19Immer1}.}
  \label{fig:comp-imga}
  \vspace{-6pt}
\end{figure}

The framework was then applied to a full-scale UH-60A Black Hawk rotorcraft in forward flight~\cite{Xu19Immer1}. Two flight speeds were considered, 30~m/s and 60~m/s, with the latter corresponding to a free-stream Mach number of $0.17$ and a blade-tip relative Mach number of $0.63$. A hybrid meshing strategy was adopted. The main rotor blade was lofted through airfoil cross-sections to obtain a clean NURBS surface that is straightforward to mesh. A boundary-fitted prismatic boundary-layer mesh was generated around it inside a sliding-interface subdomain rotating at $27$~rad/s, with the no-slip condition on the blade imposed weakly using the symmetric Nitsche method. The fuselage and landing gear, given as a complex non-watertight tessellated surface, were instead immersed in a non-boundary-fitted background mesh and treated with the non-symmetric Nitsche method. The two subdomains were coupled across the rotating interface using the sliding-interface formulation extended to compressible flow~\cite{Xu19Immer1}. Following a mesh refinement study on fuselage-only computations at both flight speeds, the IMGA predictions for normalized fuselage drag versus angle of attack agreed with the NASA quarter-scale UH-60A wind-tunnel data~\cite{Howlett81UH60A}. Full-machine computations captured the tip-vortex propagation, the unsteady fuselage pressure distribution under the rotor downwash, and the wake asymmetry at the higher speed (Figure~\ref{fig:comp-imga}).

Building on these capabilities, the compressible IMGA framework has also been extended to operate directly on point-cloud data~\cite{Balu23Direc}. By evaluating the weak BC operator on the individual point samples, aerodynamic simulations of a torpedo-shaped body represented by a point cloud were successfully carried out in both subsonic and supersonic regimes. The predicted Mach number contours and surface pressure distributions showed excellent agreement with boundary-fitted reference solutions. This demonstrates the potential of performing high-speed compressible flow analysis directly on point-cloud objects without the need to reconstruct a CAD model or surface mesh.

\section{Conclusions}
\label{sec:conclusions}

This review has traced the development of weakly enforced boundary conditions as successive applications of a single unified idea. Beginning with the scalar advection and diffusion equations, we showed that the weak BC operator is assembled from a stabilizing inflow term and the consistency, adjoint-consistency, and penalty terms of Nitsche's method. These contributions extend to the incompressible Navier--Stokes equations, where the scalar operator is recast in vector form, with the normal and tangential wall conditions penalized separately. For compressible flow, the scalar operator is generalized to a nonlinear flux, which accommodates the coupling of mass, momentum, and energy at the wall and yields the velocity, temperature, and heat-flux conditions within a single framework. We also noted, in the context of both ALE and ST methods, the slip/sliding-interface techniques underpinned by weak BCs.

From this development, two distinct roles for weak BCs emerged. In boundary-fitted computations, weak enforcement acts as a near-wall model: it allows a controlled slip at the wall that imitates an unresolved thin boundary layer, recovering accurate mean-flow, wall heat-transfer, and engineering quantities on coarse meshes, as demonstrated across many examples and applications involving turbulent boundary layers. In immersed computations, the interface arbitrarily intersects the background elements. Strong imposition of no-slip and no-penetration wall conditions becomes intractable, making weak enforcement an essential alternative methodology. This observation underlies immersogeometric analysis and the direct geometry-to-analysis frameworks that have followed, in which flow analysis is carried out directly using CAD models, point clouds, photogrammetric reconstructions, and segmented medical images.

\section*{Acknowledgments}
M.-C. Hsu and M. Jaiswal would like to acknowledge funding from the National Institutes of Health under Award Number R01HL184128 and the National Science Foundation under Award Number DMS-2436623. Y. Bazilevs would like to acknowledge partial support from the Center for Information Geometric Mechanics and Optimization (CIGMO), a PSAAP-IV Focused Investigatory Center funded by the U.S.\ Department of Energy, National Nuclear Security Administration under Award Number DE-NA0004261. We also acknowledge the Texas Advanced Computing Center (TACC) at The University of Texas at Austin for providing computational resources that have contributed to the research results reported within this paper.

\small
\bibliographystyle{vancouver-mch}
\bibliography{refs-mch}

@article{Riviere01apri,
  author  = {Rivi\`{e}re, B. and Wheeler, M. F. and Girault, V.},
  title   = {A priori error estimates for finite element methods based on discontinuous approximation spaces for elliptic problems},
  journal = {SIAM Journal on Numerical Analysis},
  volume  = {39},
  pages   = {902--931},
  year    = {2001}
}

@article{Kirby05Selec,
  author  = {Kirby, R. M. and Karniadakis, G. E.},
  title   = {Selecting the numerical flux in discontinuous {G}alerkin methods for diffusion problems},
  journal = {Journal of Scientific Computing},
  volume  = {22},
  pages   = {385--411},
  year    = {2005}
}

@article{Hartmann07Adjoi,
  author  = {Hartmann, R.},
  title   = {Adjoint consistency analysis of discontinuous {G}alerkin discretizations},
  journal = {SIAM Journal on Numerical Analysis},
  volume  = {45},
  pages   = {2671--2696},
  year    = {2007}
}

@article{Heimann13unfit,
  author  = {Heimann, F. and Engwer, C. and Ippisch, O. and Bastian, P.},
  title   = {An unfitted interior penalty discontinuous {G}alerkin method for incompressible {N}avier--{S}tokes two-phase flow},
  journal = {International Journal for Numerical Methods in Fluids},
  volume  = {71},
  pages   = {269--293},
  year    = {2013}
}

@article{Boiveau16penal,
  author  = {Boiveau, T. and Burman, E.},
  title   = {A penalty-free {N}itsche method for the weak imposition of boundary conditions in compressible and incompressible elasticity},
  journal = {IMA Journal of Numerical Analysis},
  volume  = {36},
  pages   = {770--795},
  year    = {2016}
}

@article{Dettmer16stabi,
  author  = {Dettmer, W. G. and Kadapa, C. and Peri\'{c}, D.},
  title   = {A stabilised immersed boundary method on hierarchical b-spline grids},
  journal = {Computer Methods in Applied Mechanics and Engineering},
  volume  = {311},
  pages   = {415--437},
  year    = {2016}
}

@article{Guo17param,
  author  = {Guo, Y. and Ruess, M. and Schillinger, D.},
  title   = {A parameter-free variational coupling approach for trimmed isogeometric thin shells},
  journal = {Computational Mechanics},
  volume  = {59},
  pages   = {693--715},
  year    = {2017}
}

@article{Annavarapu12robus,
  author = {C. Annavarapu and M. Hautefeuille and J. E. Dolbow}, 
  title = {{A robust Nitsche's formulation for interface problems}}, 
  journal = {Computer Methods in Applied Mechanics and Engineering}, 
  volume = {225--228}, 
  pages = {44--54}, 
  year = {2012}
}

@article{Stenberg95Onsom,
  author  = {Stenberg, R.},
  title   = {On some techniques for approximating boundary conditions in the finite element method},
  journal = {Journal of Computational and Applied Mathematics},
  volume  = {63},
  pages   = {139--148},
  year    = {1995}
}

@article{Becker03finit,
  author  = {Becker, R. and Hansbo, P. and Stenberg, R.},
  title   = {A finite element method for domain decomposition with non-matching grids},
  journal = {ESAIM: Mathematical Modelling and Numerical Analysis},
  volume  = {37},
  pages   = {209--225},
  year    = {2003}
}

@article{Hansbo05Nitsc,
  author  = {Hansbo, P.},
  title   = {{N}itsche's method for interface problems in computational mechanics},
  journal = {{GAMM}-Mitteilungen},
  volume  = {28},
  pages   = {183--206},
  year    = {2005}
}

@article{Hansbo03Nitsc,
author = {P. Hansbo and J. Hermansson}, 
title = {{Nitsche's method for coupling non-matching meshes in fluid--structure vibration problems}}, 
journal = {Computational Mechanics}, 
volume = {32}, 
pages = {134--139}, 
year = {2003}
}

@article{Hansbo02unfit,
  author  = {Hansbo, A. and Hansbo, P.},
  title   = {An unfitted finite element method, based on {N}itsche's method, for elliptic interface problems},
  journal = {Computer Methods in Applied Mechanics and Engineering},
  volume  = {191},
  pages   = {5537--5552},
  year    = {2002}
}

@article{Baumann99disco,
  author  = {Baumann, C. E. and Oden, J. T.},
  title   = {A discontinuous {\emph{hp}} finite element method for convection--diffusion problems},
  journal = {Computer Methods in Applied Mechanics and Engineering},
  volume  = {175},
  pages   = {311--341},
  year    = {1999}
}

@article{Burman12penal,
  author  = {Burman, E.},
  title   = {A penalty-free nonsymmetric {N}itsche-type method for the weak imposition of boundary conditions},
  journal = {{SIAM} Journal on Numerical Analysis},
  volume  = {50},
  pages   = {1959--1981},
  year    = {2012}
}

@article{Schillinger16nonsy,
  author  = {Schillinger, D. and Harari, I. and Hsu, M.-C. and Kamensky, D. and Stoter, S. K. F. and Yu, Y. and Zhao, Y.},
  title   = {The non-symmetric {N}itsche method for the parameter-free imposition of weak boundary and coupling conditions in immersed finite elements},
  journal = {Computer Methods in Applied Mechanics and Engineering},
  volume  = {309},
  pages   = {625--652},
  year    = {2016}
}

@article{Benzaken24Const,
  author  = {Benzaken, J. and Evans, J. A. and Tamstorf, R.},
  title   = {Constructing {N}itsche's method for variational problems},
  journal = {Archives of Computational Methods in Engineering},
  volume  = {31},
  pages   = {1867--1896},
  year    = {2024}
}

@article{Warburton03Oncon,
  author  = {Warburton, T. and Hesthaven, J. S.},
  title   = {On the constants in {\emph{hp}}-finite element trace inverse inequalities},
  journal = {Computer Methods in Applied Mechanics and Engineering},
  volume  = {192},
  pages   = {2765--2773},
  year    = {2003}
}

@article{Main18shift2,
  author  = {Main, A. and Scovazzi, G.},
  title   = {The shifted boundary method for embedded domain computations. {P}art {II}: {L}inear advection--diffusion and incompressible {N}avier--{S}tokes equations},
  journal = {Journal of Computational Physics},
  volume  = {372},
  pages   = {996--1026},
  year    = {2018}
}

@article{Karki25Direc,
  author  = {Karki, S. and Shadkhah, M. and Yang, C.-H. and Balu, A. and Scovazzi, G. and Krishnamurthy, A. and Ganapathysubramanian, B.},
  title   = {Direct flow simulations with implicit neural representation of complex geometry},
  journal = {Computer Methods in Applied Mechanics and Engineering},
  volume  = {446},
  pages   = {118248},
  year    = {2025}
}

@article{Burman14Ficti,
author = {E. Burman and P. Hansbo}, 
title = {{Fictitious domain methods using cut elements: {III}. A stabilized Nitsche method for Stokes' problem}}, 
journal = {ESAIM: Mathematical Modelling and Numerical Analysis}, 
volume = {48}, 
pages = {859--874}, 
year = {2014}
}

@article{Piomelli02Walll,
  author  = {Piomelli, U. and Balaras, E.},
  title   = {Wall-layer models for large-eddy simulations},
  journal = {Annual Review of Fluid Mechanics},
  volume  = {34},
  pages   = {349--374},
  year    = {2002}
}

@article{Bose18Wallm,
  author  = {Bose, S. T. and Park, G. I.},
  title   = {Wall-modeled large-eddy simulation for complex turbulent flows},
  journal = {Annual Review of Fluid Mechanics},
  volume  = {50},
  pages   = {535--561},
  year    = {2018}
}

@article{Burman12Ficti,
  author  = {Burman, E. and Hansbo, P.},
  title   = {Fictitious domain finite element methods using cut elements: {II}. {A} stabilized {N}itsche method},
  journal = {Applied Numerical Mathematics},
  volume  = {62},
  pages   = {328--341},
  year    = {2012}
}

@article{Ruess13Weakl,
  author  = {Ruess, M. and Schillinger, D. and Bazilevs, Y. and Varduhn, V. and Rank, E.},
  title   = {Weakly enforced essential boundary conditions for {NURBS}-embedded and trimmed {NURBS} geometries on the basis of the finite cell method},
  journal = {International Journal for Numerical Methods in Engineering},
  volume  = {95},
  pages   = {811--846},
  year    = {2013}
}

@article{dePrenter17Condi,
  author  = {de Prenter, F. and Verhoosel, C. V. and van Zwieten, G. J. and van Brummelen, E. H.},
  title   = {Condition number analysis and preconditioning of the finite cell method},
  journal = {Computer Methods in Applied Mechanics and Engineering},
  volume  = {316},
  pages   = {297--327},
  year    = {2017}
}

@article{dePrenter18note,
  author  = {de Prenter, F. and Lehrenfeld, C. and Massing, A.},
  title   = {A note on the stability parameter in {N}itsche's method for unfitted boundary value problems},
  journal = {Computers \& Mathematics with Applications},
  volume  = {75},
  pages   = {4322--4336},
  year    = {2018}
}

@article{Badia18aggre,
  author  = {Badia, S. and Verdugo, F. and Mart\'{i}n, A. F.},
  title   = {The aggregated unfitted finite element method for elliptic problems},
  journal = {Computer Methods in Applied Mechanics and Engineering},
  volume  = {336},
  pages   = {533--553},
  year    = {2018}
}

@article{Larsson22finit,
  author  = {Larsson, K. and Kollmannsberger, S. and Rank, E. and Larson, M. G.},
  title   = {The finite cell method with least squares stabilized {N}itsche boundary conditions},
  journal = {Computer Methods in Applied Mechanics and Engineering},
  volume  = {393},
  pages   = {114792},
  year    = {2022}
}

@article{dePrenter23Stabi,
  author  = {de Prenter, F. and Verhoosel, C. V. and van Brummelen, E. H. and Larson, M. G. and Badia, S.},
  title   = {Stability and conditioning of immersed finite element methods: analysis and remedies},
  journal = {Archives of Computational Methods in Engineering},
  volume  = {30},
  pages   = {3617--3656},
  year    = {2023}
}

@article{Schillinger12isoge,
author = {D. Schillinger and L. Ded{\`e} and M. A. Scott and J. A. Evans and M. J. Borden and E. Rank and T. J. R. Hughes}, 
title = {{An isogeometric design-through-analysis methodology based on adaptive hierarchical refinement of {NURBS}, immersed boundary methods, and T-spline {CAD} surfaces}}, 
journal = {Computer Methods in Applied Mechanics and Engineering}, 
volume = {249--252}, 
pages = {116--150}, 
year = {2012}
}

@article{Wassermann19Integ,
  author  = {Wassermann, B. and Kollmannsberger, S. and Yin, S. and Kudela, L. and Rank, E.},
  title   = {Integrating {CAD} and numerical analysis: ``dirty geometry'' handling using the finite cell method},
  journal = {Computer Methods in Applied Mechanics and Engineering},
  volume  = {351},
  pages   = {808--835},
  year    = {2019}
}

@article{Kudela20Direc,
  author  = {Kudela, L. and Kollmannsberger, S. and Almac, U. and Rank, E.},
  title   = {Direct structural analysis of domains defined by point clouds},
  journal = {Computer Methods in Applied Mechanics and Engineering},
  volume  = {358},
  pages   = {112581},
  year    = {2020}
}

@article{Griffith20Immer,
  author  = {Griffith, B. E. and Patankar, N. A.},
  title   = {Immersed methods for fluid--structure interaction},
  journal = {Annual Review of Fluid Mechanics},
  volume  = {52},
  pages   = {421--448},
  year    = {2020}
}

@article{Hartmann22Enfor,
  author  = {Hartmann, F. and Kollmannsberger, S.},
  title   = {Enforcing essential boundary conditions on domains defined by point clouds},
  journal = {Computers \& Mathematics with Applications},
  volume  = {113},
  pages   = {13--23},
  year    = {2022}
}

@article{Hughes84Finit,
  author  = {Hughes, T. J. R. and Tezduyar, T. E.},
  title   = {Finite element methods for first-order hyperbolic systems with particular emphasis on the compressible {E}uler equations},
  journal = {Computer Methods in Applied Mechanics and Engineering},
  volume  = {45},
  pages   = {217--284},
  year    = {1984}
}

@article{Hughes86Symme,
  author  = {Hughes, T. J. R. and Franca, L. P. and Mallet, M.},
  title   = {A new finite element formulation for computational fluid dynamics: {I}. {S}ymmetric forms of the compressible {E}uler and {N}avier--{S}tokes equations and the second law of thermodynamics},
  journal = {Computer Methods in Applied Mechanics and Engineering},
  volume  = {54},
  pages   = {223--234},
  year    = {1986}
}

@article{Hughes86Beyon,
  author  = {Hughes, T. J. R. and Mallet, M. and Mizukami, A.},
  title   = {A new finite element formulation for computational fluid dynamics: {II}. {B}eyond {SUPG}},
  journal = {Computer Methods in Applied Mechanics and Engineering},
  volume  = {54},
  pages   = {341--355},
  year    = {1986}
}

@article{Hughes86gener,
  author  = {Hughes, T. J. R. and Mallet, M.},
  title   = {A new finite element formulation for computational fluid dynamics: {III}. {T}he generalized streamline operator for multidimensional advective--diffusive systems},
  journal = {Computer Methods in Applied Mechanics and Engineering},
  volume  = {58},
  pages   = {305--328},
  year    = {1986}
}

@article{Hughes86disco,
  author  = {Hughes, T. J. R. and Mallet, M.},
  title   = {A new finite element formulation for computational fluid dynamics: {IV}. {A} discontinuity-capturing operator for multidimensional advective--diffusive systems},
  journal = {Computer Methods in Applied Mechanics and Engineering},
  volume  = {58},
  pages   = {329--336},
  year    = {1986}
}

@article{Hughes87Conve,
  author  = {Hughes, T. J. R. and Franca, L. P. and Mallet, M.},
  title   = {A new finite element formulation for computational fluid dynamics: {VI}. {C}onvergence analysis of the generalized {SUPG} formulation for linear time-dependent multidimensional advective--diffusive systems},
  journal = {Computer Methods in Applied Mechanics and Engineering},
  volume  = {63},
  pages   = {97--112},
  year    = {1987}
}

@article{Shakib91compr,
  author  = {Shakib, F. and Hughes, T. J. R. and Johan, Z.},
  title   = {A new finite element formulation for computational fluid dynamics: {X}. {T}he compressible {E}uler and {N}avier--{S}tokes equations},
  journal = {Computer Methods in Applied Mechanics and Engineering},
  volume  = {89},
  pages   = {141--219},
  year    = {1991}
}

@article{LeBeau93SUPGf,
  author  = {Le Beau, G. J. and Ray, S. E. and Aliabadi, S. K. and Tezduyar, T. E.},
  title   = {{SUPG} finite element computation of compressible flows with the entropy and conservation variables formulations},
  journal = {Computer Methods in Applied Mechanics and Engineering},
  volume  = {104},
  pages   = {397--422},
  year    = {1993}
}

@article{Tezduyar86Disco,
  author  = {Tezduyar, T. E. and Park, Y. J.},
  title   = {Discontinuity capturing finite element formulations for nonlinear convection--diffusion--reaction equations},
  journal = {Computer Methods in Applied Mechanics and Engineering},
  volume  = {59},
  pages   = {307--325},
  year    = {1986}
}

@inproceedings{LeBeau91Finit,
  author    = {Le Beau, G. J. and Tezduyar, T. E.},
  title     = {Finite element computation of compressible flows with the {SUPG} formulation},
  booktitle = {Advances in Finite Element Analysis in Fluid Dynamics},
  series    = {{FED}-Vol.~123},
  publisher = {{ASME}},
  address   = {New York},
  pages     = {21--27},
  year      = {1991}
}

@article{Aliabadi93Space,
  author  = {Aliabadi, S. K. and Tezduyar, T. E.},
  title   = {Space--time finite element computation of compressible flows involving moving boundaries and interfaces},
  journal = {Computer Methods in Applied Mechanics and Engineering},
  volume  = {107},
  pages   = {209--223},
  year    = {1993}
}

@article{Tezduyar94Massi,
  author  = {Tezduyar, T. E. and Aliabadi, S. K. and Behr, M. and Mittal, S.},
  title   = {Massively parallel finite element simulation of compressible and incompressible flows},
  journal = {Computer Methods in Applied Mechanics and Engineering},
  volume  = {119},
  pages   = {157--177},
  year    = {1994}
}

@article{Almeida96Adapt,
  author  = {Almeida, R. C. and Gale{\~a}o, A. C.},
  title   = {An adaptive {P}etrov--{G}alerkin formulation for the compressible {E}uler and {N}avier--{S}tokes equations},
  journal = {Computer Methods in Applied Mechanics and Engineering},
  volume  = {129},
  pages   = {157--176},
  year    = {1996}
}

@article{Mittal98unifi,
  author  = {Mittal, S. and Tezduyar, T. E.},
  title   = {A unified finite element formulation for compressible and incompressible flows using augmented conservation variables},
  journal = {Computer Methods in Applied Mechanics and Engineering},
  volume  = {161},
  pages   = {229--243},
  year    = {1998}
}

@article{Hauke01Simpl,
  author  = {Hauke, G.},
  title   = {Simple stabilizing matrices for the computation of compressible flows in primitive variables},
  journal = {Computer Methods in Applied Mechanics and Engineering},
  volume  = {190},
  pages   = {6881--6893},
  year    = {2001}
}

@article{Tezduyar06Stabi,
  author  = {Tezduyar, T. E. and Senga, M.},
  title   = {Stabilization and shock-capturing parameters in {SUPG} formulation of compressible flows},
  journal = {Computer Methods in Applied Mechanics and Engineering},
  volume  = {195},
  pages   = {1621--1632},
  year    = {2006}
}

@article{Hughes10Stabi,
  author  = {Hughes, T. J. R. and Scovazzi, G. and Tezduyar, T. E.},
  title   = {Stabilized methods for compressible flows},
  journal = {Journal of Scientific Computing},
  volume  = {43},
  pages   = {343--368},
  year    = {2010}
}

@article{Burman10Ghost,
  author  = {Burman, E.},
  title   = {Ghost penalty},
  journal = {Comptes Rendus Mathematique},
  volume  = {348},
  pages   = {1217--1220},
  year    = {2010}
}

@article{Burman15CutFE,
  author  = {Burman, E. and Claus, S. and Hansbo, P. and Larson, M. G. and Massing, A.},
  title   = {{CutFEM}: {D}iscretizing geometry and partial differential equations},
  journal = {International Journal for Numerical Methods in Engineering},
  volume  = {104},
  pages   = {472--501},
  year    = {2015}
}

@article{Takizawa15Space,
  author  = {Takizawa, K. and Tezduyar, T. E. and Mochizuki, H. and Hattori, H. and Mei, S. and Pan, L. and Montel, K.},
  title   = {Space--time {VMS} method for flow computations with slip interfaces ({ST-SI})},
  journal = {Mathematical Models and Methods in Applied Sciences},
  volume  = {25},
  pages   = {2377--2406},
  year    = {2015}
}

@article{Takizawa16Compu,
  author  = {Takizawa, K. and Tezduyar, T. E. and Kuraishi, T. and Tabata, S. and Takagi, H.},
  title   = {Computational thermo-fluid analysis of a disk brake},
  journal = {Computational Mechanics},
  volume  = {57},
  pages   = {965--977},
  year    = {2016}
}

@article{Takizawa17Poros,
  author  = {Takizawa, K. and Tezduyar, T. E. and Kanai, T.},
  title   = {Porosity models and computational methods for compressible-flow aerodynamics of parachutes with geometric porosity},
  journal = {Mathematical Models and Methods in Applied Sciences},
  volume  = {27},
  pages   = {771--806},
  year    = {2017}
}

@article{Tezduyar25chron1,
  author  = {Tezduyar, T. E. and Takizawa, K.},
  title   = {A chronological catalog of methods and solutions in the {S}pace--{T}ime {C}omputational {F}low {A}nalysis: {I}. {F}inite element analysis},
  journal = {Computational Mechanics},
  volume  = {75},
  pages   = {793--831},
  year    = {2025}
}

@article{Tezduyar25chron2,
  author  = {Tezduyar, T. E. and Takizawa, K.},
  title   = {A chronological catalog of methods and solutions in the {S}pace--{T}ime {C}omputational {F}low {A}nalysis: {II}. {I}sogeometric analysis},
  journal = {Computational Mechanics},
  volume  = {75},
  pages   = {833--874},
  year    = {2025}
}

@article{Takizawa16Space,
  author  = {Takizawa, K. and Tezduyar, T. E. and Asada, S. and Kuraishi, T.},
  title   = {Space--time method for flow computations with slip interfaces and topology changes ({ST-SI-TC})},
  journal = {Computers \& Fluids},
  volume  = {141},
  pages   = {124--134},
  year    = {2016}
}

@article{Takizawa17Heart,
  author  = {Takizawa, K. and Tezduyar, T. E. and Terahara, T. and Sasaki, T.},
  title   = {Heart valve flow computation with the integrated {S}pace--{T}ime {VMS}, {S}lip {I}nterface, {T}opology {C}hange and {I}sogeometric {D}iscretization methods},
  journal = {Computers \& Fluids},
  volume  = {158},
  pages   = {176--188},
  year    = {2017}
}

@article{Kuraishi19Space,
  author  = {Kuraishi, T. and Takizawa, K. and Tezduyar, T. E.},
  title   = {Space--time isogeometric flow analysis with built-in {R}eynolds-equation limit},
  journal = {Mathematical Models and Methods in Applied Sciences},
  volume  = {29},
  pages   = {871--904},
  year    = {2019}
}

@article{Kuraishi19Tire,
  author  = {Kuraishi, T. and Takizawa, K. and Tezduyar, T. E.},
  title   = {Tire aerodynamics with actual tire geometry, road contact and tire deformation},
  journal = {Computational Mechanics},
  volume  = {63},
  pages   = {1165--1185},
  year    = {2019}
}

@article{Kuraishi25Space,
  author  = {Kuraishi, T. and Xu, Z. and Takizawa, K. and Tezduyar, T. E. and Kakegami, T.},
  title   = {Space--time isogeometric analysis of tire aerodynamics with complex tread pattern, road contact, and tire deformation},
  journal = {Computational Mechanics},
  volume  = {75},
  pages   = {575--591},
  year    = {2025}
}

@article{Terahara26Space,
  author = {Terahara, T. and Miura, H. and Takizawa, K. and Tezduyar, T. E.},
  title  = {Space--time isogeometric analysis of the aortic-valve to aorta flow with high-resolution boundary-layer representation},
  journal = {Computational Mechanics},
  volume  = {78},
  pages   = {999--1013},
  year    = {2026}
}

@article{Taniguchi26Diric,
  author = {Taniguchi, Y. and Takizawa, K. and Kashiwabara, T. and Tezduyar, T. E.},
  title  = {{D}irichlet-boundary stabilization in space--time computational analysis: {D}erivation in 1{D} elastodynamics},
  year   = {2026},
  journal = {Mathematical Models and Methods in Applied Sciences}, 
  note   =  {(2026) \url{https://doi.org/10.1142/S0218202526500491}}
}

@article{Bazilevs12Isoge,
  author  = {Bazilevs, Y. and Hsu, M.-C. and Scott, M. A.},
  title   = {Isogeometric fluid--structure interaction analysis with emphasis on non-matching discretizations, and with application to wind turbines},
  journal = {Computer Methods in Applied Mechanics and Engineering},
  volume  = {249--252},
  pages   = {28--41},
  year    = {2012}
}

@article{Takizawa23Varia,
  author  = {Takizawa, K. and Otoguro, Y. and Tezduyar, T. E.},
  title   = {Variational multiscale method stabilization parameter calculated from the strain-rate tensor},
  journal = {Mathematical Models and Methods in Applied Sciences},
  volume  = {33},
  pages   = {1661--1691},
  year    = {2023}
}

@article{Xu21octre,
  author  = {Xu, S. and Gao, B. and Lofquist, A. and Fernando, M. and Hsu, M.-C. and Sundar, H. and Ganapathysubramanian, B.},
  title   = {An octree-based immersogeometric approach for modeling inertial migration of particles in channels},
  journal = {Computers \& Fluids},
  volume  = {214},
  pages   = {104764},
  year    = {2021}
}

@article{Ruberg12Subdi,
  author  = {R{\"u}berg, T. and Cirak, F.},
  title   = {Subdivision-stabilised immersed b-spline finite elements for moving boundary flows},
  journal = {Computer Methods in Applied Mechanics and Engineering},
  volume  = {209--212},
  pages   = {266--283},
  year    = {2012}
}

@article{Schott14new,
  author  = {Schott, B. and Wall, W. A.},
  title   = {A new face-oriented stabilized {XFEM} approach for {2D} and {3D} incompressible {N}avier--{S}tokes equations},
  journal = {Computer Methods in Applied Mechanics and Engineering},
  volume  = {276},
  pages   = {233--265},
  year    = {2014}
}

@article{Kozak20Optim,
  author  = {Kozak, N. and Rajanna, M. R. and Wu, M. C. H. and Murugan, M. and Bravo, L. and Ghoshal, A. and Hsu, M.-C. and Bazilevs, Y.},
  title   = {Optimizing gas turbine performance using the surrogate management framework and high-fidelity flow modeling},
  journal = {Energies},
  volume  = {13},
  pages   = {4283},
  year    = {2020}
}

@article{Kozak20High,
  author  = {Kozak, N. and Xu, F. and Rajanna, M. R. and Bravo, L. and Murugan, M. and Ghoshal, A. and Bazilevs, Y. and Hsu, M.-C.},
  title   = {High-fidelity finite element modeling and analysis of adaptive gas turbine stator--rotor flow interaction at off-design conditions},
  journal = {Journal of Mechanics},
  volume  = {36},
  pages   = {595--606},
  year    = {2020}
}

@article{Peskin02immer,
  author  = {Peskin, C. S.},
  title   = {The immersed boundary method},
  journal = {Acta Numerica},
  volume  = {11},
  pages   = {479--517},
  year    = {2002}
}

@article{Hsu12Wind,
  author  = {Hsu, M.-C. and Akkerman, I. and Bazilevs, Y.},
  title   = {Wind turbine aerodynamics using {ALE--VMS}: {V}alidation and the role of weakly enforced boundary conditions},
  journal = {Computational Mechanics},
  volume  = {50},
  pages   = {499--511},
  year    = {2012}
}

@article{Hsu14Finit,
  author  = {Hsu, M.-C. and Akkerman, I. and Bazilevs, Y.},
  title   = {Finite element simulation of wind turbine aerodynamics: validation study using {NREL} {P}hase {VI} experiment},
  journal = {Wind Energy},
  volume  = {17},
  pages   = {461--481},
  year    = {2014}
}

@article{Jaiswal26Weak,
  author  = {Jaiswal, M. and Rajanna, M. R. and Islam, M. R. and Hsu, M.-C. and Bazilevs, Y.},
  title   = {Weak wall boundary conditions for compressible flows},
  journal = {Engineering with Computers},
  volume  = {42},
  pages   = {16},
  year    = {2026}
}

@article{Saurabh21Indus,
  author  = {Saurabh, K. and Gao, B. and Fernando, M. and Xu, S. and Khanwale, M. A. and Khara, B. and Hsu, M.-C. and Krishnamurthy, A. and Sundar, H. and Ganapathysubramanian, B.},
  title   = {Industrial scale {L}arge {E}ddy {S}imulations with adaptive octree meshes using immersogeometric analysis},
  journal = {Computers \& Mathematics with Applications},
  volume  = {97},
  pages   = {28--44},
  year    = {2021}
}

@article{Codoni21Stabi,
  author  = {Codoni, D. and Moutsanidis, G. and Hsu, M.-C. and Bazilevs, Y. and Johansen, C. and Korobenko, A.},
  title   = {Stabilized methods for high-speed compressible flows: toward hypersonic simulations},
  journal = {Computational Mechanics},
  volume  = {67},
  pages   = {785--809},
  year    = {2021}
}

@article{Codoni24Heat,
  author  = {Codoni, D. and Bayram, A. and Rajanna, M. R. and Johansen, C. and Hsu, M.-C. and Bazilevs, Y. and Korobenko, A.},
  title   = {Heat flux prediction for hypersonic flows using a stabilized formulation},
  journal = {Computational Mechanics},
  volume  = {73},
  pages   = {419--426},
  year    = {2024}
}

@article{Rajanna22Finit,
  author  = {Rajanna, M. R. and Johnson, E. L. and Codoni, D. and Korobenko, A. and Bazilevs, Y. and Liu, N. and Lua, J. and Phan, N. and Hsu, M.-C.},
  title   = {Finite element methodology for modeling aircraft aerodynamics: development, simulation, and validation},
  journal = {Computational Mechanics},
  volume  = {70},
  pages   = {549--563},
  year    = {2022}
}

@article{Xu19Immer2,
  author  = {Xu, S. and Xu, F. and Kommajosula, A. and Hsu, M.-C. and Ganapathysubramanian, B.},
  title   = {Immersogeometric analysis of moving objects in incompressible flows},
  journal = {Computers \& Fluids},
  volume  = {189},
  pages   = {24--33},
  year    = {2019}
}

@article{Corpuz25Direc,
  author  = {Corpuz, A. M. and Jaiswal, M. and Du, P. and Ramachandra, A. B. and Wang, J.-X. and Hsu, M.-C.},
  title   = {Direct medical image to simulation using auto-segmentation and point cloud-based {CFD}},
  journal = {Advances in Computational Science and Engineering},
  volume  = {3},
  pages   = {95--124},
  year    = {2025}
}

@article{Jaiswal24Mesh,
  author  = {Jaiswal, M. and Corpuz, A. M. and Hsu, M.-C.},
  title   = {Mesh-driven resampling and regularization for robust point cloud-based flow analysis directly on scanned objects},
  journal = {Computer Methods in Applied Mechanics and Engineering},
  volume  = {432},
  pages   = {117426},
  year    = {2024}
}

@article{Wang23Photo,
  author  = {Wang, X. and Jaiswal, M. and Corpuz, A. M. and Paudel, S. and Balu, A. and Krishnamurthy, A. and Yan, J. and Hsu, M.-C.},
  title   = {Photogrammetry-based computational fluid dynamics},
  journal = {Computer Methods in Applied Mechanics and Engineering},
  volume  = {417},
  pages   = {116311},
  year    = {2023}
}

@article{Balu23Direc,
  author  = {Balu, A. and Rajanna, M. R. and Khristy, J. and Xu, F. and Krishnamurthy, A. and Hsu, M.-C.},
  title   = {Direct immersogeometric fluid flow and heat transfer analysis of objects represented by point clouds},
  journal = {Computer Methods in Applied Mechanics and Engineering},
  volume  = {404},
  pages   = {115742},
  year    = {2023}
}

@article{Xu19resid,
  author  = {Xu, S. and Gao, B. and Hsu, M.-C. and Ganapathysubramanian, B.},
  title   = {A residual-based variational multiscale method with weak imposition of boundary conditions for buoyancy-driven flows},
  journal = {Computer Methods in Applied Mechanics and Engineering},
  volume  = {352},
  pages   = {345--368},
  year    = {2019}
}

@article{Zhu20immer,
  author  = {Zhu, Q. and Xu, F. and Xu, S. and Hsu, M.-C. and Yan, J.},
  title   = {An immersogeometric formulation for free-surface flows with application to marine engineering problems},
  journal = {Computer Methods in Applied Mechanics and Engineering},
  volume  = {361},
  pages   = {112748},
  year    = {2020}
}

@article{Golshan15Large,
  author  = {Golshan, R. and Tejada-Mart{\'\i}nez, A. E. and Juha, M. and Bazilevs, Y.},
  title   = {Large-eddy simulation with near-wall modeling using weakly enforced no-slip boundary conditions},
  journal = {Computers \& Fluids},
  volume  = {118},
  pages   = {172--181},
  year    = {2015}
}

@article{Kamensky15immer,
  author  = {Kamensky, D. and Hsu, M.-C. and Schillinger, D. and Evans, J. A. and Aggarwal, A. and Bazilevs, Y. and Sacks, M. S. and Hughes, T. J. R.},
  title   = {An immersogeometric variational framework for fluid--structure interaction: {A}pplication to bioprosthetic heart valves},
  journal = {Computer Methods in Applied Mechanics and Engineering},
  volume  = {284},
  pages   = {1005--1053},
  year    = {2015}
}

@article{Takizawa18Stabi,
  author  = {Takizawa, K. and Tezduyar, T. E. and Otoguro, Y.},
  title   = {Stabilization and discontinuity-capturing parameters for space--time flow computations with finite element and isogeometric discretizations},
  journal = {Computational Mechanics},
  volume  = {62},
  pages   = {1169--1186},
  year    = {2018}
}

@article{Xu17Compr,
  author  = {Xu, F. and Moutsanidis, G. and Kamensky, D. and Hsu, M.-C. and Murugan, M. and Ghoshal, A. and Bazilevs, Y.},
  title   = {Compressible flows on moving domains: {S}tabilized methods, weakly enforced essential boundary conditions, sliding interfaces, and application to gas-turbine modeling},
  journal = {Computers \& Fluids},
  volume  = {158},
  pages   = {201--220},
  year    = {2017}
}

@article{Xu19Immer1,
  author  = {Xu, F. and Bazilevs, Y. and Hsu, M.-C.},
  title   = {Immersogeometric analysis of compressible flows with application to aerodynamic simulation of rotorcraft},
  journal = {Mathematical Models and Methods in Applied Sciences},
  volume  = {29},
  pages   = {905--938},
  year    = {2019}
}

@article{Wang17Rapid,
  author  = {Wang, C. and Xu, F. and Hsu, M.-C. and Krishnamurthy, A.},
  title   = {Rapid {B}-rep model preprocessing for immersogeometric analysis using analytic surfaces},
  journal = {Computer Aided Geometric Design},
  volume  = {52--53},
  pages   = {190--204},
  year    = {2017}
}

@article{Hsu16Direc,
  author  = {Hsu, M.-C. and Wang, C. and Xu, F. and Herrema, A. J. and Krishnamurthy, A.},
  title   = {Direct immersogeometric fluid flow analysis using {B}-rep {CAD} models},
  journal = {Computer Aided Geometric Design},
  volume  = {43},
  pages   = {143--158},
  year    = {2016}
}

@article{Xu16tetra,
  author  = {Xu, F. and Schillinger, D. and Kamensky, D. and Varduhn, V. and Wang, C. and Hsu, M.-C.},
  title   = {The tetrahedral finite cell method for fluids: {I}mmersogeometric analysis of turbulent flow around complex geometries},
  journal = {Computers \& Fluids},
  volume  = {141},
  pages   = {135--154},
  year    = {2016}
}

@article{Mittal05Immer,
  author  = {Mittal, R. and Iaccarino, G.},
  title   = {Immersed boundary methods},
  journal = {Annual Review of Fluid Mechanics},
  volume  = {37},
  pages   = {239--261},
  year    = {2005}
}

@article{Evans13Expli,
  author  = {Evans, J. A. and Hughes, T. J. R.},
  title   = {Explicit trace inequalities for isogeometric analysis and parametric hexahedral finite elements},
  journal = {Numerische Mathematik},
  volume  = {123},
  pages   = {259--290},
  year    = {2013}
}

@article{Embar10Impos,
  author  = {Embar, A. and Dolbow, J. and Harari, I.},
  title   = {Imposing {D}irichlet boundary conditions with {N}itsche's method and spline-based finite elements},
  journal = {International Journal for Numerical Methods in Engineering},
  volume  = {83},
  pages   = {877--898},
  year    = {2010}
}

@article{Duster08finit,
  author  = {D{\"u}ster, A. and Parvizian, J. and Yang, Z. and Rank, E.},
  title   = {The finite cell method for three-dimensional problems of solid mechanics},
  journal = {Computer Methods in Applied Mechanics and Engineering},
  volume  = {197},
  pages   = {3768--3782},
  year    = {2008}
}

@article{Schillinger15finit,
  author  = {Schillinger, D. and Ruess, M.},
  title   = {The finite cell method: {A} review in the context of higher-order structural analysis of {CAD} and image-based geometric models},
  journal = {Archives of Computational Methods in Engineering},
  volume  = {22},
  pages   = {391--455},
  year    = {2015}
}

@article{Parvizian07Finit,
  author  = {Parvizian, J. and D{\"u}ster, A. and Rank, E.},
  title   = {Finite cell method},
  journal = {Computational Mechanics},
  volume  = {41},
  pages   = {121--133},
  year    = {2007}
}

@article{Bazilevs07Varia,
  author  = {Bazilevs, Y. and Calo, V. M. and Cottrell, J. A. and Hughes, T. J. R. and Reali, A. and Scovazzi, G.},
  title   = {Variational multiscale residual-based turbulence modeling for large eddy simulation of incompressible flows},
  journal = {Computer Methods in Applied Mechanics and Engineering},
  volume  = {197},
  pages   = {173--201},
  year    = {2007}
}

@article{Bazilevs07Weak1,
  author  = {Bazilevs, Y. and Hughes, T. J. R.},
  title   = {Weak imposition of {D}irichlet boundary conditions in fluid mechanics},
  journal = {Computers \& Fluids},
  volume  = {36},
  pages   = {12--26},
  year    = {2007}
}

@article{Bazilevs07Weak2,
  author  = {Bazilevs, Y. and Michler, C. and Calo, V. M. and Hughes, T. J. R.},
  title   = {Weak {D}irichlet boundary conditions for wall-bounded turbulent flows},
  journal = {Computer Methods in Applied Mechanics and Engineering},
  volume  = {196},
  pages   = {4853--4862},
  year    = {2007}
}

@article{Bazilevs08NURBS,
  author  = {Bazilevs, Y. and Hughes, T. J. R.},
  title   = {{NURBS}-based isogeometric analysis for the computation of flows about rotating components},
  journal = {Computational Mechanics},
  volume  = {43},
  pages   = {143--150},
  year    = {2008}
}

@article{Bazilevs10Isoge,
  author  = {Bazilevs, Y. and Michler, C. and Calo, V. M. and Hughes, T. J. R.},
  title   = {Isogeometric variational multiscale modeling of wall-bounded turbulent flows with weakly enforced boundary conditions on unstretched meshes},
  journal = {Computer Methods in Applied Mechanics and Engineering},
  volume  = {199},
  pages   = {780--790},
  year    = {2010}
}

@article{Brooks82Strea,
  author  = {Brooks, A. N. and Hughes, T. J. R.},
  title   = {Streamline upwind/{P}etrov--{G}alerkin formulations for convection dominated flows with particular emphasis on the incompressible {N}avier--{S}tokes equations},
  journal = {Computer Methods in Applied Mechanics and Engineering},
  volume  = {32},
  pages   = {199--259},
  year    = {1982}
}

@article{Hughes05Isoge,
  author  = {Hughes, T. J. R. and Cottrell, J. A. and Bazilevs, Y.},
  title   = {Isogeometric analysis: {CAD}, finite elements, {NURBS}, exact geometry and mesh refinement},
  journal = {Computer Methods in Applied Mechanics and Engineering},
  volume  = {194},
  pages   = {4135--4195},
  year    = {2005}
}

@article{Kim87Turbu,
  author  = {Kim, J. and Moin, P. and Moser, R.},
  title   = {Turbulence statistics in fully developed channel flow at low {R}eynolds number},
  journal = {Journal of Fluid Mechanics},
  volume  = {177},
  pages   = {133--166},
  year    = {1987}
}

@article{Takizawa11Multi,
  author  = {Takizawa, K. and Tezduyar, T. E.},
  title   = {Multiscale space--time fluid--structure interaction techniques},
  journal = {Computational Mechanics},
  volume  = {48},
  pages   = {247--267},
  year    = {2011}
}

@article{Takizawa14Space,
  author  = {Takizawa, K. and Tezduyar, T. E. and Buscher, A. and Asada, S.},
  title   = {Space--time interface-tracking with topology change ({ST-TC})},
  journal = {Computational Mechanics},
  volume  = {54},
  pages   = {955--971},
  year    = {2014}
}

@article{Tezduyar00Finit,
  author  = {Tezduyar, T. E. and Osawa, Y.},
  title   = {Finite element stabilization parameters computed from element matrices and vectors},
  journal = {Computer Methods in Applied Mechanics and Engineering},
  volume  = {190},
  pages   = {411--430},
  year    = {2000}
}

@article{Nitsche71Uber,
  author  = {Nitsche, J.},
  title   = {\"{U}ber ein {V}ariationsprinzip zur {L}\"{o}sung von {D}irichlet-{P}roblemen bei {V}erwendung von {T}eilr\"{a}umen, die keinen {R}andbedingungen unterworfen sind},
  journal = {Abhandlungen aus dem Mathematischen Seminar der Universit\"{a}t Hamburg},
  volume  = {36},
  pages   = {9--15},
  year    = {1971}
}

@article{Arnold02Unifi,
  author  = {Arnold, D. N. and Brezzi, F. and Cockburn, B. and Marini, L. D.},
  title   = {Unified analysis of discontinuous {G}alerkin methods for elliptic problems},
  journal = {SIAM Journal on Numerical Analysis},
  volume  = {39},
  pages   = {1749--1779},
  year    = {2002}
}

@article{Harari92What,
  author  = {Harari, I. and Hughes, T. J. R.},
  title   = {What are {$C$} and {$h$}?: {I}nequalities for the analysis and design of finite element methods},
  journal = {Computer Methods in Applied Mechanics and Engineering},
  volume  = {97},
  pages   = {157--192},
  year    = {1992}
}

@article{Hauke98compa,
  author  = {Hauke, G. and Hughes, T. J. R.},
  title   = {A comparative study of different sets of variables for solving compressible and incompressible flows},
  journal = {Computer Methods in Applied Mechanics and Engineering},
  volume  = {153},
  pages   = {1--44},
  year    = {1998}
}

@article{Hsu10Impro,
  author  = {Hsu, M.-C. and Bazilevs, Y. and Calo, V. M. and Tezduyar, T. E. and Hughes, T. J. R.},
  title   = {Improving stability of stabilized and multiscale formulations in flow simulations at small time steps},
  journal = {Computer Methods in Applied Mechanics and Engineering},
  volume  = {199},
  pages   = {828--840},
  year    = {2010}
}

@book{Hughes87Finit,
  author    = {Hughes, T. J. R.},
  title     = {The Finite Element Method. Linear Static and Dynamic Finite Element Analysis},
  publisher = {Prentice-Hall},
  address   = {Englewood Cliffs, New Jersey},
  year      = {1987}
}

@book{Ciarlet78Finit,
  author    = {Ciarlet, P. G.},
  title     = {The Finite Element Method for Elliptic Problems},
  publisher = {North-Holland},
  address   = {Amsterdam},
  year      = {1978}
}

@book{Wilcox06Turbu,
  author    = {Wilcox, D. C.},
  title     = {Turbulence Modeling for {CFD}},
  publisher = {{DCW} Industries, Inc.},
  address   = {La Ca\~{n}ada, California},
  edition   = {Third},
  year      = {2006}
}

@techreport{Howlett81UH60A,
  title={{UH-60A Black Hawk engineering simulation program. Volume 1: Mathematical model}},
  author={Howlett, J. J.},
  Institution = {NASA},
  Number = {NASA-CR-16630},
  Type = {NASA Technical Report},
  year={1981}
}

\end{document}